\documentclass[
    a4paper,
    11pt,
]{article}

\usepackage[T1]{fontenc}
\usepackage{amsmath,amssymb,amsthm}
\usepackage{mathtools,dsfont}
\usepackage{paralist,enumitem,bm,placeins,url}
\usepackage{color,graphicx,subfig} 
\usepackage[margin=30mm,top=40mm]{geometry}
\usepackage[font=small]{caption}
\usepackage{cite}
\usepackage{diagbox}
\usepackage{array}
\newcolumntype{C}[1]{>{\centering\arraybackslash}p{#1}}
\usepackage{tocbasic}
\DeclareTOCStyleEntry[
beforeskip=.2em plus 1pt,
]{tocline}{section}
\usepackage{todonotes}
\usepackage{hyperref} 
\makeatletter
\hypersetup{%
	colorlinks	=true,
	linkcolor	=blue,%
	citecolor	=red,%
	urlcolor	=blue
}

\renewcommand{\vec}[1]{\mathbf{#1}}

\newcommand{\bs}[1]{{\boldsymbol{#1}}}
\newcommand{\cF}{{c_{_F}}}
\newcommand{\cFsqu}{{c^2_{_F}}}
\newcommand{\VS}{{\scriptscriptstyle{VS}}}
\newcommand{\FS}{{\scriptscriptstyle{FS}}}
\newcommand{\rd}{{\rm d}}
\newcommand{\nn}{\nonumber}

\newcommand{\bds}[1]{\boldsymbol{#1}}
 
\newcommand{\emb}{\hookrightarrow} 

\newcommand{\Om}{\Omega}

\newcommand{\R}{\mathbb{R}}
\newcommand{\N}{\mathbb{N}}

\newcommand{\dS}{\, \mathrm{d}S}

\newcommand{\ds}{\, \mathrm{d}s}
\newcommand{\dx}{\, \mathrm{d}x}
\newcommand{\dt}{\, \mathrm{d}t}

\newcommand{\drho}{\, \mathrm{d}\rho}
\newcommand{\deta}{\, \mathrm{d}\eta}

\newcommand{\del}{\partial}
\newcommand{\n}{\vec{n}}
\newcommand{\nw}{\vec{n}_w}
\newcommand{\delt}{\partial_t}

\newcommand{\abs}[1]{\left| #1 \right|}
\newcommand{\norm}[1]{\left\| #1 \right\|}
\newcommand{\ang}[2]{ \left< #1 \,{,}\, #2  \right>}
\newcommand{\bigang}[2]{ \big< #1 \,{,}\, #2  \big>}
\newcommand{\scp}[2]{ \left( #1 \,{,}\, #2  \right)}
\newcommand{\bigscp}[2]{\big( #1 \,{,}\, #2 \big)}
\newcommand{\mean}[1]{\left< #1 \right>}
\newcommand{\eps}{{\varepsilon}}
\newcommand{\mR}{{\mathcal{R}}}

\newcommand{\Grad}{\nabla}

\newcommand{\dtau}{\, d\tau}

\newcommand{\intO}{\int_{\Omega}}

\newcommand{\intQ}{\iint_Q}
\newcommand{\intGw}{\int_{\Gamma_w}}
\newcommand{\intSw}{\iint_{\Sigma_w}}

\newcommand{\p}{\vec{p}}
\newcommand{\q}{\vec{q}}

\newcommand{\ov}{\overline}
\newcommand{\revised}[1]{{#1}}

\theoremstyle{plain}
\newtheorem{thm}{Theorem}[section]

\newtheorem{cor}[thm]{Corollary}

\newtheorem{defn}[thm]{Definition}

\newtheorem{remark}[thm]{Remark}
\theoremstyle{definition}

\numberwithin{equation}{section}

{{\upshape\bfseries AMS subject classifications. }\ignorespaces}{}

\begin{document}
\title{
A diffuse-interface approach for solid-state dewetting with
anisotropic surface energies}

\author{Harald Garcke\footnotemark[1] 
    \and Patrik Knopf\footnotemark[1] 
    \and Robert N\"urnberg\footnotemark[2]
    \and Quan Zhao\footnotemark[1]}

\renewcommand{\thefootnote}{\fnsymbol{footnote}}

\footnotetext[1]{Fakult\"at f\"ur Mathematik, Universit\"at Regensburg, 93053 Regensburg, Germany, \\ \tt(\href{mailto:harald.garcke@ur.de}{harald.garcke@ur.de},
	\href{mailto:patrik.knopf@ur.de}{patrik.knopf@ur.de}
	\href{mailto:quan.zhao@ur.de}{quan.zhao@ur.de}) }
	
\footnotetext[2]{Dipartimento di Mathematica, Universit\`a di Trento,
38123 Trento, Italy, \\
\tt(\href{mailto:robert.nurnberg@unitn.it}{robert.nurnberg@unitn.it}) }

\date{}
\maketitle

\begin{abstract}
\noindent
We present a diffuse-interface model for the solid-state dewetting problem with anisotropic surface energies in $\R^d$ for $d\in\{2,3\}$. The introduced model consists of the anisotropic Cahn--Hilliard equation,
with either a smooth or a double-obstacle potential, together
with a degenerate mobility function and appropriate boundary conditions on the 
wall. 
Upon regularizing the introduced diffuse-interface model, and with the help of
suitable asymptotic expansions, we recover as the sharp-interface limit
the anisotropic surface diffusion flow for the interface together with an
anisotropic Young's law and a zero-flux condition at the contact line
of the interface with a fixed external boundary. 
Furthermore, we show the existence of weak solutions for the regularized model, for both smooth and obstacle potential. Numerical results based on an 
appropriate finite element approximation are presented to demonstrate the 
excellent agreement between the proposed diffuse-interface model and its sharp-interface limit. 
\end{abstract}


\setlength\parskip{1ex}
\renewcommand{\thefootnote}{\arabic{footnote}}

\section{Introduction}
\setlength\parindent{24pt}

Deposited \revised{solid} thin films are unstable and could dewet to form isolated islands  on the substrate in order to minimize the total surface energy \cite{leroy2016control,thompson2012solid}. This phenomenon is known as solid-state dewetting (SSD), since the thin films remain in a solid state during the process. SSD has attracted a lot of attention recently, and is emerging as a promising route to produce  patterns of arrays of particles used in sensor technology, optical and magnetic devices, and catalyst formations, see e.g.\ \cite{armelao2006recent,Benkouider2015ordered,schmidt2009silicon, Bollani19, Salvalaglio20, Backofen2017convexity}.

The dominant mass transport mechanism in SSD is surface diffusion \cite{Srolovitz86}. This evolution law was first introduced by Mullins \cite{Mullins57} to describe the mass diffusion within interfaces in polycrystalline materials. 
For surface diffusion, the normal velocity of the interface is proportional to the surface Laplacian of the mean curvature. 
In the case of SSD the evolution of the interface that separates the thin film
from the surrounding vapor also involves the motion of the contact line, i.e., the region where the film/vapor interface meets the substrate. The equilibrium contact angle is given by Young's law which prescribes a force balance along the substrate.  Many efforts have been devoted to SSD problems in recent years. For example, a large body of experiments have revealed that the pattern formations could depend highly on the crystallographic alignments, the film sizes and shapes, as well as the substrate topology, see e.g.\ \cite{YT2011, AKR2012, thompson2012solid, Naffouti17, Bollani19}. In addition, mathematical studies based on different models have been considered in \cite{Srolovitz86, Burger05, DBDLE06, FonsecaFLM12, Jiang12, Naffouti17, BRTMP2022stress, Jiang18, Wang15, JZB2020sharp, Dziwnik17anisotropic}. 

In this work, we aim to study the SSD problem with anisotropic surface energies in the diffuse-interface framework. In the isotropic case, diffuse-interface models are based on the Ginzburg--Landau energy
\begin{equation}
    \mathcal{E}_{iso}(\varphi) = \int_{\Omega}\frac{\varepsilon}{2}|\nabla\varphi|^2 + \varepsilon^{-1}F(\varphi)\dx,
    \label{eq:isoGLE}
\end{equation}
where $\Omega\subset\R^d$ is a given domain with $d\in\{2,3\}$, $\varphi:\Omega\to \R$ is the order parameter, $\varepsilon>0$ is a small parameter proportional to the thickness of the interfacial layer, and $F(\varphi)$ is the free energy density. The following three choices for $F$ are mainly used in the literature:
\begin{itemize}
\item [(i)] the smooth double-well potential \cite{Taylor94linking} 
\begin{subequations}
\begin{equation}
    F(\varphi)=\frac{1}{2}(1-\varphi^2)^2,
    \label{eq:dwF}
\end{equation}
which has two global minimum points at $\varphi=\pm 1$ and a local maximum point at $\varphi=0$;
\item [(ii)] the logarithmic potential \cite{CH1958free}
\begin{equation}
    F(\varphi)=\frac{1}{2}\theta\, [(1+\varphi)\,\log(1+\varphi) +(1-\varphi)\,\log(1-\varphi)]+\frac{1}{2}(1-\varphi^2),\label{eq:logF}
\end{equation}
where $\theta>0$ is the absolute temperature. This potential has two minima $\varphi = \pm(1-\tilde{k}(\theta))$, where $\tilde{k}(\theta)$ is a small positive real number satisfying $\tilde{k}(\theta)\to 0$ as $\theta\to 0$, and its usage enforces $\varphi$ to attain values within $(-1,1)$;
\item [(iii)] the double-obstacle potential \cite{Blowey1991cahn}
\begin{equation}
    F(\varphi)=\left\{
    \begin{array}{ll}
    \frac{1}{2}(1-\varphi^2)\quad &\mbox{if}\quad|\varphi|\leq 1,\\[0.5em]
    \infty\quad &\mbox{otherwise}.
    \end{array}\right.
    \label{eq:doF}
\end{equation}
It can be characterized via the \textit{deep quench limit} of the logarithmic potential, i.e., the limit of \eqref{eq:logF} as $\theta\to 0$.
 \end{subequations}
 \end{itemize}
The (isotropic) Cahn--Hilliard equation can be interpreted as a weighted $H^{-1}$-gradient flow of the free energy \eqref{eq:isoGLE}. It reads as
\begin{align}
\partial_t\varphi = \nabla\cdot(m(\varphi)\,\nabla\mu),\qquad\mu = 
-\varepsilon \Delta\varphi + \varepsilon^{-1} F'(\varphi),
\label{eq:isoch}
\end{align}
where $m(\varphi)$ is a mobility function,
together with Neumann boundary conditions for $\mu$ and $\varphi$. The Cahn--Hilliard equation was first introduced to study the spinodal decomposition in binary alloys \cite{CH1958free, Cahn1961spinodal} and has since then been used to model many other phenomenon, e.g., \cite{AbelGarck12, Garcke00singular, Khain2008, Bertozzi06image}. We note that the double-obstacle potential is not differentiable at $\varphi=\pm 1$, and the definition of the generalized chemical potential in this case becomes
\begin{equation}
    \mu\in -\varepsilon\Delta\varphi +\varepsilon^{-1}\partial F(\varphi),
\end{equation}
where $\partial F(\varphi)$ is the Fr\'echet sub-differential of $F$ at $\varphi$ and $\Delta\varphi$ has to be understood in a weak sense, see \cite{Blowey1991cahn,BGN2013stable}. In the case of a constant mobility $m(\varphi)\equiv 1$, \eqref{eq:isoch} converges to the Mullins--Sekerka problem \cite{ Mullins1963} as $\varepsilon\to 0$ \cite{pego1989front, Alikakos94C}. 
In order to obtain the surface diffusion equation in the sharp-interface limit, a degenerate mobility needs to be chosen. For example, it was shown in \cite{Cahn1996cahn} by a formal asymptotic analysis that the surface diffusion flow is recovered by considering a slow time scale $\tau=O(\varepsilon^{-1}t)$ of \eqref{eq:isoch} with $m(\varphi)= (1-\varphi^2)_+$ and with the potential $F(\varphi)$ 
either chosen as in \eqref{eq:doF}, or as in \eqref{eq:logF}
with $\theta = O(\varepsilon^\xi)$, $\xi>0$. 
When using the smooth double-well potential \eqref{eq:dwF} the situation is
less clear. While the limiting motion of surface diffusion is obtained
with the choice $m(\varphi)= (1-\varphi^2)^2$ \cite{Voigt2016comment, Ratz2006surface, Jiang12, Dai2014coarsen}, using the less degenerate
mobility $m(\varphi)= (1-\varphi^2)_+$ may not lead to pure surface diffusion 
in the limit $\varepsilon\to0$, since an additional bulk diffusion term 
is conjectured to be present due to the non-zero flux contributions \cite{Dai2014coarsen, Lee15deg, lee16sharp}. However, in all these cases, no rigorous proof for the sharp-interface limit or the presence of non-zero flux contributions are available so far. 

\revised{A natural} generalization of the free energy \eqref{eq:isoGLE} to the case of anisotropic surface energies is given by
\begin{align}
\mathcal{E}_\gamma(\varphi) = \int_{\Omega}\frac{\varepsilon}{2}|\gamma(\nabla\varphi)|^2 + \varepsilon^{-1} F(\varphi)\dx=\int_{\Omega}\varepsilon A(\nabla\varphi) + \varepsilon^{-1} F(\varphi)\dx \revised{,}
\label{eq:aniGL}
\end{align}
\revised{see e.g.\ \cite{Kobayashi93,Elliott97}.}
Here, $\gamma: \R^d\to [0, \infty)$ is the anisotropic density function, which is positively homogeneous of degree one, and $A: = \frac{1}{2}\gamma^2$.  This then gives rise to the anisotropic Cahn--Hilliard equation
\begin{equation}
  \partial_t\varphi = \nabla\cdot(m(\varphi)\nabla\mu),\qquad  \mu = -\varepsilon\nabla\cdot A^\prime(\nabla\varphi) +\varepsilon^{-1} F^\prime(\varphi),
  \label{eq:asoch}
\end{equation}
where $A^\prime$ represents the gradient of the map $A:\R^d\to [0, \infty)$. In contrast to the isotropic case,  diffuse-interface models based on \eqref{eq:aniGL} result in a nonuniform asymptotic interface thickness, which in fact depends on the anisotropic density function $\gamma(\nabla\varphi)$, see \cite{Wheeler1996xi,wheeler2006phase, BP96anisotropic, Elliott96limit, Alfaro2010motion}. To remedy this issue, an alternative energy of the form 
\begin{equation}
  \widetilde{\mathcal{E}}_\gamma(\varphi) = \int_{\Omega}|\nabla\varphi|^{-1}\gamma(\nabla\varphi)\left(\frac{\varepsilon}{2}|\nabla\varphi|^2+ \varepsilon^{-1} F(\varphi)\right)\dx 
    \label{eq:altasoGL}
\end{equation}
can be considered, see \cite{Torabi209,Salva15face}, so that a constant thickness of the asymptotic interface is achieved. However,  the resulting diffuse-interface models based on \eqref{eq:altasoGL} become more nonlinear and are singular at $\nabla\varphi=0$, which poses great challenges in the mathematical analysis and the stable numerical approximation. Therefore, in this work, we will restrict ourselves to the classical energy in \eqref{eq:aniGL}. We also note that to guarantee that \eqref{eq:asoch} converges to the anisotropic surface diffusion flow as $\varepsilon\to 0$, a rescaled anisotropic coefficient needs to be introduced to the degenerate mobility \cite{Ratz2006surface, Li2009geometric}. We refer to Section~\ref{sec:mod} below for the precise details.

\begin{figure}[!htp]
    \centering
    \includegraphics[width=0.55\textwidth]{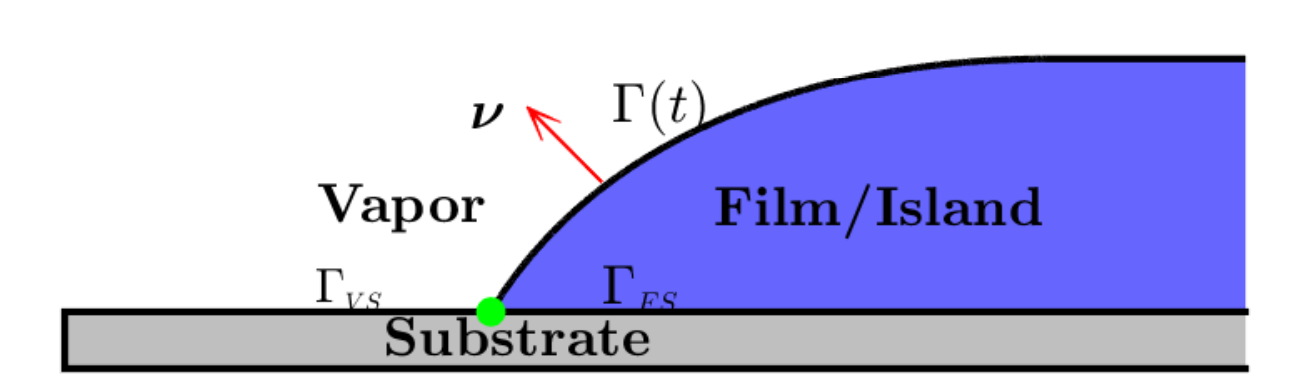}
    \caption{Sketch of the structure for SSD near the contact line (green point), where the vapor, film and substrate phases meet. }
    \label{fig:clm}
\end{figure}

When it comes to SSD, as shown in Fig.~\ref{fig:clm}, the total surface energy of the system consists of the film/vapor interface energy $\mathcal{E}_{{inf}}$ and the substrate energy $\mathcal{E}_{{ sub}}$, 
\begin{equation}
\mathcal{E}_{inf}=\int_{\Gamma(t)}\gamma(\bds{\nu})\,\dS,\qquad \mathcal{E}_{sub}= \gamma_{_\FS}\int_{\Gamma_{\FS}}\dS  + \gamma_{_\VS}\int_{\Gamma_{\VS}}\dS, 
\end{equation}
where $\Gamma(t)$ is the dynamic film/vapor interface with $\bds{\nu}$ being the interface normal pointing into the vapor phase, $\Gamma_{\FS}$ and $\Gamma_{\VS}$ are the interfaces between film/substrate and vapor/substrate, respectively, and $\gamma_{_\FS}$ and $\gamma_{_\VS}$ are the corresponding surface energy densities. 
In order to model SSD by the diffuse-interface approach, we associate the vapor phase with
$\varphi\approx1$ and the film phase with $\varphi\approx-1$.
Then the Ginzburg--Landau type energy \eqref{eq:aniGL}, up to a multiplicative constant, will approximate
the sharp interface energy $\mathcal{E}_{inf}$.
Moreover, the contribution to the wall energy $\mathcal{E}_{{ sub}}$ can be approximated by
\begin{equation}
\mathcal{E}_w(\varphi)=\int_{\Gamma_{\FS}\cup\Gamma_{\VS}}\frac{\gamma_{\VS}+ \gamma_{\FS}}{2} + (\gamma_{\VS} - \gamma_{\FS})G(\varphi)\,\dS,\label{eq:subE}
\end{equation}
where $G(\varphi)$ is a smooth function satisfying $G(\pm1)=\pm\frac{1}{2}$, see \cite{Jiang12,Dziwnik17anisotropic, Huang2019efficient,Backofen2017convexity} for SSD and \cite{jacqmin_2000,Qian2003molecular} for moving contact lines in fluid mechanics.  

There are several results on the existence of weak solutions for the degenerate Cahn--Hilliard equation \eqref{eq:isoch} with homogeneous  boundary conditions or its variants with inhomogeneous boundary conditions, see \cite{Elliott1996exist, Dai2016weak, Yin1992existence}. However, little is known about the anisotropic case except the work in \cite{dziwnik2019existence} which focuses on a particular $n$-fold anisotropy in two space dimensions. 

The main aim of this work is to develop a diffuse-interface approach to SSD in the case of anisotropic surface energies based on the energy contributions \eqref{eq:aniGL} and \eqref{eq:subE}. The obtained diffuse-interface model consists of a degenerate anisotropic Cahn--Hilliard equation with appropriate boundary conditions. We study the sharp-interface limit and show the existence of weak solutions to the diffuse-interface model.

The rest of the paper is organized as follows. In Section~\ref{sec:mod}, we review a sharp-interface model for SSD and then \revised{introduce} a diffuse-interface model based on a gradient flow approach. 
We then derive the sharp-interface limit from a regularized model with the help of asymptotic expansions in Section~\ref{se:asymptotic}. In Section~\ref{se:analysis}, we prove the existence of weak solutions to the diffuse-interface model. Numerical tests are presented in Section~\ref{se:nr}, where a comparison between sharp-interface approximations and diffuse-interface approximations is made.

\section{Modeling aspects} \label{sec:mod}

In this section, we first review a sharp-interface model for SSD with anisotropic surface energies
\revised{in two or three space dimensions}.
Then, we propose a suitable diffuse-interface model to approximate this sharp-interface model.
\revised{Here we note that there exist several works on the modelling of SSD using a diffuse-interface approach in the literature. However, these works consider either the isotropic case, e.g., \cite{Jiang18, Backofen2017convexity}, or the anisotropic case in 2d, e.g., \cite{Dziwnik17anisotropic}.}

\subsection{The sharp-interface model} \label{SUBSECT:SI}
We consider the dewetting of a solid thin film on a flat substrate in $\R^d$ with $d\in\{2,3\}$, as shown in Fig.~\ref{fig:clm}. We parameterize the interface of $\Gamma(t)$ over the initial interface as follows \[\vec x(\cdot, t)\,:\, \Gamma(0)\times[0,T]\to \R^d,\] where $T>0$ is a prescribed final time. The induced velocity is then given by
\begin{equation}
    \mathcal{\vec V}(\vec x(\vec q,t), t) =\partial_t\vec x(\vec q, t)\qquad\mbox{for all}\quad\vec q\in \Gamma(0),\quad t\in[0,T],\nn
\end{equation}
where $\Gamma(0)$ is a smooth hypersurface with boundary. 
The sharp-interface model for SSD (cf.\ \cite{Cahn1994overview, Taylor94linking, BGN2010finite,JZB2020sharp}) reads as:
\begin{subequations}
\label{eqn:sharpvk}
    \begin{align}
    \label{eq:sharpv}
        \mathcal{V}& = -\nabla_s\cdot( D(\bds{\nu})\,\nabla_s\varkappa_\gamma), \\[0.5em]
        \varkappa_\gamma &=-\nabla_s\cdot\gamma^\prime(\bds{\nu}),
        \label{eq:sharpk}
    \end{align}
\end{subequations} 
which has to hold for all $t\in[0, T]$ and all points on $\Gamma(t)$. Here, $\mathcal{V}=\mathcal{\vec V}\cdot\bds{\nu}$ is the normal velocity,  $\bds{\nu}$ is the unit normal to $\Gamma(t)$ pointing into the vapor, and $\nabla_s$ is the surface gradient operator on $\Gamma(t)$. Besides, $D(\bds{\nu})$ is an orientation 
dependent mobility (cf.~\cite{Taylor94linking}). The function $D$ needs to be defined for unit vectors, but here we extend its domain to $\R^d$ such that it is positively homogeneous
of degree one. The term $\varkappa_\gamma$ represents the anisotropic mean curvature, 
and $\gamma^\prime(\bds{\nu})$ is the Cahn--Hoffman vector, where $\gamma^\prime$ denotes the gradient of $\gamma$ (cf.~\cite{Hoffman1972vector}). The above equations are subject to the following boundary conditions at the contact line, where the film/vapor interface $\Gamma(t)$ meets the substrate: 
\begin{subequations}
\label{eqn:sharpbd}
\begin{itemize}
    \item attachment condition 
    \begin{equation}
        \mathcal{\vec V}\cdot\vec n_w = 0,
        \label{eq:sharpbd1}
    \end{equation}
    \item contact angle condition
\begin{equation}
    \gamma^\prime(\bds{\nu})\cdot\vec n_w + \sigma =0,
    \label{eq:sharpbd2}
\end{equation}
\item zero-flux condition
\begin{equation}
  D(\bds{\nu})\,\nabla_s\varkappa_\gamma\cdot \vec n_c=0,
  \label{eq:sharpbd3}
\end{equation}
\end{itemize}
\end{subequations}
where 
\begin{equation} \label{eq:sigma}
\sigma = \gamma_{_\VS} - \gamma_{_\FS} 
\end{equation}
denotes the difference of the substrate energy densities across the contact line. Here, $\vec n_w$ is the unit normal to the substrate and points in the direction of the substrate, and $\vec n_c$ is the conormal vector of $\Gamma(t)$, i.e., it is the outward unit normal to $\partial\Gamma(t)$ and it lies within the tangent plane of $\Gamma(t)$.
We observe that \eqref{eq:sharpbd2} enforces an angle condition between the
Cahn--Hoffman vector $\gamma^\prime(\bds{\nu})$ and the substrate unit normal 
$\vec n_w$. For example, in the isotropic case, $\gamma(\vec p) = |\vec p|$,
the Cahn--Hoffman vector reduces to the normal $\bds{\nu}$, and so if
$\sigma=0$ the condition \eqref{eq:sharpbd2} encodes a $90^\circ$ contact angle
between the film/vapor interface and the substrate.

We assume that the anisotropy function $\gamma$ belongs to $C^2\big(\R^d \setminus \{\vec 0\}\big) \cap C(\R^d,\R_{\geq0})$, is convex and satisfies $\gamma>0$ on $\R^d\setminus\{\mathbf 0\}$. We further assume that $\gamma$ is positively homogeneous of degree one, meaning that
\begin{align*}
    \gamma(\lambda\vec p) = \lambda \gamma(\vec p) \quad\text{for all $\lambda>0$, $\vec p\in\R^d$}.
\end{align*}
This immediately implies $\gamma(\mathbf 0) = 0$ and 
the gradient of $\gamma(\vec p)$ satisfies
\begin{align}
    \label{eq:gamma}
    \gamma'(\p)\cdot \p = \gamma(\p) \quad\text{for all $\p\in\R^d\setminus \{\vec 0\}$}.
\end{align}
Similarly, the orientation dependent mobility function $D\in C^2\big(\R^d \setminus \{\vec 0\}\big) \cap C(\R^d,\R_{\geq0})$ is assumed to satisfy 
$D>0$ on $\R^d\setminus \{\vec 0\}$ and
\begin{align*}
    D(\lambda\vec p) = \lambda D(\vec p) \quad\text{for all $\lambda>0$, $\vec p\in\R^d$}.
\end{align*}
Consequently, for the map 
\begin{align}
\label{DEF:A}
    A:\R^d \to \R, \quad \vec p \mapsto \tfrac 12 \gamma^2(\vec p)
\end{align}
introduced in \eqref{eq:aniGL}, we have $A\in C^2\big(\R^d \setminus \{\vec 0\}\big) \cap C(\R^d,\R_{\geq0})$. It also follows directly from \eqref{eq:gamma} that the relations 
\begin{subequations}
\label{eqn:useid}
\begin{alignat}{3}
        A(\lambda\p) &= \lambda^2 A(\p),
        \qquad
	    &A'(\vec p) &= \gamma(\vec p) \gamma'(\vec p), 
	    \qquad
	    &A'(\vec p) \cdot \vec p &= 2 A(\vec p),
	    \\[1ex]
	    A'(\lambda\vec p) &= \lambda A'(\vec p),
	    \qquad
	    &A''(\lambda\vec p) &= A''(\vec p),
	    \qquad
	    &A''(\vec p) \p &= A'(\vec p)
	\end{alignat}
\end{subequations}
hold for all $\p\in\R^d\setminus \{\vec 0\}$ and all $\lambda>0$. Here, $A^{\prime}$ and $A^{\prime\prime}$ denote the gradient and the Hessian of $A$, respectively.

\subsection{The diffuse-interface model}
\label{ssec:model}
\begin{figure}[!htp]
\centering
\includegraphics[width=0.6\textwidth]{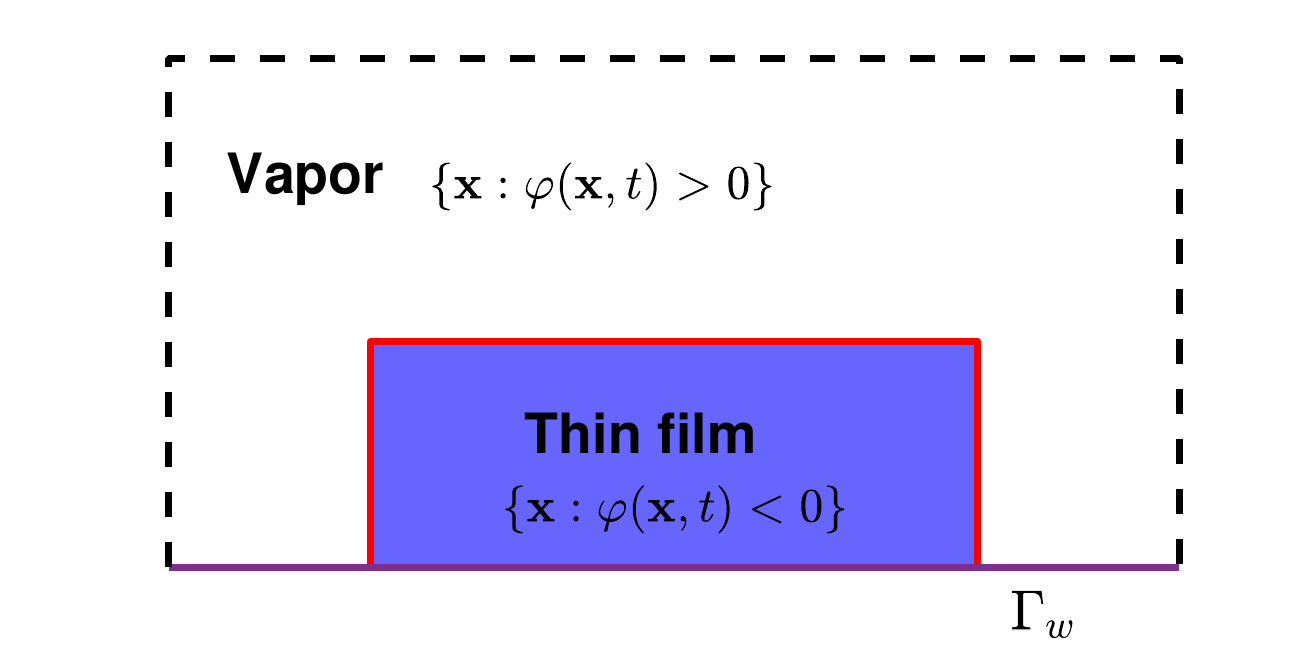}
\caption{Geometric setup for SSD in a bounded domain $\Omega$ with $\Omega=\overline{\Omega_-}(t)\cup\overline{\Omega_+}(t)$, where $\Omega_{-}(t):=\{\vec x\in\Omega: \varphi(\vec x,t)<0\}$ and $\Omega_{+}(t):=\{\vec x\in\Omega: \varphi(\vec x,t)>0\}$.}
\label{fig:model}
\end{figure}

Let $\varphi:\Omega\times[0, T]\to\R$ be an order parameter such that the zero level set $\{\vec x\in\Omega: \varphi(\vec x, t) =0\}$ approximates the film/vapor interface $\Gamma(t)$, $\{\vec x\in\Omega: \varphi(\vec x, t) < 0\}$ corresponds to the region occupied by the thin film at time $t$, whereas $\{\vec x\in\Omega: \varphi(\vec x, t) > 0\}$ represents the region occupied by the vapor at time $t$ (see~Fig.~\ref{fig:model}). In addition, $\Gamma_w\subset\partial\Omega$ models the boundary of the substrate.
As a combination of \eqref{eq:aniGL} and \eqref{eq:subE}, the total free energy of the system is given by 
\begin{align}
\mathcal{E}(\varphi):=&\frac{1}{\cF}\mathcal{E}_\gamma(\varphi) + \mathcal{E}_w(\varphi) - \gamma_{_\VS}|\Gamma_w|  \nn\\ =&\frac{1}{\cF}\int_{\Omega}\varepsilon A(\nabla\varphi) + \varepsilon^{-1} F(\varphi)\dx+\sigma\,\int_{\Gamma_w} G(\varphi)\,\dS - \frac{\sigma}{2}|\Gamma_w|,
\label{eq:Energy}
\end{align}
where $\cF=\int_{-1}^1\sqrt{2F(s)}\,\rd s$ and $|\Gamma_w|= \int_{\Gamma_w}\dS$. This choice of $\cF$ ensures that
\[\frac{1}{\cF}\mathcal{E}_\gamma(\varphi)\approx \int_{\Gamma(t)}\gamma(\bs{\nu})\,\dS\]
for sufficiently small $\eps>0$. Besides, the constant term $-\gamma_{_\VS}|\Gamma_w|$ was added to the total energy such that $\mathcal{E}(\varphi)$ now only depends on the single parameter $\sigma$ (see \eqref{eq:sigma}) instead of on $\gamma_{_\VS}$ and $\gamma_{_\FS}$. We next derive the diffuse-interface model. To this end, we use the smooth double-well potential
\begin{equation}
\label{DEF:F}
F(\varphi) = \frac{1}{2}(1-\varphi^2)^2.
\end{equation}
This implies
\begin{equation*}
    \cF=\int_{-1}^1\sqrt{2F(s)}\,\rd s = \frac{4}{3}.
\end{equation*}
We further choose
\begin{equation}
    G(\varphi) = \frac{1}{4}(3\varphi - \varphi^3),
    \label{DEF:G}
\end{equation}
which yields $G(\pm 1) = \pm \frac{1}{2}$ and $G^\prime(\pm 1) = 0$.
Let $\psi:\Omega\to\R$ be a sufficiently smooth function. Then the first variation of the total free energy \eqref{eq:Energy} in the direction of $\psi$ can be computed as 
\begin{align}
&\lim_{\delta\to 0}\frac{\mathcal{E}(\varphi+ \delta\psi)- \mathcal{E}(\varphi)}{\delta}= \frac{1}{\cF}\int_{\Omega} \varepsilon A'(\nabla\varphi)\cdot\nabla\psi  + \varepsilon^{-1} F^\prime(\varphi)\,\psi\dx  + \sigma\,\int_{\Gamma_w} G^\prime(\varphi)\,\psi\dS\nn\\
&\hspace{1cm} = \frac{1}{\cF}\int_{\Omega} [\varepsilon^{-1} F^\prime(\varphi) - \varepsilon\nabla\cdot A^\prime(\nabla\varphi)]\,\psi\dx + \frac{1}{\cF}\int_{\partial\Omega\setminus\Gamma_w} \varepsilon A^\prime(\nabla\varphi)\cdot\vec n\,\psi\,\dS\nn\\ &\hspace{1.5cm}+ \frac{1}{\cF}\int_{\Gamma_w}\left[\varepsilon A^\prime(\nabla\varphi)\cdot\vec n_w + \cF\,\sigma\, G^\prime(\varphi)\right]\,\psi\dS,
\label{eq:variation}
\end{align}
where $\vec n$ is the outward unit normal to $\partial\Omega\setminus\Gamma_w$ and $\vec n_w$ is the outward unit normal to $\Gamma_w$, as defined previously. 
The following diffuse-interface model for SSD can be interpreted as a weighted $H^{-1}$-gradient flow of the energy functional \eqref{eq:Energy}:
\begin{subequations}
\label{eqn:model}
\begin{alignat}{2}\label{eq:modela}
\alpha\,\partial_t\varphi &=\,\eps^{-1}\, \nabla\cdot\left(m(\varphi)\,\beta(\nabla\varphi)\,\nabla\mu\right),\quad&&\mbox{in}\quad Q=\Omega\times(0,T],\\ \mu &= -\varepsilon \nabla\cdot A'(\nabla\varphi) + \varepsilon^{-1} F^\prime(\varphi),\quad&&\mbox{in}\quad Q.
\label{eq:modelb}
\end{alignat}
\end{subequations}
Here, $\alpha>0$ is a time scaling coefficient, $m(\varphi)$ is the degenerate mobility given by 
\begin{equation}
m(\varphi) = (1-\varphi^2)_+^2=\left\{\begin{array}{ll}
2\,F(\varphi)\quad&\mbox{if}\;|\varphi|\leq 1,\\[0.5em]
0\quad&\mbox{otherwise},
\end{array}\right. 
\end{equation}
and $\beta(\nabla\varphi)$ is defined as
\begin{equation}
\beta(\nabla\varphi) = \frac{D(\nabla\varphi)}{\gamma(\nabla\varphi)},
\end{equation}
and so is positively homogeneous of degree zero.

We now write $\Sigma = \partial\Omega\times(0,T]$ and $\Sigma_w = \Gamma_w\times(0,T]$. On $\Sigma_w$, we impose the boundary conditions
\begin{subequations}
\label{eqn:bd}
\begin{align}\label{eq:bda}
m(\varphi)\,\beta(\nabla\varphi)\nabla\mu\cdot\vec n_w = 0,\qquad \varepsilon A'(\nabla\varphi)\cdot\vec n_w + \cF\,\sigma\, G^\prime(\varphi)=0.
\end{align}
Here the first equation is the zero-flux condition on the boundary, whereas the second equation guarantees the integral over $\Gamma_w$ in \eqref{eq:variation} vanishes.
Moreover, on $\Sigma\setminus\Sigma_w$, we impose the natural boundary conditions
\begin{align}\label{eq:bdb}
m(\varphi)\,\beta(\nabla\varphi)\nabla\mu\cdot\vec n = 0,\qquad A'(\nabla\varphi)\cdot\vec n = 0.
\end{align}
\end{subequations}

\begin{remark}
It is also possible to consider the double-obstacle potential \eqref{eq:doF} 
along with the mobility $m(\varphi) = (1-\varphi^2)_+$. 
Then the corresponding diffuse-interface model consists of \eqref{eq:modela} and the variational inclusion
\begin{equation}
\mu\in -\eps \nabla\cdot A^\prime(\nabla\varphi) + \eps^{-1}\partial F(\varphi).
\label{eq:muobexp}
\end{equation}
instead of \eqref{eq:modelb}.
\end{remark}

\section{The sharp-interface limit}
\label{se:asymptotic}

We consider the smooth double-well potential introduced in \eqref{DEF:F} and regularize the coefficients $m(\varphi)$ and $\beta(\nabla\varphi)$ of the diffuse-interface model \eqref{eqn:model} with the help of the interfacial parameter $\eps$ by defining 
\begin{align}
    \label{DEF:ME}
    m^\eps(\varphi) 
    &:= \eps^r + m(\varphi) = \eps^r + (1-\varphi^2)_+^2,
    \\
    \label{DEF:BE}
    \beta^\eps(\nabla\varphi)&:= \sqrt{ \frac{d_1^2\eps^{\revised{r}} + D^2(\Grad\varphi)}{\gamma_0^2 \eps^{\revised{r}} + \gamma^2(\Grad\varphi)}}
    \revised{,}
\end{align}
where \revised{$r\geq 2$. The regularized diffuse-interface model is then given by}
\begin{subequations}\label{eqn:regmod}
\begin{alignat}{2}
    &\label{eq:pm1}
\alpha\,\partial_t\varphi^\eps = \eps^{-1}\nabla\cdot\left(m^\eps(\varphi^\eps)\,\beta^\eps(\nabla\varphi^\eps)\,\nabla\mu^\eps\right)
    &&\quad\text{in $Q$},\\
    \label{eq:pm2}
    &\mu^\eps  = -\varepsilon \nabla\cdot A'(\nabla\varphi^\eps) + \varepsilon^{-1} F^\prime(\varphi^\eps)
    &&\quad\text{in $Q$},\\
    \label{eq:pbd1}
    &m^\eps(\varphi^\eps)\,\beta^\eps(\Grad\varphi^\eps) \Grad\mu^\eps\cdot \nw = 0
    &&\quad\text{on $\Sigma_w$},\\
    \label{eq:pbd2}
    &\eps A'(\Grad\varphi^\eps)\cdot \nw + \cF\,\sigma\, G^\prime(\varphi^\eps) = 0
    &&\quad\text{on $\Sigma_w$},\\
    \label{eq:pbd3}
    &m^\eps(\varphi^\eps)\,\beta^\eps(\Grad\varphi^\eps) \Grad\mu^\eps\cdot \n = 0
    &&\quad\text{on $\Sigma\setminus\Sigma_w$},\\
    \label{eq:pbd4}
    &\eps A'(\Grad\varphi^\eps)\cdot \n  = 0
    &&\quad\text{on $\Sigma\setminus\Sigma_w$}.
\end{alignat}
\end{subequations}
\revised{We note that the introduction of the three regularization terms $\eps^r$ in
\eqref{DEF:ME} and \eqref{DEF:BE} allows for a mathematical analysis of \eqref{eqn:regmod} in
Section~\ref{se:analysis} below. In fact, on defining}
\begin{alignat*}{2}
    \gamma_0 &:= \underset{|\vec p|=1}{\min}\; \gamma(\vec p) > 0,
    &\qquad
    \gamma_1 &:= \underset{|\vec p|=1}{\max}\; \gamma(\vec p) > 0,
    \\
    d_0 &:= \underset{|\vec p|=1}{\min}\; D(\vec p) > 0,
    &\qquad
    d_1 &:= \underset{|\vec p|=1}{\max}\; D(\vec p) > 0,
\end{alignat*}
\revised{we} have
\begin{align*}
    \eps^r \le m^\eps(\varphi) \le \eps^r + 1
    \qquad\text{and}\qquad
    \frac{d_0}{\gamma_1} \le \beta^\eps(\Grad\varphi) \le \frac{d_1}{\gamma_0}.
\end{align*}
\revised{Moreover, by choosing $r\geq2$ we ensure that the sharp interface
limit of \eqref{eqn:regmod} is unchanged compared to the limit of \eqref{eqn:model}.}

We now formally derive the sharp-interface limit of the regularized model \eqref{eqn:model} via the method of matched asymptotic expansions. We suppose that for $\eps>0$, $(\varphi^\eps, \mu^\eps)$ is the solution of the regularized diffuse-interface model \revised{\eqref{eqn:regmod}. Then} we write
\begin{align}\label{eq:Geps}
\Gamma^\varepsilon(t):=\bigl\{\vec x\in\Omega\;|\; \varphi^\eps(\vec x, t)=0\bigr\}\qquad\mbox{and}\qquad
\Lambda^\varepsilon(t):=\Gamma^\eps(t)\cap\Gamma_w
\end{align}
to denote the interface and the contact line, respectively.
We further assume that their limits as $\varepsilon\to 0$ are given by $\Gamma(t)$ and $\Lambda(t)$, respectively. 
We introduce a local parameterization for $\Gamma(t)$ on an open subset $\mathcal{O}\subset\R^{d-1}$ by
\begin{equation}
\vec r(\vec s, t):\mathcal{O}\times[0,T]\to \mathbb{R}^d.
\end{equation}
\revised{Our asymptotic analysis for the interface dynamics will follow similar techniques in the literature for degenerate Cahn-Hilliard equations, see e.g.\ \cite{Cahn1996cahn, Dai2014coarsen} for the isotropic case and \cite{Ratz2006surface, Dziwnik17anisotropic} for the anisotropic case in 2d.}

\subsection{Outer expansions}
Away from the interface and the contact line, we assume that the following ansatz holds
\begin{subequations}\label{eqn:outexp}
\begin{align}
\varphi^\eps(\vec x, t) &= \varphi_0(\vec x, t) +\eps \varphi_1(\vec x, t) + \eps^2\varphi_2(\vec x, t) + \cdots,\\
\mu^\eps(\vec x, t) &= \eps^{-1}\,\mu_{-1}(\vec x, t)+\mu_0(\vec x, t) + \eps\mu_1(\vec x, t) + \eps^2\mu_2(\vec x, t) + \cdots.
\end{align}
\end{subequations}
Moreover, in view of \eqref{DEF:ME} and \eqref{DEF:BE}, we know that
\begin{subequations}\label{eqn:mb}
\begin{align}
    m^\eps(\varphi^\eps) &= m(\varphi_0) +\eps\, m^\prime(\varphi_0)\,\varphi_1 + O(\eps^2),\\
    \beta^\eps(\nabla\varphi^\eps)&= \beta(\nabla\varphi_0) +  \eps\,
    \beta^\prime(\nabla\varphi_0)\cdot\nabla\varphi_1 + O(\eps^2),
\end{align}
\end{subequations}
\revised{since $r\geq2$ and $\beta^\eps(\vec p) = \beta(\vec p) + O(\eps^r)$}, where $\beta^\prime$ denotes the gradient of $\beta$. Plugging the expansions \eqref{eqn:outexp} and \eqref{eqn:mb} into \eqref{eq:pm1} and \eqref{eq:pm2} gives
\begin{align*}
    0 = \nabla\cdot(\beta(\nabla\varphi_0)\,m(\varphi_0)\,\nabla\mu_{-1}),\qquad \mu_{-1} = F^\prime(\varphi_0).
\end{align*}
As the energy \eqref{eq:aniGL} is expected to be bounded at leading order, it needs to hold $F(\varphi_0) = 0$. This means that $\varphi_0$ attains only the values $-1$ and $1$. Hence, $\mu_{-1}=0$.
We now define 
\begin{align}
\Omega_+(t):=\bigl\{\vec x\in\Omega\;|\;\varphi_0(\vec x, t) = 1\bigr\},\qquad
\Omega_-(t):=\bigl\{\vec x\in\Omega\;|\;\varphi_0(\vec x, t) = -1\bigr\},\nn
\end{align}
as the outer regions, meaning that $\varphi_0 = \pm 1$ in $\Omega_\pm(t)$.

\subsection{Inner expansions}
In the inner region near the interface $\Gamma(t)$, we introduce the annular neighbourhood
\begin{align}
\mathcal{N}(t):=\Bigl\{\vec x\in\Omega\; :\;|{\rm d}(\vec x, t)|<\delta \Bigr\},\qquad \delta >0,\nn
\end{align}
where  ${\rm d}(\vec x, t)$ represents the signed distance of $\vec x $ to $\Gamma(t)$, defined to be positive in $\Omega_+(t)$. Assuming $\Gamma(t)$ to be sufficiently smooth, we find a $\delta>0$ such that for every $\vec x\in\mathcal{N}(t)$, there exist unique vectors $\vec r(\vec x, t)$ and $\vec s(\vec x, t)$ such that
\begin{align}
\vec x = \vec r(\vec s(\vec x, t), t) + \rd(\vec x, t)\,\bs{\nu}(\vec s(\vec x, t), t).
\end{align}
Here $\bs{\nu}(\vec s ,t)$ is the unit normal vector on $\Gamma(t)$ 
at $\vec r(\vec s, t)$ pointing into $\Omega_+(t)$. 

Due to rapid changes of $\varphi^\eps$ in normal direction, we introduce the stretched variable $\rho(\vec x, t) = \eps^{-1} \rd(\vec x, t)$. Any scalar function $b(\vec x, t)$ can be expressed in the new coordinate system as $b(\vec x, t) = \overline{b}(\vec s(\vec x, t),\rho(\vec x, t),t)$.
For any vector field $\vec b(\vec x, t)$, we use an analogous notation.
Without loss of generality, we assume that $\{\vec t_1,~\vec t_2,\cdots, \vec t_{d-1}\}$ forms an orthonormal basis of the tangent space of $\Gamma(t)$ at the point $\vec r(\vec s, t)$ such that
\begin{equation}
 \vec t_i\cdot\vec t_j = \delta_{ij},\qquad    \vec t_j = \partial_{s_j}\vec r,\qquad \partial_{s_j}\vec t_j = - \kappa_j\,\bs{\nu}\quad\mbox{for}\quad \vec s = (s_1,s_2,\cdots,s_{d-1})^T,\nn
\end{equation}
where $\kappa_j$ is the principal curvature of $\Gamma(t)$ at the point $\vec r(\vec s, t)$ in the direction of $\vec t_j$. As in \cite{Dai2014coarsen}, we obtain the identities 
\begin{equation*}
\nabla\rd = \eps\,\nabla\rho= \bs{\nu}(\vec s, t) \quad\text{and}\quad
\nabla s_j = \frac{1}{1+\eps\,\rho\,\kappa_j(\vec s, t)}\vec t_j(\vec s, t) 
,\quad\ 1\leq j\leq d-1, \quad \text{in } \mathcal{N}(t).
\end{equation*}
Therefore, using the new coordinates, we calculate
\begin{subequations}\label{eqn:coordT}
\begin{align}
\partial_t b &= \partial_t\overline{b} + \sum_{j=1}^{d-1}\partial_{s_j}\overline{b}\,\partial_t s_j + \partial_\rho\overline{b}\,\partial_t\rho =  \partial_t^\Gamma\overline{b}-\eps^{-1}\mathcal{V}\,\partial_\rho\overline{b},\\
\nabla b &=\nabla\rho\,\partial_\rho\overline{b} + \sum_{j=1}^{d-1}\partial_{s_j}\overline{b}\nabla\,s_j=\eps^{-1}\,\bs{\nu}\,\partial_\rho\overline{b}  + \nabla_{s}\overline{b}+O(\eps),\\
\nabla\cdot\vec b &=\nabla\rho\cdot\partial_\rho\overline{\vec b}+\sum_{j=1}^{d-1}\partial_{s_j}\overline{\vec b}\cdot\nabla\,s_j=\eps^{-1}\bs{\nu}\cdot\partial_\rho\overline{\vec b} + \nabla_s\cdot\overline{\vec b} +O(\eps),
\end{align}
\end{subequations}
where $\nabla_s=\sum_{j=1}^{d-1}\vec t_j\,\partial_{s_j}$ denotes the surface gradient operator on $\Gamma(t)$,  
$$\partial_t^\Gamma\overline{b}= \partial_t\overline{b} + \sum_{j=1}^{d-1}\partial_{s_j}\overline{b}\,\partial_t s_j,$$ and $\mathcal{V}$ is the velocity of $\Gamma(t)$ in the direction of $\bs{\nu}$, i.e., $\mathcal{V} = -\partial_t{\rm d} = - \eps\,\partial_t\rho$.

In the inner region, we assume the following expansions
\begin{subequations}\label{eqn:inner}
\begin{align}
\label{eq:inner1}
\varphi^\eps &= \Phi_0(\vec s, \rho, t) + \eps\Phi_1(\vec s, \rho, t) +\eps^2 \Phi_2(\vec s, \rho, t) +\cdots,\\
\mu^\eps&=\eps^{-1}\,M_{-1}(\vec s, \rho, t) + M_0(\vec s, \rho, t) + \eps M_1(\vec s, \rho, t) +\eps^2 M_2(\vec s, \rho, t) + \cdots.
\label{eq:inner2}
\end{align}
\end{subequations}
In particular, on assuming $\partial_\rho\Phi_0>0$, we have, \revised{similarly to \eqref{eqn:mb}, that}
\begin{subequations}\label{eqn:mb2}
\begin{align}
    m^\eps(\nabla\varphi^\eps) = m(\Phi_0) + \eps\,m^\prime(\Phi^0)\,\Phi_1 + O(\eps^2),
    \\
    \beta^\eps(\nabla\varphi^\eps) = \beta(\bs{\nu}) + \eps\, \beta^\prime(\bs{\nu})\cdot\nabla_s\Phi_0 + O(\eps^2),
\end{align}
\end{subequations}
where we have used the fact that \revised{$\beta$} is positively homogeneous of order zero.

Plugging \eqref{eqn:inner} and \eqref{eqn:mb2} into \eqref{eq:pm1}, we obtain the leading order term
\begin{equation}
   0 = \partial_\rho\big(\beta(\bs{\nu}) m(\Phi_0) \partial_\rho M_{-1}\big), 
\end{equation}
which implies that $m(\Phi_0)\partial_\rho M_{-1}$ is independent of $\rho$, i.e., it can be expressed as
\begin{equation*}
    m(\Phi_0)\,\partial_\rho\,M_{-1} = J(\vec s,t) .
\end{equation*}
In addition, using the matching condition 
\begin{align}
\lim_{\rho\to \pm\infty} \Phi_0(\rho) = \pm 1,\label{eq:phimat}
\end{align}
we infer $J(\vec s,t) = 0$ due to the degenerate mobility $m(\Phi_0)$. Since $m(s)>0$ if $s\in(-1,1)$, we thus conclude that $M_{-1}$ is independent of $\rho$. 
By the matching condition  $\lim_{\rho\to\pm\infty}M_{-1}(\vec s,t) = \mu_{-1}$, we obtain 
\begin{equation*}
    M_{-1} = M_{-1}(\vec s,t) \equiv 0.
\end{equation*}
For the terms of order $O(\frac{1}{\eps^3})$, we obtain
\begin{equation}
0 = \partial_\rho\left(\beta(\bs{\nu}) m(\Phi_0) \partial_\rho M_0\right).\nn
\end{equation}
Repeating the above line of argument, we deduce
\begin{equation}
\partial_\rho M_0 =0,\qquad M_0 = M_0(\vec s, t).
\end{equation}
Using the fact that $M_{-1}=0$ and $\partial_\rho M_0 = 0$, we then have the following expansions 
\begin{align}
&\nabla\cdot(\beta^\eps(\nabla\varphi^\eps)\,m^\eps(\varphi^\eps)\nabla\mu^\eps) 
    \notag\\
    &\quad =
    \frac{1}{\eps}\,\partial_\rho(\beta(\bs{\nu})\,m(\Phi_0)\,\partial_\rho M_1)\label{eq:expans}
    \notag\\ &\qquad
+\partial_\rho(\beta^\prime(\bs{\nu})\cdot\nabla_s\Phi_0\,m(\Phi_0)\,\partial_\rho M_1 + \beta(\bs{\nu})m^\prime(\Phi_0)\Phi_1\,\partial_\rho M_1)
    \notag\\ &\qquad
+\partial_\rho(\beta(\bs{\nu})\,m(\Phi_0)\,\partial_\rho\,M_2) + \nabla_s\cdot(\beta(\bs{\nu})\,m(\Phi_0)\,\nabla_s M_0)+O(\eps).
\end{align}
Considering the order $O(\frac{1}{\eps^2})$ of \eqref{eq:pm1}, we obtain that 
\begin{equation}
0 = \partial_\rho\left(\beta(\bs{\nu})\,m(\Phi_0)\,\partial_\rho M_1\right).
\end{equation}
Similarly, by using the matching conditions we arrive at
\begin{equation}
M_1 = M_1(\vec s,t).
\end{equation}
At $O(\frac{1}{\eps})$, using $\partial_\rho M_1=0$ and \eqref{eq:expans}, we have
\begin{equation}
-\alpha\,\mathcal{V}\partial_\rho\Phi_0 = \partial_\rho\left(\beta(\bs{\nu})\,m(\Phi_0)\,\partial_\rho M_2\right) + \nabla_s\cdot\left(\beta(\bs{\nu})\,m(\Phi_0)\,\nabla_s M_0\right).
\label{eq:vel}
\end{equation}

 We next consider the expansion of \eqref{eq:pm2}. Using the identities in \eqref{eqn:useid} and assuming $\partial_\rho\Phi_0>0$, we expand the anisotropic term $A^\prime(\nabla\Phi^\eps)$ as follows:
\begin{align}
A'(\nabla\Phi^\eps) &= A^\prime\left(\frac{1}{\eps}\partial_\rho\Phi^\eps\,\bs{\nu} + \nabla_s\Phi^\eps + O(\eps)\right)= \frac{1}{\eps}\partial_\rho\Phi^\eps\,A^\prime(\bs{\nu}) + A^{\prime\prime}(\bs{\nu})\nabla_s\Phi^\eps + O(\eps).\nn
\end{align}
This then yields 
\begin{align}
&\nabla\cdot A^\prime(\nabla\Phi^\eps) = \frac{1}{\eps}\partial_\rho[A^\prime(\nabla\Phi^\eps)]\cdot\bs{\nu} + \nabla_s\cdot A^\prime(\nabla\Phi^\eps)\nn\\
&\hspace{1cm}= \frac{1}{\eps}\partial_\rho\left(\frac{1}{\eps}\partial_\rho\Phi^\eps\,2A(\bs{\nu}) + A^\prime(\bs{\nu})\cdot\nabla_s\Phi^\eps\right) + \nabla_s\cdot\left( \frac{1}{\eps}\partial_\rho\Phi^\eps A^\prime(\bs{\nu}) + A^{\prime\prime}(\bs{\nu})\nabla_s\Phi^\eps \right)+ O(\eps)\nn\\
& \hspace{1cm}= \frac{2}{\eps^2}\partial_{\rho\rho}\Phi^\eps\,A(\bs{\nu}) + \frac{1}{\eps}\left(A^\prime(\bs{\nu})\cdot\partial_\rho(\nabla_s\Phi^\eps) + \nabla_s\cdot(\partial_\rho\Phi^\eps\,A^\prime(\bs{\nu}))\right)+ O(1).\nn
\end{align}
Now, plugging \eqref{eqn:inner} into \eqref{eq:pm2}, we obtain for the leading order term that
\begin{equation}
2A(\bs{\nu})\partial_{\rho\rho}\Phi_0 - F^\prime(\Phi_0) =M_{-1} =0. \label{eq:phie}
\end{equation}
Using the translation identity $\Phi_0(0) = 0$, we then obtain 
\begin{equation}
\label{eq:phi0}
\Phi_0(\rho) = \tanh\left(\frac{\rho}{\gamma(\bs{\nu})}\right),\quad-\infty<\rho<+\infty.
\end{equation}
Similarly, the $O(1)$ term resulting from \eqref{eq:pm2} implies
\begin{equation}
\label{eq:m0}
2A(\bs{\nu})\partial_{\rho\rho}\Phi_1 + A^\prime(\bs{\nu})\cdot\partial_\rho(\nabla_s\Phi_0) +\nabla_s\cdot(\partial_\rho\Phi_0\,A^\prime(\bs{\nu})) - F^{\prime\prime}(\Phi_0)\Phi_1 = - M_0(\vec s,t).
\end{equation}
Multiplying \eqref{eq:m0} by $\partial_\rho\Phi_0$ and then integrating from $-\infty$ to $\infty$ with respect to $\rho$ yields 
\begin{align}
& \int_{-\infty}^{+\infty}\left(A^\prime(\bs{\nu})\cdot\partial_\rho(\nabla_s\Phi_0) +\nabla_s\cdot(\partial_\rho\Phi_0\,A^\prime(\bs{\nu}))\right)\partial_\rho\Phi_0\drho\nn\\
 &\hspace{1cm}+\int_{-\infty}^{+\infty}\left(2A(\bs{\nu})\partial_{\rho\rho}\Phi_1\partial_\rho\Phi_0 - F^{\prime\prime}(\Phi_0)\Phi_1\partial_\rho\Phi_0\right)\rd\rho= -M_0\int_{-\infty}^{+\infty}\partial_\rho\Phi_0\drho. 
\label{eq:mi}  
\end{align}
Differentiating \eqref{eq:phie} with respect to $\rho$ gives 
\begin{equation}
2A(\bs{\nu})\partial_{\rho\rho\rho}\Phi_0 - F^{\prime\prime}(\Phi_0)\partial_\rho\Phi_0 = 0.\nn
\end{equation}
Therefore, since $\lim_{\rho\to\pm\infty}\partial_\rho\Phi_0=0$ and $\lim_{\rho\to\pm\infty}\Phi_1=0$, we compute 
\begin{align}
&\int_{-\infty}^{+\infty}\left(2A(\bs{\nu})\partial_{\rho\rho}\Phi_1\partial_\rho\Phi_0 - F^{\prime\prime}(\Phi_0)\Phi_1\partial_\rho\Phi_0\right)\rd\rho\nn
\\ &\quad
= \int_{-\infty}^{+\infty} (2A(\bs{\nu})\partial_{\rho\rho\rho}\Phi_0 - F^{\prime\prime}(\Phi_0)\partial_\rho\Phi_0)\,\Phi_1\drho =0\nn
\end{align}
via integration by parts.
Then, using \eqref{eq:phi0} and the matching condition in \eqref{eq:phimat}, we can reformulate \eqref{eq:mi} as
\begin{equation}
\int_{-\infty}^{+\infty}\left(A^\prime(\bs{\nu})\cdot\partial_\rho(\nabla_s\Phi_0) +\nabla_s\cdot(\partial_\rho\Phi\,A^\prime(\bs{\nu}))\right)\partial_\rho\Phi_0\drho 
= -  2M_0(\vec s,t). 
\end{equation}
It further follows from \eqref{eq:phi0} that $\partial_\rho(\nabla_s\Phi_0)= \nabla_s(\partial_\rho\Phi_0)$. We thus have
\begin{equation}
\int_{-\infty}^{+\infty}\nabla_s\cdot[A^\prime(\bs{\nu})(\partial_\rho\Phi_0)^2]\drho 
= -  2M_0(\vec s,t),\nn
\end{equation}
which yields
\begin{align}
M_0(s,t) = -\frac{1}{2}\nabla_s\cdot\left(A^\prime(\bs{\nu})\int_{-\infty}^{+\infty}(\partial_\rho\Phi_0)^2\drho\right) 
= -\frac{1}{2}\cF\nabla_s\cdot\gamma^\prime(\bs{\nu}) 
= \frac{1}{2}\cF\varkappa_\gamma,
\end{align}
where $\varkappa_\gamma = - \nabla_s\cdot\gamma^\prime(\bs{\nu})$ is the weighted mean curvature defined in \eqref{eq:sharpk}.

We now return to \eqref{eq:vel} and integrate it with respect to $\rho$ from $-\infty$ to $+\infty$. Using the fact that $\lim_{\rho\to\pm\infty}m(\Phi_0)\partial_\rho M_2= 0$, we get
\begin{equation}
- 2\,\alpha\,\mathcal{V} = \nabla_s\cdot\left(\beta(\bs{\nu})\int_{-\infty}^\infty m(\Phi_0)\drho\,\nabla_s M_0\right) = \cF\nabla_s\cdot[D(\bs{\nu})\nabla_sM_0],\nn
\end{equation}
where we recall \eqref{eq:phi0} and also use the identities 
\begin{equation}
\beta(\bs{\nu})=\frac{D(\bs{\nu})}{\gamma(\bs{\nu})}
\quad\text{and}\quad
\int_{-\infty}^{+\infty}m(\Phi_0)\drho =\int_{-\infty}^{+\infty}2\,F(\Phi_0(\rho))\drho = \cF\gamma(\bs{\nu}).\nn
\end{equation}
We thus obtain
\begin{align}
\mathcal{V} = -\frac{\cFsqu}{4\,\alpha}\nabla_s\cdot[D(\bs{\nu})\nabla_s\varkappa_\gamma]\quad\mbox{with}\quad\varkappa_\gamma = -\nabla_s\cdot\gamma^\prime(\bs{\nu}).
\label{eq:anisf}
\end{align}

\subsection{Expansions near the intersection with the substrate}

 We next study the expansions near the intersection with the substrate using the technique discussed in \cite{Dziwnik17anisotropic, owen1990minimizers}. 
 
 \subsubsection{The boundary layer near the wall}
 In the boundary layer near $\Gamma_w$, we first introduce the variable $\eta = \eps^{-1}\,{\rm d}_w(\vec x)$, where ${\rm d}_w(\vec x)$ represents the distance from $\vec x$ to the wall $\Gamma_w$. Then for a scalar function $b(\vec x, t)$, we can write it as $b(\vec x, t) = \widehat{b}(\eta, \vec y, t)$, 
where $\vec y$ is the $(d-1)$-dimensional coordinate system that is orthogonal to $\eta$. This implies \[\nabla b = \nabla_{\vec y}\widehat{b} - \eps^{-1}\partial_\eta\widehat{b}\,\vec n_w.\] We consider the expansions 
\begin{align}
\varphi^\eps &= \widehat{\varphi}_0(\eta, \vec y, t) + \eps \widehat{\varphi}_1(\eta, \vec y, t) + \eps^2\widehat{\varphi}_2(\eta, \vec y, t)+ \cdots,\\
\mu^\eps &= \widehat{\mu}_0(\eta, \vec y,t) + \eps \widehat{\mu}_1(\eta, \vec y,t) + \eps^2\widehat{\mu}_2(\eta, \vec y,t)+\cdots,
\end{align}
and plug them into \eqref{eq:pm1} and \eqref{eq:pm2}. The leading order terms  yield
\begin{subequations}
\begin{align}
\label{eq:suba}
\partial_\eta\left(\widehat{\beta}_0\,m(\widehat{\varphi}_0)\,\partial_\eta \widehat{\mu}_0\right) = 0,\\[0.3em]
\partial_\eta[A^\prime(-\partial_\eta\widehat{\varphi}_0\,\vec n_w)]\cdot\vec n_w + F^\prime(\widehat{\varphi}_0)=0,
\label{eq:subb}
\end{align}
\end{subequations}
where $\widehat{\beta}_0 =\beta(-\partial_\eta\widehat{\varphi}_0\,\vec n_w)$.  At the boundary $\eta=0$, it holds
\begin{subequations}
\begin{align}
\label{eq:ssuba}
- A^\prime(-\partial_\eta\widehat{\varphi}_0\,\vec n_w)\cdot\vec n_w + \cF\,\sigma\, G^\prime(\widehat{\varphi}_0) = 0,\\ \widehat{\beta}_0\,m(\widehat{\varphi}_0)\partial_\eta \widehat{\mu}_0 = 0.
\label{eq:ssubb}
\end{align}
\end{subequations}
Thus from \eqref{eq:suba} and \eqref{eq:ssubb} we obtain
\begin{equation}
m(\widehat{\varphi}_0)\partial_\eta \widehat{\mu}_0 = 0.\nn
\end{equation}
Multiplying \eqref{eq:subb} by $\partial_\eta\widehat{\varphi}_0$ and using the identities in \eqref{eqn:useid}, we arrive at
\begin{align}
0&=-\partial_{\eta\eta}\,\widehat{\varphi}_0\,\vec n_w\cdot A^{\prime\prime}(-\partial_\eta\widehat{\varphi}_0\,\vec n_w)\,\partial_\eta\widehat{\varphi}_0\,\vec n_w + F^\prime(\widehat{\varphi}_0)\partial_\eta\widehat{\varphi}_0\nn\\
 &= \partial_{\eta\eta}\widehat{\varphi}_0\,\vec n_w\cdot A^\prime(-\partial_\eta\widehat{\varphi}_0\,\vec n_w) + F^\prime(\widehat{\varphi}_0)\partial_\eta\widehat{\varphi}_0.
\label{eq:subi}
\end{align}
 Integrating \eqref{eq:subi} over $\eta$ leads to 
 \begin{equation}
\label{eq:bldf}
 F(\widehat{\varphi}_0)  = A (-\partial_\eta\widehat{\varphi}_0\,\vec n_w) + c(\vec y, t)= (\partial_\eta\widehat{\varphi}_0)^2\,A(-\vec n_w),
 \end{equation}
 where $c(\vec y, t)=0 $ due to the matching condition when $\eta\to\infty$. This implies
\begin{equation}
\label{eq:bldd}
\partial_\eta\widehat{\varphi}_0 =\left\{\begin{array}{ll}
 -\sqrt{\frac{F(\widehat{\varphi}_0)}{A(-\vec n_w)}}\quad\mbox{if}\quad \partial_\eta\widehat{\varphi}_0<0,\\[0.4em]
  +\sqrt{\frac{F(\widehat{\varphi}_0)}{A(-\vec n_w)}}\quad\mbox{if}\quad \partial_\eta\widehat{\varphi}_0>0.
 \end{array}\right.
\end{equation}

\subsubsection{The inner layer near the contact line} 
We assume that a local parameterization of the contact line $\Lambda(t)$ is given by
\begin{equation}
 \vec r_w(s_w,t): \mathcal{O}_w\times[0,T]\to\mathbb{R}^d,
\end{equation}
where in the case $d=2$, we simply set $\mathcal{O}_w = \{0\}$. For a contact point $\vec x_c\in\Lambda(t)$, we then introduce an interior layer near it. Precisely, for any $\vec x$ in the plane that contains $\vec x_c$ and is spanned by $\vec n_s$ and $\vec n_w$, we write 
\begin{align}
\xi = \eps^{-1}\,(\vec x - \vec x_c)\cdot\vec n_s, \qquad \eta = -\eps^{-1}\,(\vec x -\vec x_c)\cdot\vec n_w, \nn
\end{align}
where $\vec n_s$ is the unit normal to $\Lambda(t)$ on the wall $\Gamma_w$ and pointing into $\Omega_+(t)$. For a scalar function $b(\vec x, t)$, we can rewrite it as  $b(\vec x, t) =  \widetilde{b}(s_w, \xi, \eta, t)$. In a similar manner to \eqref{eqn:coordT}, we compute 
\begin{align}
\partial_t b &= -\eps^{-1}\partial_\xi\widetilde{b}\, (\partial_t\vec x_c\cdot\vec n_s) + \partial_t^\Lambda \widetilde{b},\\
\nabla b &= \eps^{-1}\,(\partial_\xi\widetilde{b}\,\vec n_s - \partial_\eta\widetilde{b}\,\vec n_w) + \nabla_{s_w} \widetilde{b} + O(\eps),\\
\nabla\cdot\vec b &= \eps^{-1}\,(\partial_\xi\widetilde{\vec b}\cdot\vec n_s - \partial_\eta\widetilde{\vec b}\cdot\vec n_w) + \nabla_{s_w}\cdot\widetilde{\vec b} + O(\eps),
\end{align}
where $\nabla_t^\Lambda\widetilde{b} = \partial_t\widetilde{b} +\partial_ts_w\cdot\nabla_{s_w}\widetilde{b}$.
We then consider the expansions 
\begin{align}
\varphi^\eps = \widetilde{\varphi}_0(s_w,\xi,\eta,t) + \eps \widetilde{\varphi}_1(s_w,\xi,\eta,t) + \eps^2\widetilde{\varphi}_2(s_w,\xi,\eta,t) + \cdots,\\
\mu^\eps  = \widetilde{\mu}_0(s_w,\xi,\eta,t) + \eps \widetilde{\mu}_1(s_w,\xi,\eta,t) + \eps^2\widetilde{\mu}_2(s_w,\xi,\eta,t) + \cdots,
\end{align}
and plug them into \eqref{eq:pm1} and \eqref{eq:pm2}. By defining $\nabla_c = \vec n_s\,\partial_\xi-\vec n_w\,\partial_\eta$, the leading order term yields 
\begin{subequations}
\begin{align}
\label{eq:wallld1}
\nabla_c\cdot\left(\widetilde{\beta}_0\,m(\widetilde{\varphi}_0)\nabla_c \widetilde{\mu}_0\right) &=0,\\[0.3em]
\label{eq:wallld2}
\partial_\xi\bigl( A^\prime(\nabla_c\widetilde{\varphi}_0)\cdot\vec n_s\bigr) - \partial_\eta\bigl(A^\prime(\nabla_c\widetilde{\varphi}_0)\cdot\vec n_w\bigr)&= F^\prime(\widetilde{\varphi}_0),
\end{align}
\end{subequations}
where  $\widetilde{\beta}_0=\beta(\nabla_c\widetilde{\varphi}_0)$. Similarly, the leading order terms of the boundary conditions \eqref{eq:pbd1} and \eqref{eq:pbd2} give
\begin{subequations}
\begin{align}
\label{eq:wbd1}
\widetilde{\beta}_0\,m(\widetilde{\varphi}_0)\partial_\eta \widetilde{\mu}_0 = 0,\\
 A^\prime(\nabla_c\widetilde{\varphi}_0)\cdot\vec n_w + \cF\,\sigma\,G^\prime(\widetilde{\varphi}_0)=0.
\label{eq:wbd2}
\end{align}
\end{subequations}
Besides, we have the matching condition 
\begin{align}
\label{eq:matchbl}
\lim_{\xi\to\pm\infty}\widetilde{\varphi}_0 = \lim_{\vec y\to \vec y (x_c^\pm)} \widehat{\varphi}_0(\vec y,\eta)=\widehat{\varphi}_0^\pm.
\end{align}
Now, multiplying \eqref{eq:wallld2} by $\partial_\xi\widetilde{\varphi}_0$ and 
integrating the resulting equation in a box $R:=[-\xi_1,\xi_1]\times[0,\eta_1]$, we get
\begin{align}
\int_{-\xi_1}^{\xi_1}\int_{0}^{\eta_1}\partial_{\xi}\widetilde{\varphi}_0\,\left[\partial_\xi\bigl( A^\prime(\nabla_c\widetilde{\varphi}_0)\cdot\vec n_s\bigr) - \partial_\eta\bigl(A^\prime(\nabla_c\widetilde{\varphi}_0)\cdot\vec n_w\bigr)\right]\deta\rd\xi = \int_{-\xi_1}^{\xi_1}\int_{0}^{\eta_1}\partial_\xi\widetilde{\varphi}_0\, F^\prime(\widetilde{\varphi}_0)\deta\rd\xi,\nn
\end{align}
which can be rewritten as
\begin{align}
 &\int_{0}^{\eta_1}\int_{-\xi_1}^{\xi_1}\partial_{\xi}\left[F(\widetilde{\varphi}_0) + A(\nabla_c\widetilde{\varphi}_0) - \partial_{\xi}\widetilde{\varphi}_0\,A^\prime(\nabla_c\widetilde{\varphi}_0)\,\cdot\vec n_s\right]\rd\xi\rd\eta\nn\\ &\hspace{2cm}+\int_{-\xi_1}^{\xi_1}\int_{0}^{\eta_1}\partial_{\eta}\left[\partial_{\xi}\widetilde{\varphi}_0\,A^\prime(\nabla_c\widetilde{\varphi}_0)\cdot\vec n_w\right]\rd\eta\rd\xi
= 0,\label{eq:wallint3}
\end{align}
by using the identity
\begin{align}
- \partial_\xi\widetilde{\varphi}_0\,\partial_\xi [A^\prime(\nabla_c\widetilde{\varphi}_0)\cdot\vec n_s] -\partial_{\xi}\partial_\eta\widetilde{\varphi}_0\, A^\prime(\nabla_c\widetilde{\varphi}_0)\cdot\vec n_w = \partial_\xi\left[A(\nabla_c\widetilde{\varphi}_0)  - \partial_\xi\widetilde{\varphi}_0 \,A^\prime(\nabla_c\widetilde{\varphi}_0)\cdot\vec n_s\right].\nn
\end{align}

For the first integral in \eqref{eq:wallint3}, applying Gauss's theorem and using the matching condition in \eqref{eq:matchbl} as well as the fact $\lim_{\xi\to+\infty}\partial_{\xi}\widetilde{\varphi}_0=0$, we have
\begin{align}
&\lim_{\xi_1, \eta_1\to+\infty} \int_0^{\eta_1}\bigl[F(\widetilde{\varphi}_0) + A(\nabla_c\widetilde{\varphi}_0) -\partial_{\xi}\widetilde{\varphi}_0\, A^\prime(\nabla_c\widetilde{\varphi}_0)\,\cdot\vec n_s\bigr]_{-\xi_1}^{\xi_1}\deta\nn\\
&\qquad =\int_{0}^{+\infty} F(\widehat{\varphi}_0^+) + A(-\partial_\eta\widehat{\varphi}_0^+\,\vec n_w)\deta - \int_0^{+\infty} F(\widehat{\varphi}_0^-)+A(-\partial_\eta\widehat{\varphi}_0^-\,\vec n_w)\deta\nn\\
&\qquad=2\int_0^{+\infty}F(\widehat{\varphi}_0^+)\deta - 2\int_0^{+\infty}F(\widehat{\varphi}_0^-)\deta\nn\\
&\qquad= 2\,\sqrt{A(-\vec n_w)}\Bigl(\int_0^{+\infty}\sqrt{F(\widehat{\varphi}_0^+)}\partial_\eta\Phi^+_0\deta +  \int_0^{+\infty}\sqrt{F(\widehat{\varphi}_0^-)}\partial_\eta\Phi^-_0\deta\Bigr)=0,
\label{eq:asyca1}
\end{align}
where we have used \eqref{eq:bldf} and \eqref{eq:bldd}.

We then apply Gauss's theorem to the second integral in \eqref{eq:wallint3}. Recalling the boundary condition \eqref{eq:wbd2}, we obtain
\begin{align}
&\int_{-\xi_1}^{\xi_1} \big[\partial_{\xi}\widetilde{\varphi}_0\,A^\prime(\nabla_c\widetilde{\varphi}_0)\cdot\vec n_w\big]\big|_{\eta_1}\rd\xi - \int_{-\xi_1}^{\xi_1}\big[\partial_{\xi}\widetilde{\varphi}_0\,A^\prime(\nabla_c\widetilde{\varphi}_0)\cdot\vec n_w\big]\big|_0\,\rd\xi \nn\\
&\qquad=\int_{-\xi_1}^{\xi_1} \big[\partial_{\xi}\widetilde{\varphi}_0\,A^\prime(\nabla_c\widetilde{\varphi}_0)\cdot\vec n_w\big]\big|_{\eta_1}\rd\xi + \cF\,\sigma\,\int_{-\xi_1}^{\xi_1}\partial_\xi\widetilde{\varphi}_0\,G^\prime(\widetilde{\varphi}_0)\,\rd\xi = I + II.
\end{align}
Sending $\xi_1\to+\infty$ and recalling \eqref{DEF:G}, we obtain
\begin{equation}
\lim_{\xi_1\to+\infty}II=\lim_{\xi_1\to+\infty}\cF\,\sigma\,\int_{-\xi_1}^{\xi_1}\partial_\xi\widetilde{\varphi}_0\,g(\widetilde{\varphi}_0)\,\rd\xi = \cF\,\sigma\,(G(1) - G(-1))=\cF\,\sigma.
\label{eq:asyca2}
\end{equation}

\begin{figure}[!htp]
\centering
\includegraphics[width=0.6\textwidth]{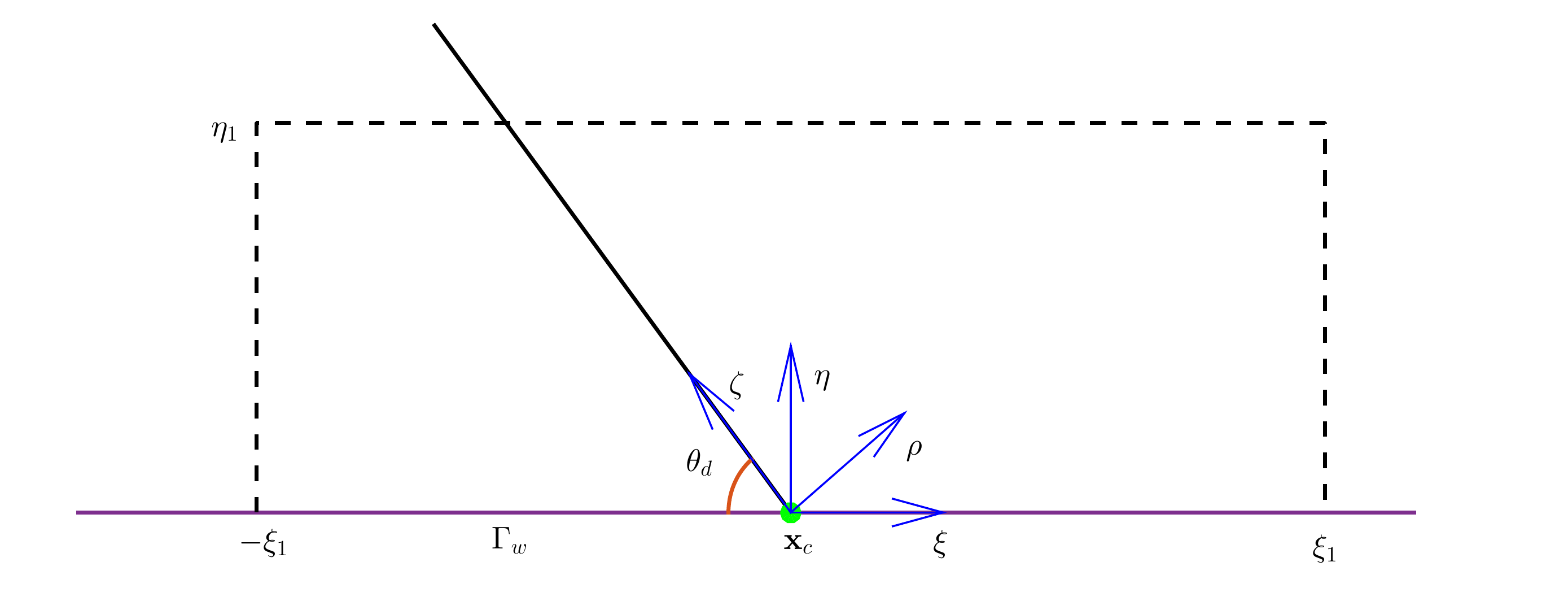}
\caption{Sketch of the local coordinates $(\xi,\eta)$ and $(\rho,\zeta)$ at a contact point $\vec x_c$, where $\theta_d\in(0,\pi)$ is the contact angle.}
\label{fig:tf}
\end{figure}

Next we rewrite the term $I$ in terms of the new coordinate system $(\rho,\zeta)$, which can be regarded as a transformation from $(\eta,\xi)$ with a counterclockwise rotation of $\theta_d$ in the plane (see Fig.~\ref{fig:tf}). Precisely, it holds that
\begin{subequations}\label{eqn:transf}
\begin{align}
\rho = \xi \sin\theta_d + \eta\cos\theta_d,\qquad
\zeta = -\xi\cos\theta_d + \eta\sin\theta_d,
\end{align}
and thus 
\begin{align}
    \partial_\xi = \partial_\rho\sin\theta_d -\partial_\zeta\cos\theta_d,\qquad \partial_\eta=\partial_\rho\cos\theta_d+\partial_\zeta\sin\theta_d.
\end{align}
Moreover, we have
\begin{equation}
   \nabla_c =\vec n_s\,\partial_\xi - \vec n_w\,\partial_\eta =\bs{\nu}\,\partial_\rho -\vec n_c\,\partial_\zeta, \qquad \vec n_c = \sin\theta_d\,\vec n_w + \cos\theta_d\,\vec n_s,
\end{equation}
\end{subequations}
where $\vec n_c$ is the conormal vector of $\Gamma(t)$ at $\vec x_c$. By \eqref{eqn:transf},  we can recast the term $I$ as
\begin{align}
I = \int_{-\xi_1 + \eta_1\cos\theta_d}^{\xi_1\sin\theta_d+\eta_1\cos\theta_d}\bigl[\partial_\rho\widetilde{\varphi}_0 -\partial_\zeta\widetilde{\varphi}_0\cot\theta_d\bigr]\,A^\prime(\partial_\rho\widetilde{\varphi}_0\bs{\nu}- \partial_{\zeta}\widetilde{\varphi}_0\,\vec n_c)\cdot\vec n_w\drho.
\end{align}
By the matching condition $\lim_{\zeta\to+\infty}\widetilde{\varphi}_0 = \Phi_0(\rho)$, we have $\lim_{\zeta\to+\infty}\partial_{\zeta}\widetilde{\varphi}_0=0$. Then it follows directly that
\begin{equation}
\lim_{\xi_1,\eta_1\to+\infty}I=\int_{-\infty}^{+\infty} (\partial_\rho\Phi_0)^2A^\prime(\bs{\nu})\cdot\vec n_w\drho = \frac{\cF\,A^\prime(\bs{\nu})\cdot\vec n_w}{\gamma(\bs{\nu})} = \cF\,\gamma^\prime(\bs{\nu})\cdot\vec n_w. 
\label{eq:asyca3}
\end{equation}

Collecting the results in \eqref{eq:asyca1}, \eqref{eq:asyca2} and \eqref{eq:asyca3} yields that
\begin{equation}
\gamma^\prime(\bs{\nu})\cdot\vec n_w + \sigma=0,
\label{eq:aniY}
\end{equation}
which is exactly the anisotropic Young's law in \eqref{eq:sharpbd2}.

We next derive the zero-flux condition. Similarly to the above, we integrate \eqref{eq:wallld1} over the box $R$. Applying Gauss's theorem and using the boundary condition \eqref{eq:wbd1} gives rise to
\begin{align}
0&=\int_{0}^{\eta_1}\int_{-\xi_1}^{\xi_1}\partial_\xi[\widetilde{\beta}_0\, m(\widetilde{\varphi}_0)\,\partial_\xi\widetilde{\mu}_0] + \partial_\eta[\widetilde{\beta}_0\,m(\widetilde{\varphi}_0)\,\partial_\eta\widetilde{\mu}_0]\,\rd\xi\rd\eta \nn\\
&=\int_0^{\eta_1} \left[\widetilde{\beta}_0\,m(\widetilde{\varphi}_0)\,\partial_\xi\widetilde{\mu}_0\right]\big|_{ -\xi_1}^{\xi_1}\deta + \int_{-\xi_1}^{\xi_1}\left[\widetilde{\beta}_0\,m(\widetilde{\varphi}_0)\,\partial_\eta\widetilde{\mu}_0\right]\big|_{\eta_1}\rd\xi = III+IV.\label{eq:fluxExp}
\end{align}
Taking $\xi_1\to\infty$ and using fact $\lim_{\xi\to\pm\infty}\widetilde{\varphi}_0 = \pm 1$ as well as $m(\widetilde{\varphi}_0)=0$, we get $III=0$. On recalling \eqref{eqn:transf} as well as the matching conditions
\begin{equation*}
  \lim_{\zeta\to+\infty}\widetilde{\varphi}_0 = \Phi_0(\rho)=\tanh\bigl(\frac{\rho}{\gamma(\bs{\nu})}\bigr),\qquad  \lim_{\zeta\to+\infty}\widetilde{\mu}_0 = M_0(\vec s, t)=\varkappa_\gamma,\qquad \lim_{\zeta\to+\infty}\partial_{\zeta}\widetilde{\varphi}_0=0,
\end{equation*}
 we get in the case of $\xi_1, \eta_1\to\infty$ that
\begin{equation}
0=\int_{-\infty}^\infty\left[m(\Phi_0)\,\beta(\bs{\nu})\partial_\zeta\varkappa_\gamma\right]\drho= -\beta(\bs{\nu})\,\int_{-\infty}^{\infty}m(\Phi_0)\drho\;\vec n_c\cdot\nabla_s\varkappa_\gamma.\nn
\end{equation}
This yields the zero-flux condition 
\begin{equation}
D(\bs{\nu})\,\vec n_c\cdot\nabla_s\varkappa_\gamma=0.
\label{eq:aninf}
\end{equation}
In addition, the attachment condition in \eqref{eq:sharpbd1} follows naturally.

In summary, we thus obtain the following system of equations as the sharp-interface limit of the regularized diffuse interface model \eqref{eqn:regmod}:
\begin{subequations}
\label{eqn:silim}
\begin{alignat}{2}
\label{eqn:silim:1}
&\mathcal{V} = -\frac{\cFsqu}{4\,\alpha}\nabla_s\cdot[D(\bs{\nu})\nabla_s\varkappa_\gamma]
\quad\mbox{with}\quad
\mathcal{V} = \vec V \cdot \bs{\nu}
\quad\mbox{and}\quad
\varkappa_\gamma = -\nabla_s\cdot\gamma^\prime(\bs{\nu})
&&\quad\mbox{on}\quad\Gamma(t),
\\
\label{eqn:silim:2}
&\mathcal{\vec V}\cdot\vec n_w = 0,\quad 
\gamma^\prime(\bs{\nu})\cdot\vec n_w + \sigma=0,\quad 
D(\bs{\nu})\,\vec n_c\cdot\nabla_s\varkappa_\gamma=0
&&\quad\mbox{on}\quad\Lambda(t).
\end{alignat}
\end{subequations}

\begin{remark}
In the case of the double-obstacle potential \eqref{eq:doF} and the degenerate mobility $m(\varphi)=(1-\varphi^2)_+$, we could obtain \eqref{eqn:silim:1} in a similar manner. But the leading order inner solution \eqref{eq:phi0} should be replaced by 
\begin{equation}
    \Phi_0(\rho) =\left\{\begin{array}{lll}
    &\sin(\frac{\rho}{\gamma(\bds{\nu})})\quad&\mbox{if}\quad\rho\in[-\frac{\pi}{2}\gamma(\bds{\nu}), \frac{\pi}{2}\gamma(\bds{\nu})],\\[0.4em]
    &-1\quad&\mbox{if}\quad\rho<-\frac{\pi}{2}\gamma(\bds{\nu}),\\[0.4em]
    &1\quad&\mbox{if}\quad\rho>\frac{\pi}{2}\gamma(\bds{\nu}).
   \end{array}\right.
\end{equation}
This yields $\cF=\frac{\pi}{2}$ in \eqref{eqn:silim:1}. The boundary conditions in \eqref{eqn:silim:2} can be derived similarly. It is also possible to consider the logarithmic potential (see \eqref{eq:logF}) along with the mobility $m(\varphi)=(1-\varphi^2)_+$. 
If $\theta = O(\eps^\xi)$ for some $\xi>0$, it can be shown by means of the techniques from
\cite{Cahn1996cahn} that the same desired
sharp interface limit is obtained.
\end{remark}

\section{Analysis of the diffuse interface model}
\label{se:analysis}

In this section, we analyze a general class of diffuse interface models of the type
\begin{subequations}
\label{eq:diffint}
\begin{alignat}{2}
    &\alpha\,\delt \varphi = \Grad \cdot \big(M(\Grad\varphi,\varphi) \Grad\mu \big)
    &&\quad\text{in $Q$},\\
    \label{eq:muexp}
    &\mu = -\eps \Grad\cdot A'(\Grad\varphi) + \eps^{-1} F'(\varphi)
    &&\quad\text{in $Q$},\\
    &\Grad\mu\cdot \n = 0
    &&\quad\text{on $\Sigma$},\\
    &\eps A'(\Grad\varphi)\cdot \n + \cF \,\sigma\, G'(\varphi) = 0
    &&\quad\text{on $\Sigma_w$},\\
    &\eps A'(\Grad\varphi)\cdot \n  = 0
    &&\quad\text{on $\Sigma\setminus\Sigma_w$},\\
    &\varphi\vert_{t=0} = \varphi_0
    &&\quad\text{in $\Omega$},
\end{alignat}
\end{subequations} 
where $\alpha,\eps,\cF \in \R_{>0}$ and $\sigma \in \R$ are given constants.
In contrast to the previous sections, the potential $F:\R\to\R$ as well as $G:\R\to\R$, $A:\R^d\to\R$ and $M:\R^d \times \R\to\R$ are general functions satisfying certain conditions that will be specified in Subsection~\ref{ssec:pre}.
If $A$, $F$, $G$, $m^\eps$ and $\beta^\eps$ are chosen as in \eqref{DEF:A}, \eqref{DEF:F}, \eqref{DEF:G}, \eqref{DEF:ME} and \eqref{DEF:BE}, respectively, and if $M$ is defined by $M(\p,s):= \beta^\eps(\p) m^\eps(s)$ for all $\p\in\R^d$ and $s\in\R$, then the system \eqref{eq:diffint} is exactly the model \eqref{eqn:regmod} that was introduced in Section~\ref{se:asymptotic}. The total free energy functional $\mathcal{E}:H^1(\Omega) \to \R$ 
associated with the system \eqref{eq:diffint},
up to an additive constant, reads as
\begin{align}\label{DEF:E}
\mathcal{E}(\varphi) := \frac{1}{\cF}\int_{\Omega} \varepsilon A(\nabla\varphi) 
    + \varepsilon^{-1} F(\varphi) \dx 
    + \sigma\,\int_{\Gamma_w} G(\varphi) \dS.
\end{align}
It is also possible to consider the system \eqref{eq:diffint} for $F$ being the double-obstacle potential, which can be expressed as
\begin{align}
    \label{DEF:F:OBST}
    F:\R\to [0,\infty], \qquad
    F(s) = F_0(s) + I_{[-1,1]}(s),
\end{align}
where the function
\begin{align}
    \label{DEF:F:OBST:REG}
    F_0:\R\to [0,\infty), \qquad
    F_0(s) = \tfrac 12 (1-\varphi^2),
\end{align}
represents its regular part, and 
\begin{align}
    \label{DEF:F:OBST:IND}
    I:\R\to [0,\infty], \qquad
    I_{[-1,1]}(s) = 
    \begin{cases}
        0 &\text{if $\abs{s}\le 1$},\\
        +\infty &\text{if $\abs{s}> 1$}
    \end{cases}
\end{align}
denotes the indicator functional of the interval $[-1,1]$. In this case, \eqref{eq:muexp} needs to be represented by a variational inequality, see \eqref{DEF:WS*:WF2}.

\subsection{Notation and preliminaries} \label{ssec:pre}
\paragraph{Notation.} In this section, we use the following notation:  For any $1 \leq p \leq \infty$ and $k \geq 0$, the standard Lebesgue and Sobolev spaces on $\Omega$ are denoted by $L^p(\Omega)$ and $W^{k,p}(\Omega)$. Their standard norms are written as $\norm{\,\cdot\,}_{L^p(\Omega)}$ and $\norm{\,\cdot\,}_{W^{k,p}(\Omega)}$.  
In the case $p = 2$, these spaces are Hilbert spaces, and we write $H^k(\Omega) = W^{k,2}(\Omega)$. Here, we identify $H^0(\Omega)$ with $L^2(\Omega)$. 
For the Lebesgue and Sobolev spaces on $\del\Omega$, we use an analogous notation.  For any Banach space $X$, its dual space is denoted by $X'$, and the associated duality pairing by $\ang{\cdot}{\cdot}_X$.  
If $X$ is a Hilbert space, we write $(\cdot, \cdot)_X$ to denote its inner product. We further define
\begin{align*}
\mean{f}_\Omega := 
\frac{1}{\abs{\Omega}} \ang{f}{1}_{H^1(\Omega)} \quad \text{ for } f \in H^1(\Omega)' 
\end{align*}
as the generalized spatial mean of $f$, where $\abs{\Omega}$ denotes the $d$-dimensional Lebesgue measure of $\Omega$.  
With the usual identification $L^1(\Omega) \subset H^1(\Omega)'$ it holds that
$\mean{f}_\Omega = \frac{1}{\abs{\Omega}} \int_\Omega f \dx$ if 
$f \in L^1(\Omega)$.
In addition, we introduce 
\begin{align*}
    H^1_{(m)}(\Omega) &:= \big\{ u \in H^1(\Omega) \,\big\vert\, \mean{u}_\Om = m \big\}
    \qquad\text{for any $m\in\R$},
    \\
    H^{-1}_{(0)}(\Omega) &:= \big\{ f \in \big(H^1(\Omega)\big)' \,\big\vert\, \mean{f}_\Om = 0 \big\}.
\end{align*}
We point out that for every $m\in\R$, $H^1_{(m)}(\Omega)$ is an affine subspace of the Hilbert space $H^1(\Omega)$. In the case $m=0$, it is even a closed linear subspace, meaning that $H^1_{(0)}(\Omega)$ is also a Hilbert space.
\bigskip

\paragraph{General assumptions.} We make the following general assumptions that are supposed to hold throughout this section.
\begin{enumerate}[label=$\mathbf{A \arabic*}$, ref = $\mathbf{A \arabic*}$]
	\item \label{ass:dom} The set $\Omega \subset \R^d$ with $d \in \{2,3\}$ is a bounded Lipschitz domain. 
	Moreover, $T>0$ denotes an arbitrary final time.
	\item  \label{ass:G}
	The function $G:\R\to\R$ is non-negative and twice continuously differentiable. 
	Moreover, there exists an exponent $q\in[2,4)$ as well as 
	positive constants $C_{G}$ and $C_{G'}$ such that
		\begin{align*}
		G(s) \le C_{G}(1+\abs{s}^q) 
		\quad\text{and}\quad
		\abs{G'(s)} \le C_{G'}(1+\abs{s}^{q-1})
		\end{align*}
		for all $s\in\R$.
	\item \label{ass:A} The function $A:\R^d\to \R$ is continuously differentiable and there exist constants $A_0,A_1\in\R$ with $0<A_0\le A_1$ such that 
	\begin{align*}
	    A_0  \abs{\p}^2 \le A(\p) \le A_1  \abs{\p}^2 
	    \quad\text{for all $\p\in\R^d$}.
	\end{align*}
	The gradient $A':\R^d\to\R^d$ is strongly monotone, i.e., there exists a constant $a_0>0$ such that
	\begin{align*}
	    \big(A'(\p)-A'(\q)\big)\cdot(\p-\q) \ge a_0 \abs{\p-\q}^2 
	    \quad\text{for all $\p,\q\in\R^d$},
	\end{align*}
	which implies that $A$ is strongly convex and thus strictly convex. 
	Moreover, there exists a constant $a_1>0$ such that
	\begin{align*}
	    \abs{A'(\p)} \le a_1 \abs{\p}
	    \quad\text{for all $\p\in\R^d$}.
	\end{align*}
    \item \label{ass:M} The function $M:\R^d\times\R\to\R$ is continuous and there exist constants $M_0,M_1\in\R$ with $0<M_0\le M_1$ such that
	\begin{align*}
	    M_0 \le M(\p,s) \le M_1
	    \quad\text{for all $\p\in\R^d$ and $s\in\R$}.
	\end{align*}
\end{enumerate}

\bigskip
\begin{remark}	\label{rem:ass} 
\begin{enumerate}[label = \textnormal{(\alph*)}]
\item \label{rem:ass:a}
We point out that the choices
	\begin{alignat*}{3}
		G(s) &:= \frac{1}{4}(3s - s^3)
		&&\quad\text{for all $s\in\R$} 
		&&\quad\big(\text{cf.\ \eqref{DEF:G}}\big),
		\\
		M(\p,s) &:= \beta^\eps(\p)\, m^\eps(s)
	\end{alignat*} 
	with
	\begin{alignat*}{3}
		m^\eps(s) &:= \eps^r + (1-s^2)_+^2
		&&\quad\text{for all $s\in\R$} 
		&&\quad\big(\text{cf.\ \eqref{DEF:ME}}\big),
		\\
		\beta^\eps(\p) &:= \sqrt{ \frac{d_1^2\,\eps^{2r} + D^2({\p})}{\gamma_0^2\eps^{2r} + \gamma^2(\p)} } 
		&&\quad\text{for all $\p\in\R^d$} 
		&&\quad\big(\text{cf.\ \eqref{DEF:BE}}\big),
	\end{alignat*} 
	are admissible as they satisfy the conditions imposed in \ref{ass:G} (with $q=3$) and \ref{ass:M}.
\item \label{rem:ass:b}
    Suppose that the function $\gamma$ that was introduced in Subsection~\ref{SUBSECT:SI} additionally satisfies the following convexity condition: There exists a constant $\alpha_0>0$ such that
    \begin{align}
    \label{SCON}
    \gamma^{\prime\prime}(\vec p)\vec q\cdot\vec q\geq \alpha_0|\vec q|^2
    \qquad
    \mbox{for all}
\quad \vec p,\vec q\in\R^d\;\mbox{with}\quad|\vec p|=1\quad\mbox{and}\quad \vec p\cdot\vec q = 0,
    \end{align}
    where $\gamma^{\prime\prime}$ represents the Hessian of $\gamma$. 
    Thus, the function 
    \begin{align*}
        A:\R^d\to\R,\quad A(\p) = \tfrac{1}{2} \gamma^2(\p)
    \end{align*}
    is admissible as it satisfies all conditions imposed in assumption \ref{ass:A}.
    In particular, as shown in \cite{Graser2013time}, the convexity condition \eqref{SCON} ensures that $A'$ is strongly monotone.
\end{enumerate}
\end{remark}

\paragraph{A special inner product on $H^{-1}_{(0)}(\Omega)$.}
We now introduce a certain inner product on the function space $H^{-1}_{(0)}(\Omega)$ based on the solution operator of a suitable elliptic problem. 
Therefore, let $a\in L^\infty(\Omega)$ be a uniformly positive function, i.e., there exist $a_0, a_1 \in \R$ with $0<a_0<a_1$ such that 
\begin{align*}
    a_0 \le a \le a_1
    \quad\text{a.e.~in $\Omega$}.
\end{align*}
Then, for every $f\in H^{-1}_{(0)}(\Omega)$, there exists a unique weak solution $u_f\in H^{1}_{(0)}(\Omega)$ of the elliptic boundary value problem
\begin{subequations}
\label{ELL:SF}
\begin{alignat}{2}
    - \Grad\cdot \big( a \Grad  u \big) &= f 
    &&\quad\text{in $\Omega$},
    \\
    \Grad u \cdot \n &= 0
    &&\quad\text{on $\del\Omega$},
\end{alignat}
\end{subequations}
meaning that
\begin{align}
    \label{ELL:WF}
    \intO a \Grad u_f \cdot \Grad \zeta \dx
    = \bigang{f}{\zeta}_{H^1}
    \qquad\text{for all $\zeta\in H^1(\Omega)$.}
\end{align}
We can thus define a solution operator 
\begin{align}
    \label{DEF:SALPHA}
    S_a: H^{-1}_{(0)}(\Omega) \to H^{1}_{(0)}(\Omega), \quad
    f \mapsto S_a(f) := u_f.
\end{align}
We next define the bilinear form
\begin{align}
    \label{DEF:INN:ALPHA}
    \bigscp{\cdot}{\cdot}_{S_a}: H^{-1}_{(0)}(\Omega)\times H^{-1}_{(0)}(\Omega) \to \R, 
    \quad
    \bigscp{f}{g}_{S_a}
    := \intO a \; \Grad S_a(f) \cdot \Grad S_a(g) \dx,
\end{align}
which defines an inner product on $H^{-1}_{(0)}(\Omega)$ since $a$ is uniformly positive and $\Grad S_a(f) = 0$ a.e.~in $\Omega$ already implies $f=0$ a.e.~in $\Omega$. Its induced norm is given by
\begin{align}
    \label{DEF:NORM:ALPHA}
    \norm{\,\cdot\,}_{S_a}: H^{-1}_{(0)}(\Omega) \to \R, 
    \quad
    \norm{f}_{S_a}:= \bigscp{f}{f}_{S_a}^{1/2}.
\end{align}
We point out that on the space $H^{-1}_{(0)}(\Omega)$, the norm $\norm{\,\cdot\,}_{S_a}$ is equivalent to the standard operator norm $\norm{\,\cdot\,}_{(H^1(\Omega))'}\,$.
The bilinear form $\scp{\cdot}{\cdot}_{S_a}$ also defines an inner product on the space $H^{1}_{(0)}(\Omega)$. Moreover, $\norm{\,\cdot\,}_{S_a}$ is also a norm on $H^{1}_{(0)}(\Omega)$ but the space is not complete with respect to this norm.

\subsection{Existence of weak solutions}

For ease of presentation, in what follows we simply fix 
$\alpha=\eps=\sigma=\cF=1$,
since the precise choice of these values has no impact on the 
mathematical analysis. 

\subsubsection{Weak solutions for smooth potentials}
In this subsection, we make the following assumption on the potential $F$:
\begin{enumerate}[label=$\mathbf{F \arabic*}$, ref = $\mathbf{F \arabic*}$]
	\item \label{ass:F1} 
	The potential $F:\R\to\R$ is continuously differentiable. 
	Moreover, there exists an exponent $p\in[2,6)$ as well as 
	non-negative constants $B_F$, $C_{F}$ and $C_{F'}$ such that
		\begin{align*}
		- B_{F} \le F(s) \le C_{F}(1+\abs{s}^p) 
		\quad\text{and}\quad
		\abs{F'(s)} \le C_{F'}(1+\abs{s}^{p-1}).
		\end{align*}
	for all $s\in\R$.
\end{enumerate}
Obviously, the smooth double-well potential introduced in \eqref{eq:dwF} fulfills \ref{ass:F1} with $p=4$. However, the logarithmic potential (see \eqref{eq:logF}) and the double-obstacle potential (see \eqref{eq:doF}) do not satisfy this assumption.

A weak solution of the general diffuse interface model \eqref{eq:diffint} is then defined as follows.
\begin{defn} \label{DEF:WS}
    Suppose that the assumptions \ref{ass:dom}--\ref{ass:M} and \ref{ass:F1} are fulfilled, and let ${\varphi_0\in H^1(\Omega)}$ be any initial datum.
    Then, the pair $(\varphi,\mu)$ is called a weak solution to system \eqref{eq:diffint} if the following properties hold:
    \begin{enumerate}[label=\textnormal{(\roman*)}, ref = \textnormal{(\roman*)}]
        \item \label{DEF:WS:REG} 
        The functions $\varphi$ and $\mu$ have the following regularity:
        \begin{alignat*}{2}
            \varphi &\in C^{0,1/4}\big([0,T];L^2(\Omega)\big) 
                \cap L^\infty\big(0,T;H^1(\Omega)\big)
                \cap H^1\big(0,T;H^1(\Omega)'\big),
            \\
            \mu &\in L^2\big(0,T;H^1(\Omega)\big).
        \end{alignat*}
        \item \label{DEF:WS:WEAK} 
        The pair $(\varphi,\mu)$ satisfies the weak formulations
        \begin{subequations}
        \label{DEF:WS:WF}
        \begin{alignat}{2}
            \label{DEF:WS:WF1}
            \bigang{\delt\varphi}{\zeta}_{H^1(\Omega)} 
            &= - \intO M(\Grad\varphi,\varphi)\, \Grad \mu \cdot \Grad \zeta \dx,
            \\
            \label{DEF:WS:WF2}
            \intO \mu\, \eta \dx 
            &= \intO A'(\Grad\varphi)\cdot \Grad \eta + F'(\varphi)\,\eta \dx
                + \intGw G'(\varphi) \, \eta \dS
        \end{alignat}
        \end{subequations}
        a.e.~on $[0,T]$ for all test functions $\zeta,\eta\in H^1(\Omega)$. Moreover, $\varphi$ satisfies the initial condition 
        \begin{align}
            \label{DEF:WS:INI}
            \varphi(0) = \varphi_0 \qquad\text{a.e.~in $\Omega$.}
        \end{align}
        \item \label{DEF:WS:ENERGY}
        The pair $(\varphi,\mu)$ satisfies the weak energy dissipation law
        \begin{align}
            \label{DEF:WS:DISS}
            \mathcal{E}\big(\varphi(t)\big) 
            + \frac 12 \int_0^t \intO M(\Grad\varphi,\varphi)\, \abs{\Grad\mu}^2 \dx\dt
            \le \mathcal{E}(\varphi_0)
            \qquad\text{for almost all $t\in[0,T]$.}
        \end{align}
    \end{enumerate}
\end{defn}

\bigskip

The existence of such a weak solution is ensured by the following theorem.

\begin{thm} \label{THM:REGPOT}
    Suppose that the assumptions \ref{ass:dom}--\ref{ass:M} and \ref{ass:F1} are fulfilled, and let ${\varphi_0\in H^1(\Omega)}$ be any initial datum. Then there exists a weak solution $(\varphi,\mu)$ to the system \eqref{eq:diffint} in the sense of Definition~\ref{DEF:WS}.
\end{thm}

The proof of this theorem is presented in Section~\ref{SECT:PROOFS}.

\revised{In the next subsection, we intend to prove the existence of a weak solution to the diffuse-interface model \eqref{eq:diffint} for the double-obstacle potential \eqref{eq:doF}. Our strategy is to approximate the double-obstacle potential by a sequence of regular potentials. To this end, in Corollary~\ref{COR:REGPOT}, we will present an additional uniform estimate for $F'(\varphi)$, where $(\varphi,\mu)$ is a weak solution to \eqref{eq:diffint} with a regular potential $F$ satisfying the following assumption:}

\begin{enumerate}[label=$\mathbf{F \arabic*}$, ref = $\mathbf{F \arabic*}$, start=2]
	\item \label{ass:F2} 
    The potential $F:\R\to\R$ is twice continuously differentiable
	and there exist constants $c_0,c_1 \ge 0$ such that 
    \begin{align}
    \label{COND:F''}
        - c_0 \le F''(s) \le c_1 \quad\text{for all $s\in\R$}.
    \end{align}
\end{enumerate}
We point out that if \ref{ass:F2} is fulfilled, then \ref{ass:F1} holds with $p=2$. 

\begin{cor} \label{COR:REGPOT}
    Suppose that the assumptions \ref{ass:dom}--\ref{ass:M}
    and \ref{ass:F2} are fulfilled. 
    Let $\varphi_0\in H^1(\Omega)$ be any initial datum satisfying 
    $\abs{\mean{\varphi_0}_\Omega} \le 1 - \kappa$ for some $\kappa\in(0,1]$,
    and let $(\varphi,\mu)$ be a corresponding weak solution. Then there exists a constant $c>0$ depending only on $\varphi_0$, $\mathcal E(\varphi_0)$, $c_0$ and the constants in \ref{ass:dom}--\ref{ass:M}, but not on $c_1$, such that
    \begin{align}
    \label{EST:F'}
        \norm{F'(\varphi)}_{L^2(Q)}^2 \le \frac{c}{\kappa^2}
        \big( 1 + \norm{F}_{L^\infty([-R,R])}^2 \big),
    \end{align}
    where $ R:= \abs{\mean{\varphi_0}_\Omega} + \frac{\kappa}{2} < 1$. 
\end{cor}

\begin{remark}
    In Corollary~\ref{COR:REGPOT}, the assumption $\abs{\mean{\varphi_0}_\Omega} \le 1 - \kappa$ is made in order to ensure $R\le 1$, which is crucial for later use.
    However, without this assumption a similar estimate can be derived if $R>1$ is allowed. For instance, choosing $R:=\abs{\mean{\varphi_0}_\Omega}+1$, we obtain the estimate
    \begin{align}
    \label{EST:F''}
        \norm{F'(\varphi)}_{L^2(Q)}^2 \le c\big( 1 + \norm{F}_{L^\infty([-R,R])}^2 \big)
    \end{align}
    instead of \eqref{EST:F'} even without the mean value assumption.
\end{remark}

\subsubsection{Weak solutions for the double-obstacle potential}
In this subsection, we assume that $F = F_0 + I_{[-1,1]}$ is the double-obstacle potential as introduced in \eqref{DEF:F:OBST}.
Then a weak solution of the general diffuse interface model \eqref{eq:diffint} is defined as follows.
\begin{defn} \label{DEF:WS*}
    Suppose that the assumptions \ref{ass:dom}--\ref{ass:M} are fulfilled, and let $\varphi_0\in H^1(\Omega)$ be any initial datum satisfying 
$\abs{\varphi_0}\leq 1$ a.e.~in $\Omega$.
    Then, the pair $(\varphi,\mu)$ is called a weak solution to system \eqref{eq:diffint} 
    if the following properties hold:
    \begin{enumerate}[label=\textnormal{(\roman*)}, ref = \textnormal{(\roman*)}]
        \item \label{DEF:WS*:REG} 
        The functions $\varphi$ and $\mu$ have the following regularity:
        \begin{alignat*}{2}
            \varphi &\in C^{0,1/4}\big([0,T];L^2(\Omega)\big) 
                \cap L^\infty\big(0,T;H^1(\Omega)\big)
                \cap H^1\big(0,T;H^1(\Omega)'\big),
            \\
            \mu &\in L^2\big(0,T;H^1(\Omega)\big).
        \end{alignat*}
        \item \label{DEF:WS*:WEAK} 
        It holds that $\abs{\varphi}\le 1$ a.e.~in $Q$ and the pair $(\varphi,\mu)$ satisfies the weak formulation
        \begin{subequations}
        \label{DEF:WS*:WF}
        \begin{align}
            \label{DEF:WS*:WF1}
            \bigang{\delt\varphi}{\zeta}_{H^1(\Omega)} 
            &= - \intO M(\Grad\varphi,\varphi)\, \Grad \mu \cdot \Grad \zeta \dx,
        \end{align}
        for all $\zeta \in H^1(\Omega)$ as well as the variational inequality
        \begin{align}
            \label{DEF:WS*:WF2}
            \iint_Q \mu\, (\varphi-\eta) \dx\dt
            &\ge 
            \iint_Q A'(\Grad\varphi)\cdot(\Grad\varphi-\Grad\eta)
                +F_0'(\varphi)(\varphi-\eta) \dx\dt
            \notag\\
            &\qquad 
            + \intSw G'(\varphi) \, (\varphi-\eta) \dS\dt 
        \end{align}
        \end{subequations}
        for all $\eta\in L^2(0,T;H^1(\Omega))$ with $\abs{\eta}\le 1$ a.e.~in $Q$.
        Moreover, $\varphi$ satisfies the initial condition 
        \begin{align}
            \label{DEF:WS*:INI}
            \varphi(0) = \varphi_0 \qquad\text{a.e.~in $\Omega$.}
        \end{align}
        \item \label{DEF:WS*:ENERGY}
        The pair $(\varphi,\mu)$ satisfies the weak energy dissipation law
        \begin{align}
            \label{DEF:WS*:DISS}
            \mathcal{E}\big(\varphi(t)\big) 
            + \frac 12 \int_0^t \intO M(\Grad\varphi,\varphi)\, \abs{\Grad\mu}^2 \dx\dt
            \le \mathcal{E}(\varphi_0)
            \qquad\text{for almost all $t\in[0,T]$.}
        \end{align}
    \end{enumerate}
\end{defn}

\bigskip

The existence of such a weak solution is ensured by the following theorem.

\begin{thm} \label{THM:OBST}
    Suppose that the assumptions \ref{ass:dom}--\ref{ass:M} are fulfilled, and let $\varphi_0\in H^1(\Omega)$ be any initial datum satisfying 
$\abs{\varphi_0}\leq 1$ a.e.~in $\Omega$ and
$\abs{\mean{\varphi_0}_\Omega} \le 1-\kappa$ for some $\kappa\in (0,1]$.
    Then, there exists a weak solution $(\varphi,\mu)$ to the system \eqref{eq:diffint} in the sense of Definition~\ref{DEF:WS*}.
\end{thm}

\revised{The idea behind the proof of Theorem~\ref{THM:OBST} is to approximate the double-obstacle potential by a sequence $(F_n)_{n\in\N}$ of regular potentials where for each $n\in\N$, $F_n$ is a regular potential fulfilling the condition \ref{ass:F2}. Therefore, Corollary~\ref{COR:REGPOT} can be applied to derive a suitable uniform bound on the terms involving $F_n'$.
We point out that the same strategy could be used to construct a weak solution to the diffuse-interface model \eqref{eq:diffint} in the case that $F$ is the logarithmic potential \eqref{eq:logF}. }

\subsection{Proofs} \label{SECT:PROOFS}

\subsubsection{Proof of Theorem~\ref{THM:REGPOT}}
    The proof is divided into five steps. 
    
    {\bf \textit{Step 1: Implicit time discretization.}}
    Let $N\in\N$ be arbitrary. We define $\tau := T/N$ as our time step size. Let now $n\in\{0,...,N-1\}$ be arbitrary.
    We now define functions $\varphi^n$ with $n=0,...,N$ by the following recursion: 
    \begin{itemize}[leftmargin=*]
    \item The zeroth iterate is defined as the initial datum, i.e., $\varphi^0 := \varphi_0$. 

    \item If for some $n\in\{0,...,N-1\}$ the $n$-th iterate $\varphi^n$ is already constructed, we 
    choose $\varphi^{n+1}\in H^1_{(m)}$ as a minimizer of the functional
    \begin{align}
    \label{DEF:JN}
        J_n: H^1_{(m)}(\Omega) \to \R, \quad
        J_n(\varphi) := \frac{1}{2\tau}\norm{\varphi-\varphi^n}_{S_a}^2 + \mathcal{E}(\varphi).
    \end{align}
    Here, $\mathcal{E}$ is the energy functional defined in \eqref{DEF:E}, with $\eps=\sigma=\cF=1$, and $\norm{\cdot}_{S_a}$ is the norm defined in \eqref{DEF:NORM:ALPHA} with $a$ being chosen as
    \begin{align}
        \label{DEF:ALPHA}
        a := M(\Grad\varphi^n,\varphi^n).
    \end{align}
    This choice is actually possible since the function $M$ is assumed to be bounded and uniformly positive (see \ref{ass:M}). 
    The existence of a minimizer of the functional $J_n$ will be established in Step~2.
    \end{itemize}
    The idea behind this construction is that the first variation of the functional $J_n$ at the point $\varphi^{n+1}$ is zero since $\varphi^{n+1}$ is a minimizer of $J_n$. This means that
    \begin{align}
    \label{EQ:EL}
        \scp{\frac{\varphi^{n+1}-\varphi^n}{\tau}}{\eta}_{S_a} 
            + \intO A'(\Grad\varphi^{n+1})\cdot \Grad \eta + F'(\varphi^{n+1})\, \eta \dx
            + \intGw G'(\varphi^{n+1}) \, \eta \dS
        = 0
    \end{align}
    for all test functions $\eta\in H^1_{(0)}(\Omega)$. We now define
    \begin{align}
        \label{DEF:MUN}
        \mu^{n+1} := S_a\left(- \frac{\varphi^{n+1}-\varphi^n}{\tau}\right) + c^{n+1} \in H^1(\Omega),
    \end{align}
    with
    \begin{align}
        c^{n+1} := \frac{1}{\abs{\Omega}} \left( \intO F'(\varphi^{n+1}) \dx + \intGw G'(\varphi^{n+1}) \dS \right)
    \end{align}
    and $a$ being chosen as in \eqref{DEF:ALPHA}. 
    Recalling the definition of the inner product $\scp{\cdot}{\cdot}_{S_a}$ (see \eqref{DEF:INN:ALPHA}), we infer from \eqref{EQ:EL} that
    \begin{align}
        \label{DISC:2}
        \intO \mu^{n+1} \, \eta \dx
        = \intO A'(\Grad\varphi^{n+1})\cdot \Grad \eta + F'(\varphi^{n+1})\, \eta \dx
            + \intGw G'(\varphi^{n+1}) \, \eta \dS
    \end{align}
    for all $\eta\in H^1_{(0)}(\Omega)$. 
    Due to the choice of the constant $c^{n+1}$, a straightforward computation reveals that \eqref{DISC:2} remains true even for all test functions $\eta\in H^1(\Omega)$. 
    This means that for every $n\in\{0,...,N-1\}$, the pair $(\varphi^{n+1},\mu^{n+1})$ satisfies the equations
    \begin{subequations}
    \label{WF:DISC}
    \begin{align}
        \label{WF:DISC:1}
        \ang{\frac{\varphi^{n+1}-\varphi^n}{\tau}}{\zeta}_{H^1(\Omega)}
        &= - \intO M(\Grad\varphi^n,\varphi^n) \Grad \mu^{n+1} \cdot \Grad \zeta \dx,
        \\[1ex]
        \label{WF:DISC:2}
        \intO \mu^{n+1} \, \eta \dx
        &= \intO A'(\Grad\varphi^{n+1})\cdot \Grad \eta + F'(\varphi^{n+1})\, \eta \dx
            + \intGw G'(\varphi^{n+1}) \, \eta \dS
    \end{align}
    \end{subequations}
    for all test functions $\zeta,\eta \in H^1(\Omega)$.
    Here, \eqref{WF:DISC:1} follows directly from the construction of $\mu^{n+1}$ in \eqref{DEF:MUN} and the definition of the solution operator $S_a$ (see \eqref{DEF:SALPHA}). The system \eqref{WF:DISC} can be interpreted as a time-discrete approximation of the weak formulation \eqref{DEF:WS:WF}.
    
    The time-discrete approximate solution now needs to be extended onto the whole time interval $[0,T]$. The \textit{piecewise constant extension} $(\varphi_N,\mu_N)$ is defined as
    \begin{alignat}{2}
        \label{DEF:EXT:CONST}
        \big(\varphi_N,\mu_N\big)(\cdot,t) :=
        \begin{cases}
        (\varphi_0,\mu_0) 
        &\text{if $t\le 0$}, 
        \\
        (\varphi^{n},\mu^{n}) 
        &\text{if $t\in\big((n-1)\tau,n\tau\big]$, $n\in\{1,...,N\}$}, 
        \end{cases}
    \end{alignat}
    whereas the \textit{piecewise linear extension} $(\ov\varphi_N,\ov\mu_N)$ is defined as
    \begin{align}
        \label{DEF:EXT:LIN}
        (\ov\varphi_N,\ov\mu_N)(\cdot,t)
        := \lambda (\varphi^{n},\mu^{n}) + (1-\lambda) (\varphi^{n-1},\mu^{n-1})
    \end{align}
    for $t = \lambda n\tau + (1-\lambda)(n-1)\tau$ with $n\in\{1,...,N\}$ and $\lambda\in[0,1]$.
    
    Henceforth, the letter $C$ will denote generic positive constants that may depend only on $\varphi_0$ and the constants introduced in \ref{ass:G}--\ref{ass:M} and \ref{ass:F1} but not on $n$, $N$ or $\tau$. These constants may also change their value from line to line.
    
    {\bf \textit{Step 2: Existence of a minimizer to the functional $J_n$.}}
    We now prove that the functional $J_n$ introduced in \eqref{DEF:JN} actually possesses a minimizer. Therefore, we employ the direct method of the calculus of variations.
    
    For any $\varphi \in H^1_{(m)}(\Omega)$, we obtain
    \begin{align*}
        \norm{\varphi}_{L^2(\Omega)} 
        &\le \norm{\mean{\varphi}_\Omega}_{L^2(\Omega)}
            + \norm{\varphi - \mean{\varphi}_\Omega }_{L^2(\Omega)} 
        \le C\big( 1 + \norm{\Grad\varphi}_{L^2(\Omega)}  \big)
    \end{align*}
    by means of Poincar\'e's inequality. This directly implies
    \begin{align}
        \label{EST:H1}
        \norm{\varphi}_{H^1(\Omega)} \le c_* \big( 1 + \norm{\Grad\varphi}_{L^2(\Omega)}  \big)
        \quad\text{for all $\varphi \in H^1_{(m)}(\Omega)$}
    \end{align}
    for some positive constant $c_*$ depending only on $m$ and $\Omega$.
    Recalling the assumptions on $A$ (see \ref{ass:A}), that $F\ge -B_F$ (see \ref{ass:F1}) and that $G\ge 0$ (see \ref{ass:G}), we use Poincar\'e's inequality to derive the estimate
    \begin{align}
        \label{EST:COERC}
        J_n(\varphi) 
        &\ge \intO A(\Grad\varphi) \dx - B_F |\Omega|
        \ge A_0 \norm{\Grad\varphi}_{L^2(\Om)}^2 - B_F |\Omega|
        \notag\\
        &\ge \frac{A_0}{c_*^2} \norm{\varphi}_{H^1(\Om)}^2 - A_0 - B_F |\Omega|
    \end{align}
    for all $\varphi\in H^1_{(m)}(\Omega)$. This means that $J_n$ is coercive and bounded from below. Hence, the infimum
    \begin{align*}
        \mathcal I := \underset{H^1_{(m)}(\Omega)}{\inf}\; J_n
    \end{align*}
    exists, and consequently, there also exists a corresponding minimizing sequence $(\varphi_k)_{k\in\N}$ with
    \begin{align*}
        J_n(\varphi_k) \to \mathcal I \quad\text{as $k\to\infty$}
        \quad\text{and}\quad
        J_n(\varphi_k) \le \mathcal I + 1 \quad\text{for all $k\in\N$}.
    \end{align*}
    Now, \eqref{EST:COERC} directly implies that $(\varphi_k)_{k\in\N}$ is bounded in $H^1_{(m)}(\Omega)$. Using the Banach--Alaoglu theorem, the compact embeddings $H^1_{(m)}(\Omega) \emb L^p(\Omega)$ and $H^1_{(m)}(\Omega) \emb L^q(\partial\Omega)$, we infer that there exists a function $\ov\varphi\in H^1_{(m)}(\Omega)$ such that
    \begin{align}
    \label{CONV:MIN:PHI}
        \begin{aligned}
        \varphi_k \to \ov\varphi \quad 
            &\text{weakly in $H^1_{(m)}(\Omega)$, 
            strongly in $L^p(\Omega)$ and in $L^q(\partial\Omega)$,}
            \\
            &\quad\text{pointwise a.e.~in $\Omega$, and pointwise a.e.~on $\del\Omega$}
        \end{aligned}
    \end{align}
    along a non-relabeled subsequence. Since $A$ is continuous and convex (see \ref{ass:A}), we infer
    \begin{align}
        \label{CONV:MIN:A}
        \intO A(\Grad \ov\varphi) \dx
        \le \underset{k\to\infty}{\lim\inf}\; \intO A(\Grad \varphi_k) \dx
    \end{align}
    due to weak lower semicontinuity.
    Recalling the growth conditions on $F$ and $G$ (see \ref{ass:F1} and \ref{ass:G}) and the convergences in \eqref{CONV:MIN:PHI}, we apply Lebesgue's general convergence theorem (see \cite[Section~3.25]{Alt}) to conclude 
    \begin{align}
        \label{CONV:MIN:FG}
        \intO F(\varphi_k) \dx \to \intO F(\ov\varphi) \dx
        \quad\text{and}\quad
        \intGw G(\varphi_k) \dS \to \intGw G(\ov\varphi) \dS
    \end{align}
    as $k\to\infty$. Combining \eqref{CONV:MIN:A} and \eqref{CONV:MIN:FG}, we obtain
    \begin{align*}
        J_n(\ov\varphi) \le \underset{k\to\infty}{\lim\inf} J_n(\varphi_k) = \mathcal I.
    \end{align*}
    This proves that $\ov\varphi$ is a minimizer of the functional $J_n$.
    
    {\bf \textit{Step 3: A priori estimates for the piecewise constant extension.}}
    We now claim that the piecewise constant extension $(\varphi_N,\mu_N)$ fulfills the uniform  priori estimate
    \begin{align}
        \label{EST:AP}
        \norm{\varphi_N}_{L^\infty(0,T;H^1(\Omega))} + \norm{\mu_N}_{L^2(0,T;H^1(\Omega))} \le C.
    \end{align}
    
    To prove \eqref{EST:AP}, we exploit the recursive construction of the time-discrete approximate solution. Since for any $n\in\{0,...,N-1\}$, $\varphi^{n+1}$ was chosen to be a minimizer of the functional $J_n$, we have
    \begin{align}
        \label{EST:AP:1}
         \frac{1}{2\tau}\norm{\varphi^{n+1}-\varphi^n}_{S_a}^2 + \mathcal{E}(\varphi^{n+1})
         = J_n(\varphi^{n+1}) 
         \le J_n(\varphi^{n}) 
         = \mathcal{E}(\varphi^{n})
    \end{align}
    for all $n\in\{0,...,N-1\}$. By a simple induction, we thus infer
    \begin{align}
        \label{EST:AP:2}
         \mathcal{E}(\varphi^n)
         \le \mathcal{E}(\varphi_0)
         \quad\text{for all $n\in\{0,...,N-1\}$.}
    \end{align}
    Recalling the assumptions on $A$ (see \ref{ass:A}) and that the potentials $F$ and $G$ are bounded from below (see \ref{ass:G} and \ref{ass:F1}), we use estimate \eqref{EST:H1} and \eqref{EST:AP:2} to obtain
    \begin{align}
        \label{EST:AP:3}
         \norm{\varphi^{n+1}}_{H^1(\Omega)}^2
         \le C + C \norm{\Grad\varphi^{n+1}}_{L^2(\Omega)}^2
         \le C + C \intO A(\Grad\varphi^{n+1})
         \le C + C\mathcal{E}(\varphi^{n+1})
         \le C
    \end{align}
    for all $n\in\{0,...,N-1\}$. By the definition of $\varphi_N$, this directly implies
    \begin{align}
        \label{EST:AP:PHIN}
        \norm{\varphi_N}_{L^\infty(0,T;H^1(\Omega))} \le C.
    \end{align}
    For any $n\in\{1,...,N\}$, we now set $t_n := n\tau$.
    By the definition of the piecewise constant extension, we have
    \begin{align}
        \label{PROP:EXT:CONST}
        \varphi_N(t) = \varphi(t_n) = \varphi^{n}
        \quad\text{and}\quad
        \mu_N(t) = \mu(t_n) = \mu^{n}
    \end{align}
    for all $t\in (t_{n-1},t_n]$.
    Recalling the priori estimate \eqref{EST:AP:1} and the definition of $\mu^{n}$ (see \eqref{DEF:MUN}), we obtain
    \begin{align*}
        &\mathcal{E}\big(\varphi_N(t_n)\big) 
            + \frac 12 \int_{t_{n-1}}^{t_n}\intO M\big(\Grad\varphi_N(s-\tau),\varphi_N(s-\tau)\big) \, \abs{\Grad\mu_N(s)}^2 \dx\ds
        \\
        &\quad = \mathcal{E}\big(\varphi_N(t_n)\big) 
            + \frac{1}{2\tau^2} \int_{t_{n-1}}^{t_n} 
            \norm{\varphi_N(s) - \varphi_N(s-\tau)}_{S_a}^2 \ds
        \\
        &\quad = \mathcal{E}\big(\varphi_N(t_n)\big) 
            + \frac{1}{2\tau} 
            \norm{\varphi_N(t_n) - \varphi_N(t_n-\tau)}_{S_a}^2 
        \\[1ex]
        &\quad\le \mathcal{E}\big(\varphi_N(t_{n-1})\big)
    \end{align*}
    for all $n\in\{0,...,N-1\}$. Hence, by induction, we get
    \begin{align*}
        \mathcal{E}\big(\varphi_N(t_n)\big) 
            + \frac 12 \int_{0}^{t_n} \intO M\big(\Grad\varphi_N(s-\tau),\varphi_N(s-\tau)\big) \, \abs{\Grad\mu_N(s)}^2 \dx \ds
        \le \mathcal{E}(\varphi_0)
    \end{align*}
    for all $n\in\{0,...,N-1\}$. 
    Now, for any $t\in (0,T]$ we find an index $n\in\{0,...,N-1\}$ such that $t\in (t_{n-1},t_n]$.
    Recalling \eqref{PROP:EXT:CONST}, we eventually conclude that
    \begin{align}
    \label{ENERGY:DISC}
        \mathcal{E}\big(\varphi_N(t)\big) 
            + \frac 12 \int_{0}^{t} \intO M\big(\Grad\varphi_N(s-\tau),\varphi_N(s-\tau)\big) \, \abs{\Grad\mu_N(s)}^2 \dx \ds
        \le \mathcal{E}(\varphi_0)
    \end{align}
    for all $t\in[0,T]$. In particular, choosing $t=T$, we obtain the uniform bound
    \begin{align}
        \label{EST:AP:GRADMUN}
        \norm{\Grad\mu_N}_{L^2(0,T;L^2(\Omega))}^2 \le C.
    \end{align}
    We now test \eqref{WF:DISC:2} with the constant function $\eta\equiv 1/\abs{\Omega}$. Using the growth assumptions from \ref{ass:F1}, the continuous embeddings $H^1(\Omega) \emb L^5(\Omega)$ and $H^1(\Omega) \emb L^3(\partial\Omega)$ as well as the uniform bound \eqref{EST:AP:PHIN}, we derive the estimate
    \begin{align*}
        \abs{\mean{\mu_N(t)}_\Omega} 
        &\le \frac{1}{\abs{\Omega}} \left( \intO \abs{F'\big(\varphi_N(t)\big)} \dx 
            + \intGw \abs{G'\big(\varphi_N(t)\big)} \dS \right)
        \\[1ex]
        &\le C \big( 1 + \norm{\varphi_N(t)}_{L^5(\Omega)}^5 
            + \norm{\varphi_N(t)}_{L^3(\partial\Omega)}^3 \big)
        \\[1ex]
        &\le C \big( 1 + \norm{\varphi_N}_{L^\infty(0,T;H^1(\Omega))}^5 
            + \norm{\varphi_N}_{L^\infty(0,T;H^1(\Omega))}^3\big)
        \le C.
    \end{align*}
    Applying Poincar\'e's inequality, we thus obtain
    \begin{align}
        \label{EST:AP:MUNT}
        \norm{\mu_N(t)}_{L^2(\Omega)} 
        &\le \norm{\mean{\mu_N(t)}_\Omega}_{L^2(\Omega)}
            + \norm{\mu_N - \mean{\mu_N(t)}_\Omega}_{L^2(\Omega)} 
        \notag\\
        &\le C\big( 1 + \norm{\Grad \mu_N(t)}_{L^2(\Omega)} \big).
    \end{align}
    Combining \eqref{EST:AP:GRADMUN} and \eqref{EST:AP:MUNT}, this yields
    \begin{align}
        \label{EST:AP:MUN}
        \norm{\mu_N}_{L^2(0,T;H^1(\Omega))} \le C.
    \end{align}
    Due to \eqref{EST:AP:PHIN} and \eqref{EST:AP:MUN}, the a priori estimate \eqref{EST:AP} is now established.
    
    {\bf \textit{Step 4: A priori estimate for the piecewise linear extension.}}
    We next claim that for all $s,t\in[0,T]$,%
    \begin{subequations}
    \label{EST:LIN}
    \begin{align}
        \label{EST:LIN:1}
        \norm{\ov\varphi_N(t) - \ov\varphi_N(s)}_{L^2(\Omega)} 
        &\le C\abs{t-s}^{\frac 14},
        \\[0.75ex]
        \label{EST:LIN:2}
        \norm{\ov\varphi_N(t) - \varphi_N(t)}_{L^2(\Omega)} 
        &\le C\tau^{\frac 14},
        \\[1ex]
        \label{EST:LIN:3}
        \norm{\delt \ov\varphi_N}_{L^2(0,T;H^1(\Omega)')} 
        &\le C.
    \end{align}
    \end{subequations}
    In particular, the first estimate means that
    the piecewise linear extension $\ov\varphi_N$ is Hölder continuous in time. 
    
    To prove these inequalities, we first infer from \eqref{WF:DISC:1} and the definition of the piecewise linear extension (see \eqref{DEF:EXT:LIN}) that
    \begin{align}
        \label{WF:LIN}
        \ang{\delt\ov\varphi_N(\tau)}{\zeta}_{H^1(\Omega)}
        &= - \intO \Grad \mu_N(\tau) \cdot \Grad \zeta \dx
    \end{align}
    for almost all $\tau\in[0,T]$ and all $\zeta\in H^1(\Omega)$. 
    Let now $\xi \in L^2(0,T;H^1(\Omega))$ be arbitrary.
    We test \eqref{WF:LIN} with $\xi(\tau)$
    and integrate the resulting equation with respect to $\tau$ from $0$ to $T$. Then, using Hölder's inequality as well as the a priori estimate \eqref{EST:AP}, we obtain
    \begin{align}
        \label{EST:LIN:1*}
        \abs{ \int_0^T \ang{\delt\ov\varphi_N(\tau)}{\xi}_{H^1(\Omega)} \dt }
        \le \norm{\mu_N}_{L^2(0,T;H^1(\Omega))}\, 
            \norm{\xi}_{L^2(0,T;H^1(\Omega))}
        \le C\, \norm{\xi}_{L^2(0,T;H^1(\Omega))}.
    \end{align}
    Taking the supremum over all $\xi \in L^2(0,T;H^1(\Omega))$ with $\norm{\xi}_{L^2(0,T;H^1(\Omega))} \le 1$, this proves estimate \eqref{EST:LIN:3}.
    
    Next, let $s,t\in [0,T]$ be arbitrary. Without loss of generality, we assume $s<t$. Integrating \eqref{WF:LIN} with respect to $\tau$ from $s$ to $t$, choosing $\zeta = \ov\varphi_N(t) - \ov\varphi_N(s)$, and using Hölder's inequality, we derive the estimate
    \begin{align}
        \norm{\ov\varphi_N(t) - \ov\varphi_N(s)}_{L^2(\Omega)}^2 
        &\le \norm{\Grad\ov\varphi_N(t) - \Grad\ov\varphi_N(s)}_{L^2(\Omega)}
            \int_s^t \norm{\Grad\mu_N(\tau)}_{L^2(\Omega)} \dtau
        \notag\\
        &\le 2 \, \norm{\varphi_N}_{L^\infty(0,T;H^1(\Omega))} \,
            \norm{\mu_N}_{L^2(0,T;L^2(\Omega))} \abs{s-t}^{\frac 12}.
    \end{align}
    In view of the a priori estimate \eqref{EST:AP}, this proves \eqref{EST:LIN:1}.
    
    Let now $t\in[0,T]$ be arbitrary. Then, we find $\lambda\in [0,1]$ and $n\in\{1,...,N\}$ such that $t=\lambda n \tau + (1-\lambda) (n-1) \tau$. We thus obtain 
    \begin{align*}
        \norm{\ov\varphi_N(t) - \varphi_N(t)}_{L^2(\Omega)}
        &= \norm{\lambda \varphi^n + (1-\lambda) \varphi^{n-1} 
            - \varphi^n}_{L^2(\Omega)}
        \notag\\
        &= (1-\lambda) \norm{\varphi^n - \varphi^{n-1}}_{L^2(\Omega)}
        \notag\\
        &= (1-\lambda) \norm{\varphi_N \big( n\tau \big) -         
            \varphi_N \big( (n-1)\tau \big)}_{L^2(\Omega)}.
    \end{align*}
    Applying \eqref{EST:LIN:1} with $t=n\tau$ and $s=(n-1)\tau$, we conclude \eqref{EST:LIN:2}. This means that all estimates in  \eqref{EST:LIN} are established.
    
    {\bf \textit{Step 5: Convergence to a weak solution.}}
    In view of the uniform a priori estimate \eqref{EST:AP}, the Banach--Alaoglu theorem implies the existence of functions $\varphi\in L^\infty(0,T;H^1(\Omega))$ and $\mu\in L^2(0,T;H^1(\Omega))$
    such that
    \begin{alignat}{2}
        \label{CONV:1}
        \varphi_N &\to \varphi
        &&\quad\text{weakly-$^*$ in $L^\infty(0,T;H^1(\Omega))$,}
        \\
        \label{CONV:2}
        \mu_N &\to \mu
        &&\quad\text{weakly in $L^2(0,T;H^1(\Omega))$},
    \end{alignat}
    as $N\to\infty$, along a non-relabeled subsequence. 
    We further know that 
    \begin{align*}
        \norm{\ov\varphi_N}_{L^\infty(0,T;H^1(\Omega))}
        \le \norm{\varphi_N}_{L^\infty(0,T;H^1(\Omega))}
        \le C.
    \end{align*}
    In combination with the uniform estimate \eqref{EST:LIN:3}, we use the Banach--Alaoglu theorem to infer $\varphi\in H^1(0,T;H^1(\Omega)')$ with
    \begin{align}
        \label{CONV:3}
        \ov\varphi_N \to \varphi 
        \quad \text{weakly in $H^1(0,T;H^1(\Omega)')$}
    \end{align}
    as $N\to\infty$, up to subsequence extraction. Moreover, due to the compact embeddings $H^1(\Omega) \emb L^p(\Omega)$ and $H^1(\Omega) \emb L^q(\partial\Omega)$, we apply the Aubin--Lions lemma to obtain
    \begin{align}
        \label{CONV:4}
        \ov\varphi_N \to \varphi 
        \quad \text{strongly in $C([0,T];L^p(\Omega))
            \cap C([0,T];L^q(\partial\Omega))$}.
    \end{align}
    By passing to the limit in estimate \eqref{EST:LIN:1}, we conclude $\varphi\in C^{0,1/4}([0,T],L^2(\Omega))$. This means that the functions $\varphi$ and $\mu$ satisfy the regularity conditions of Definition~\ref{DEF:WS}\ref{DEF:WS:REG}.
    Using the estimate \eqref{EST:LIN:2}, we directly deduce from \eqref{CONV:4} that
    \begin{align}
        \label{CONV:5}
        \begin{aligned}
        \varphi_N \to \varphi
        &\quad\text{strongly in $L^\infty(0,T;L^p(\Omega)) \cap L^\infty(0,T;L^q(\partial\Omega))$},
        \\
        &\qquad\text{a.e.~in $\Omega$, and a.e.~on $\partial\Omega$},
        \end{aligned}
    \end{align}
    as $N\to\infty$, after another subsequence extraction.
    
    From the time-discrete weak formulation \eqref{WF:DISC}, we infer that the piecewise constant extension $(\varphi_N,\mu_N)$ and the piecewise linear extension $(\ov\varphi_N,\ov\mu_N)$ satisfy the approximate weak formulation
    \begin{subequations}
        \label{WF:APR}
    \begin{align}
        \label{WF:APR:1}
        &\int_0^T \ang{\delt\ov\varphi_N(t)}{\xi}_{H^1(\Omega)} \dt
        = - \intQ M\big(\Grad\varphi_N(t-\tau),\varphi_N(t-\tau)\big)\,
            \Grad \mu_N(t) \cdot \Grad \xi \dx\dt,
        \\[1ex]
        \label{WF:APR:2}
        &\intQ \mu_N \, \vartheta \dx\dt
        = \intQ A'(\Grad\varphi_N)\cdot \Grad \vartheta 
            + F'(\varphi_N)\, \vartheta \dx\dt
        + \intSw G'(\varphi_N) \, \vartheta  \dS\dt
    \end{align}
    \end{subequations}
    for all test functions $\xi,\vartheta \in L^2(0,T;H^1(\Omega))$.
    Recalling the growth conditions on $F'$ and $G'$ from 
\ref{ass:F1} and \ref{ass:G} as well as the priori estimate \eqref{EST:AP}, we infer that the sequence $(F'(\varphi_N))_{N\in\N}$ is bounded in $L^\infty(0,T;L^{6/5}(\Omega))$ and the sequence $(G'(\varphi_N))_{N\in\N}$ is bounded in $L^\infty(0,T;L^{4/3}(\partial\Omega))$. Hence, according to the Banach--Alaoglu theorem, there exist functions $f^* \in L^\infty(0,T;L^{6/5}(\Omega))$ and $g^* \in L^\infty(0,T;L^{4/3}(\partial\Omega))$ such that
    \begin{alignat*}{2}
        F'(\varphi_N) &\to f^*
        &&\quad\text{weakly-$^*$ in $L^\infty(0,T;L^{6/5}(\Omega))$,}
        \\
        G'(\varphi_N) &\to g^*
        &&\quad\text{weakly-$^*$ in $L^\infty(0,T;L^{4/3}(\partial\Omega))$,}
    \end{alignat*}
    as $N\to\infty$, along a non-relabeled subsequence. Moreover, the convergences in \eqref{CONV:5} directly imply $F'(\varphi_N) \to F'(\varphi)$ a.e.~in $\Omega$ and $G'(\varphi_N) \to G'(\varphi)$ a.e.~on $\partial\Omega$. By a convergence principle based on Egorov's theorem (see \cite[Proposition~9.2c]{DiBenedetto}), we now infer $f^* = F'(\varphi)$ a.e.~in $\Omega$ and $g^* = G'(\varphi)$ a.e.~on $\partial\Omega$.
    This means that 
    \begin{alignat}{2}
        \label{CONV:f}
        F'(\varphi_N) &\to F'(\varphi)
        &&\quad\text{weakly-$^*$ in $L^\infty(0,T;L^{6/5}(\Omega))$,}
        \\
        \label{CONV:g}
        G'(\varphi_N) &\to G'(\varphi)
        &&\quad\text{weakly-$^*$ in $L^\infty(0,T;L^{4/3}(\partial\Omega))$,}
    \end{alignat}
    as $N\to\infty$. 
    Testing the approximate weak formulation \eqref{WF:APR:2} with $\vartheta = \varphi_N - \varphi$ and employing the strong monotonicity condition on $A'$ from \ref{ass:A}, we obtain
    \begin{align}
        &a_0 \norm{\Grad\varphi_N - \Grad\varphi}_{L^2(Q)}^2
        \le \intQ \big( A'(\Grad\varphi_N) - A'(\Grad\varphi) \big)
            \cdot \big( \Grad\varphi_N - \Grad\varphi \big) \dx\dt
        \notag\\
        &\quad= \intQ \mu_N\, (\varphi_N - \varphi) \dx\dt
            - \intQ F'(\varphi_N)\, (\varphi_N - \varphi) \dx\dt
        \notag\\
        &\qquad - \intSw G'(\varphi_N)\, (\varphi_N - \varphi) \dS\dt
            - \intQ A'(\Grad\varphi) \cdot \big( \Grad\varphi_N - \Grad\varphi \big) \dx\dt.
    \end{align}
    Using the convergences \eqref{CONV:2}, \eqref{CONV:5}, \eqref{CONV:f} and \eqref{CONV:g} along with the weak-strong convergence principle, we infer that the right-hand side of the above estimate tends to zero. We thus conclude that
    \begin{align}
        \label{CONV:GPHI}
        \Grad\varphi_N \to \Grad\varphi
        \quad\text{strongly in $L^2(Q)$ and a.e.~in $\Omega$}
    \end{align}
    as $N\to\infty$, up to subsequence extraction. 
    In view of the growth condition on $A'$ from \ref{ass:A}, Lebesgue's general convergence theorem further reveals that
    \begin{align}
        \label{CONV:A}
        A'(\Grad\varphi_N) \to A'(\Grad\varphi)
        \quad\text{strongly in $L^2(Q;\R^d)$}.
    \end{align}
    Due to the convergences \eqref{CONV:2}, \eqref{CONV:f}, \eqref{CONV:g} and \eqref{CONV:A}, we can now pass to the limit in \eqref{WF:APR:2} to conclude that
    \begin{align}
        \label{WF:2}
        \intQ \mu \, \vartheta \dx\dt
        = \intQ A'(\Grad\varphi)\cdot \Grad \vartheta 
            + F'(\varphi)\, \vartheta \dx\dt
        + \intSw G'(\varphi) \, \vartheta  \dS\dt
    \end{align}
    holds for all $\vartheta\in L^2(0,T;H^1(\Omega))$. 
    
    We now fix an arbitrary time $t_0\in (0,T]$. Since $\tau = T/N \to 0$ as $N\to\infty$, we may assume (without loss of generality) that $N$ is chosen large enough to ensure $t-\tau \in [0,T]$ for all $t\in [t_0,T]$.
    We have
    \begin{align}
    \label{EST:TAUN}
        &\norm{\Grad\varphi_N(t - \tau) - \Grad\varphi(t) }_{L^2(t_0,T;L^2(\Omega))}^2
        \notag\\[1ex]
        &\quad\le 
            C \norm{\Grad\varphi_N(t - \tau) 
                - \Grad\varphi(t - \tau) }_{L^2(t_0,T;L^2(\Omega))}^2
            + C \norm{\Grad\varphi(t - \tau)  
                - \Grad\varphi(t) }_{L^2(t_0,T;L^2(\Omega))}^2
        \notag\\[1ex]
        &\quad\le 
            C \norm{\Grad\varphi_N(t) 
                - \Grad\varphi(t) }_{L^2(0,T;L^2(\Omega))}^2
            + C \norm{\Grad\varphi(t - \tau)  
                - \Grad\varphi(t) }_{L^2(t_0,T;L^2(\Omega))}^2
    \end{align}
    for almost all $t\in[t_0,T]$.
    Here, from the second to the third line, we used the change of variables $s = t-\tau$ and the fact that $[t_0-\tau,T-\tau]\subset [0,T]$ to estimate the first summand. Now, as $N\to\infty$, the first summand in the third line of \eqref{EST:TAUN} tends to zero because of \eqref{CONV:GPHI}, whereas the second summand tends to zero since due to mean-continuity in $L^p(Q)$ (see, e.g., \cite[Section~4.15]{Alt}). This proves
    \begin{align}
        \label{CONV:GPHI:2*}
        \Grad\varphi_N(\,\cdot\,,\,\cdot - \tau) \to \Grad\varphi
        \quad\text{strongly in $L^2(\Omega\times[t_0,T])$}
    \end{align}
    as $N\to\infty$. Since $t_0\in(0,T]$ was arbitrary, we deduce
    \begin{align}
        \label{CONV:GPHI:2}
        \Grad\varphi_N(\,\cdot\,,\,\cdot - \tau) \to \Grad\varphi
        \quad\text{a.e. in $Q$}
    \end{align}
    as $N\to\infty$, after extracting a subsequence.
    Proceeding similarly, and using the strong convergence $\varphi_N\to\varphi$ in $L^2(Q)$ (which directly follows from \eqref{CONV:5}), we further obtain
    \begin{align}
        \label{CONV:PHI:2}
        \varphi_N(\,\cdot\,,\,\cdot - \tau) \to \varphi
        \quad\text{a.e.~in $Q$}
    \end{align}
    as $N\to\infty$.
    Using \eqref{CONV:GPHI:2} and \eqref{CONV:PHI:2} along with Lebesgue's dominated convergence theorem, we infer 
    \begin{align}
    \label{CONV:MBETA}
        M\big(\Grad\varphi_N(\,\cdot\,,\,\cdot - \tau),\varphi_N(\,\cdot\,,\,\cdot - \tau)\big) \,
        \Grad\zeta
        \;\;\to\;\; 
        M(\Grad\varphi,\varphi) \,
        \Grad\zeta
    \end{align}
    strongly in $L^2(Q)$, as $N\to\infty$, up to subsequence extraction. Employing the weak-strong convergence principle, we can thus pass to the limit $N\to\infty$ in the approximate weak formulation \eqref{WF:APR:1} to obtain
    \begin{align}
        \label{WF:1}
        \int_0^T \ang{\delt\varphi}{\zeta}_{H^1(\Omega)} \dt
        = - \intQ M(\Grad\varphi,\varphi)\,
            \Grad \mu \cdot \Grad \zeta \dx\dt
    \end{align}
    for all $\zeta\in L^2(0,T;H^1(\Omega))$. Combining \eqref{WF:2} and \eqref{WF:1}, we eventually conclude that the pair $(\varphi,\mu)$ satisfies the weak formulation \eqref{DEF:WS:WF}. Moreover, as a direct consequence of the convergence \eqref{CONV:4}, $\varphi$ satisfies the initial condition \eqref{DEF:WS:INI}. This means that all conditions of Definition~\ref{DEF:WS}\ref{DEF:WS:WEAK} are fulfilled.
    
    Recalling the growth conditions on $F$ and $G$ from \ref{ass:F1} and \ref{ass:G} as well as the convergences in \eqref{CONV:5}, we apply 
    Lebesgue's general convergence theorem (see \cite[Section~3.25]{Alt}) to conclude 
    \begin{alignat}{2}
        \label{CONV:F}
        F(\varphi_N) &\to F(\varphi)
        &&\quad\text{strongly in $L^{1}(Q)$,}
        \\
        \label{CONV:G}
        G(\varphi_N) &\to G(\varphi)
        &&\quad\text{strongly in $L^{1}(Q)$.}
    \end{alignat}
    Then, from the convergences \eqref{CONV:GPHI},\eqref{CONV:F} and \eqref{CONV:G}, we infer that
    \begin{align}
        \label{CONV:E}
        \mathcal{E}\big(\varphi_N(t)\big) \to \mathcal{E}\big(\varphi(t)\big) 
        \quad\text{for almost all $t\in[0,T]$,}
    \end{align}
    as $N\to\infty$.  
    Recalling \eqref{CONV:2} and \eqref{CONV:MBETA}, we use the weak-strong convergence principle to infer
    \begin{align}
        \label{CONV:DISS}
        &\sqrt{M\big(\Grad\varphi_N(\,\cdot\,,\,\cdot - \tau),\varphi_N(\,\cdot\,,\,\cdot - \tau)\big)}\; \Grad\mu_N
        \to 
        \sqrt{M(\Grad\varphi,\varphi)}\; \Grad\mu
        &&\text{weakly in $L^2(Q)$}
    \end{align}
    as $N\to\infty$.
    We now use the convergences \eqref{CONV:E} and \eqref{CONV:DISS}, the weak lower semicontinuity of the $L^2(Q)$-norm as well as the discrete energy inequality \eqref{ENERGY:DISC} to derive the estimate
    \begin{align}
    \label{IEQ:E}
        &\mathcal{E}\big(\varphi(t)\big) + \frac 12 \int_0^t\intO M\big(\Grad\varphi(s),\varphi(s)\big) \abs{\Grad\mu(s)}^2 \dx\ds
        \notag\\[1ex]
        &\le \underset{N\to\infty}{\lim\inf}\; \mathcal{E}\big(\varphi_N(t)\big) 
        + \underset{N\to\infty}{\lim\inf}\; \frac 12 \int_0^t\intO M\big(\Grad\varphi_N(s-\tau),\varphi_N(s-\tau)\big) \abs{\Grad\mu_N(s)}^2 \dx\ds
        \notag\\[1ex]
        &\le \underset{N\to\infty}{\lim\inf} \Bigg[\mathcal{E}\big(\varphi_N(t)\big)
        + \frac 12 \int_0^t\intO M\big(\Grad\varphi_N(s-\tau),\varphi_N(s-\tau)\big) \abs{\Grad\mu_N(s)}^2 \dx\ds \Bigg]
        \notag\\[1ex]
        &\le \mathcal{E}(\varphi_0)
    \end{align}
    for almost all $t\in[0,T]$.
    This proves the weak energy dissipation law \eqref{DEF:WS:DISS} and thus, the condition in Definition~\ref{DEF:WS}\ref{DEF:WS:ENERGY} is fulfilled.
    
    We eventually conclude that the pair $(\varphi,\mu)$ is a weak solution to system \eqref{eq:diffint} in the sense of Definition~\ref{DEF:WS}. Hence, the proof is complete. \hfill$\Box$

\subsubsection{Proof of Corollary~\ref{COR:REGPOT}}

Let $(\varphi,\mu)$ be a weak solution to the system \eqref{eq:diffint},
whose existence is guaranteed by Theorem~\ref{THM:REGPOT}.
    By a straightforward computation, we notice that
    \begin{align}
    \label{EST:COR:0}
        \norm{F'(\varphi)}_{L^2(Q)}^2
        = \int_0^T \intO (F'(\varphi))^2 \dx\dt 
        \le 2 I_1 + \frac{2}{\abs{\Omega}} I_2,
    \end{align}
    where
    \begin{align*}
        I_1 := \int_0^T \intO \big(F'(\varphi) - \mean{F'(\varphi)}_\Omega \big)^2 \dx\dt
        \quad\text{and}\quad
        I_2 := \int_0^T \left( \intO \abs{F'(\varphi)} \dx \right)^2 \dt.
    \end{align*}
    Hence, in the following, we intend to prove \eqref{EST:F'} by deriving suitable bounds on the terms $I_1$ and $I_2$. The letter $C$ will denote generic positive constants depending only on $\varphi_0$, $\mathcal E(\varphi_0)$, $c_0$ and the constants in \ref{ass:dom}--\ref{ass:M}, but not on $c_1$.
    
    Let $\eta\in H^1(\Omega)$ be arbitrary. Since $A$ is convex (see \ref{ass:A}), we know that
    \begin{align*}
        A'(\Grad \varphi)\cdot \Grad(\eta-\varphi) \le A(\Grad \eta) - A(\Grad \varphi)
        \quad\text{a.e. in $Q$.}
    \end{align*}
    Testing the weak formulation \eqref{DEF:WS:WF2} with $\eta-\varphi$ instead of $\eta$ and using the above estimate, we thus infer that the variational inequality
    \begin{align}
        \label{VAR}
        \intO F'(\varphi) (\varphi-\eta) \dx
        &\le \intO \mu (\varphi - \eta) \dx 
            + \intO A(\Grad\eta) - A(\Grad\varphi) \dx
        - \intGw G'(\varphi) (\varphi-\eta) \dS  
    \end{align}
    holds a.e.~in $[0,T]$ for all $\eta\in H^1(\Omega)$.
    Moreover, since $(\varphi,\mu)$ is a weak solution of \eqref{eq:diffint}, it satisfies the weak energy inequality \eqref{DEF:WS:DISS}. Using Poincar\'e's inequality, we infer
    \begin{align}
        \label{EST:WS}
        \norm{\varphi}_{L^\infty(0,T;H^1(\Omega))}
            + \norm{\mu}_{L^2(0,T;H^1(\Omega))}
            \le C.
    \end{align} 
    
    {\bf \textit{Step 1:}} We first derive an estimate for the term $I_1$. Therefore, we choose
    \begin{align}
        \eta := \varphi - \delta \big(F'(\varphi) - \mean{F'(\varphi)}_\Omega \big)
    \end{align}
    for sufficiently small $\delta>0$ which ensures $1 - \delta F''(\varphi) > 0$.
    Since $F'(\varphi) \in L^\infty(0,T;H^1(\Omega))$ due to \eqref{COND:F''},
    we know that $\eta \in L^\infty(0,T;H^1(\Omega))$.
    Recalling that $A$ is positively homogeneous of degree $2$, we obtain
    \begin{align}
        \label{ID:A}
        &A\big(\Grad\eta\big) - A(\Grad\varphi)
        = A\big(\Grad \varphi - \delta F''(\varphi) \Grad\varphi \big) - A(\Grad\varphi)
        \notag\\
        &\quad = \big( 1-\delta F''(\varphi) \big)^2 A(\Grad\varphi) - A(\Grad\varphi)
        = \big( -2F''(\varphi) + \delta^2 F''(\varphi)^2 \big) A(\Grad\varphi)
    \end{align}
    a.e.~in $Q$.
    We now test the variational inequality \eqref{VAR} with $\eta$.
    After dividing the resulting inequality by $\delta$, we use \eqref{ID:A} to deduce
    \begin{align*}
        &\intO \big(F'(\varphi) - \mean{F'(\varphi)}_\Omega \big)^2 \dx
        = \intO F'(\varphi) \big(F'(\varphi) - \mean{F'(\varphi)}_\Omega \big) \dx
        \notag\\
        &\quad \le \intO \big( \mu - \mean{\mu}_\Omega \big) 
            \big(F'(\varphi) - \mean{F'(\varphi)}_\Omega \big) \dx
        - \intGw G'(\varphi) \big(F'(\varphi) - \mean{F'(\varphi)}_\Omega \big) \dS 
        \notag\\
        &\quad\qquad + \intO \big( -2F''(\varphi) + \delta^2 F''(\varphi)^2 \big) A(\Grad\varphi) \dx
    \end{align*}
    a.e.~in $[0,T]$. Recalling that \ref{ass:F2} implies that \ref{ass:F1} holds with $p=2$, we derive the estimate
    \begin{align*}
         \abs{\mean{F'(\varphi)}_\Omega} 
         \le C + C\intO \abs{\varphi} \dx
         \le C + C\norm{\varphi}_{L^\infty(0,T;L^1(\Omega))} \le C
    \end{align*}
    a.e.~in $[0,T]$.
    Hence, using the growth condition on $G'$ from \ref{ass:G} and the continuous embedding $H^1(\Omega)\emb L^4(\partial\Omega)$, we deduce
    \begin{align}
        &\abs{\intGw G'(\varphi) \big(F'(\varphi) - \mean{F'(\varphi)}_\Omega \big) \dS}
        \le \intGw \abs{G'(\varphi)}\big(\abs{F'(\varphi)} + C\big) \dS
        \notag\\
        &\quad \le \intGw \big(C+C\abs{\varphi}^3\big)\big(C+C\abs{\varphi}\big) \dS
        \le C + C \intGw \abs{\varphi}^4 \dS
        \notag\\[1ex] 
        &\quad \le C + C\norm{\varphi}_{L^4(\partial\Omega)}^4 
        \le C + C\norm{\varphi}_{H^1(\Omega)}^4
    \end{align}
    a.e.~in $[0,T]$. Sending $\delta\to 0$ and using the growth condition from \ref{ass:A}, the condition $-F'' \le c_0$ (cf.~\eqref{COND:F''}) as well as Poincar\'e's inequality and Young's inequality, we infer
    \begin{align*}
        \intO \big(F'(\varphi) - \mean{F'(\varphi)}_\Omega \big)^2 \dx
        \le C \big( 1 + \norm{\varphi}_{H^1(\Omega)}^4 + \norm{\mu}_{H^1(\Omega)}^2 \big) 
    \end{align*}
    a.e.~in $[0,T]$. 
    Integrating this inequality with respect to time from $0$ to $T$,
    and using estimate \eqref{EST:WS}, we eventually conclude the bound
    \begin{align}
        \label{EST:I1}
        I_1 \le C.
    \end{align}
    
   {\bf  \textit{Step 2:}} We now derive a suitable estimate for the term $I_2$. Let $\lambda \in L^\infty([0,T])$ be any function that will be fixed later. We set
    \begin{align}
        \eta := \varphi - \delta (\varphi - \mean{\varphi}_\Omega )
    \end{align}
    for some $\delta>0$. Testing the variational formulation with this $\eta$, dividing the resulting equation by $\delta$, and recalling that $A$ is positively homogeneous of degree $2$, we derive the estimate
    \begin{align}
        \label{EST:COR:1}
        &\intO F'(\varphi) (\lambda - \mean{\varphi}_\Omega ) \dx
        \notag\\
        &\quad = \intO F'(\varphi) (\lambda - \varphi ) \dx
            + \intO F'(\varphi) (\varphi - \mean{\varphi}_\Omega ) \dx
        \notag\\
        &\quad \le \intO F'(\varphi) (\lambda - \varphi ) \dx
            + \intO (\mu - \mean{\mu}_\Omega) \varphi \dx 
            - \intGw G'(\varphi) (\varphi - \mean{\varphi}_\Omega ) \dS
        \notag\\
        &\quad\qquad + \delta(\delta-2) \intO A(\Grad\varphi) \dx.
    \end{align}
    Since $F'' +c_0 \ge 0$ due to \eqref{COND:F''}, we know that the function $s\mapsto F(s) + \tfrac 12 c_0s^2 $ is convex. We thus have 
    \begin{align*}
        F(\lambda) + \tfrac 12 c_0 \lambda^2
        &\ge F(\varphi) + \tfrac 12 c_0 \varphi^2 
            + \big(F'(\varphi) + c_0\varphi \big)(\lambda - \varphi)
        \notag\\
        &\ge \big(F'(\varphi) + c_0\varphi \big)(\lambda - \varphi)
    \end{align*}
    a.e.~in $Q$. Using this estimate as well as Young's inequality, we now get
    \begin{align}
        \label{EST:COR:2}
        \intO F'(\varphi) (\lambda - \varphi ) \dx
        \le \intO F(\lambda) \dx + \intO \tfrac 32 c_0 \varphi^2 + c_0\lambda^2 \dx
    \end{align}
    almost everywhere in $[0,T]$.
    Sending $\delta\to 0$ in \eqref{EST:COR:1} and using the above estimate, the growth conditions from \ref{ass:G} and \ref{ass:A}, the continuous embedding $H^1(\Omega)\emb L^4(\partial\Omega)$ as well as Poincar\'e's inequality, we infer
    \begin{align}
        \label{EST:COR:3}
        \intO F'(\varphi) (\lambda - \mean{\varphi}_\Omega ) \dx
        \le C \norm{F(\lambda)}_{L^\infty([0,T])}
            + C \big( 1 + \norm{\varphi}_{H^1(\Omega)}^4
            + \norm{\mu}_{H^1(\Omega)}\norm{\varphi}_{H^1(\Omega)} \big).
    \end{align}
    We now fix $\lambda$ as
    \begin{align}
        \lambda(t) :=
        \begin{cases}
            \mean{\varphi(t)}_\Omega + \tfrac{\kappa}{2}
            &\text{if $\mean{F'\big(\varphi(t)\big)}_\Omega \ge 0$,}
            \\
            \mean{\varphi(t)}_\Omega - \tfrac{\kappa}{2}
            &\text{if $\mean{F'\big(\varphi(t)\big)}_\Omega < 0$.}
        \end{cases}
    \end{align}
    for all $t\in[0,T]$.
    Testing \eqref{DEF:WS:WF1} with $\zeta\equiv 1$ and integrating the resulting equation with respect to time, we infer $\mean{\varphi(t)}_\Omega = \mean{\varphi_0}_\Omega$ for all $t\in[0,T]$.
    In view of \eqref{EST:COR:3}, we thus get
    \begin{align}
        \label{EST:COR:4}
        \frac{\kappa}{2} \intO \abs{F'(\varphi)} \dx
        \le C \norm{F(\lambda)}_{L^\infty([0,T])}
            + C \big( 1 + \norm{\varphi}_{H^1(\Omega)}^4
            + \norm{\mu}_{H^1(\Omega)}\norm{\varphi}_{H^1(\Omega)} \big).
    \end{align}
    a.e.~in $[0,T]$.
    We now multiply this estimate by $\tfrac 2\kappa$ and take the square on both sides. Integrating the resulting inequality with respect to time and using the uniform estimate \eqref{EST:WS}, we eventually conclude
    \begin{align}
        \label{EST:I2}
        I_2 
        &\le C\kappa^{-2} \norm{F}_{L^\infty([-R,R])}^2
            + C\kappa^{-2}\big( 1 + \norm{\varphi}_{H^1(\Omega)}^8
            + \norm{\mu}_{H^1(\Omega)}^2\norm{\varphi}_{H^1(\Omega)}^2 \big) 
        \notag\\
        &\le C\kappa^{-2} \norm{F}_{L^\infty([-R,R])}^2 + C.
    \end{align}
    
    We finally plug the estimates \eqref{EST:I1} for $I_1$ and \eqref{EST:I2} for $I_2$ into \eqref{EST:COR:0}. This proves \eqref{EST:F'} and thus, the proof of Corollary~\ref{COR:REGPOT} is complete. \hfill$\Box$

\subsubsection{Proof of Theorem~\ref{THM:OBST}}
The proof is split into three steps.

{\bf \textit{ Step 1: Approximation of the double-obstacle potential by smooth potentials.}}
To prove the assertion, we approximate the double-obstacle potential $F$ by a sequence of regular potentials $(F_n)_{n\in\N}$. Therefore, we define the function 
\begin{align*}
    J:\R\to [0,\infty), \quad
    &s\mapsto 
    \begin{cases}
        6s^2 + 20s + 17 &\text{if $s\le -2$},\\
        (s+1)^4 &\text{if $-2<s<-1$},\\
        0 &\text{if $-1\le s \le 1$},\\
        (s-1)^4 &\text{if $1<s<2$},\\
        6s^2 - 20s + 17 &\text{if $s\ge 2$},\\
    \end{cases}
\end{align*}
and for any $n\in\N$, we set
\begin{align}
    \label{DEF:FN}
    F_n:\R\to [0,\infty), \quad
    s\mapsto F_0(s) + n J(s).
\end{align}
By this construction, we have $J\in C^2(\R;[0,\infty))$, $J$ is convex, and $F_n = F_0$ on $[-1,1]$ for all $n\in\N$.
It is straightforward to check that for all $n\in\N$, the approximate potential $F_n$ satisfies the assumption \ref{ass:F2} with 
$c_0 = 1$ and $c_1 = 12n$. It thus follows that \ref{ass:F1} is satisfied with $p=2$ and $B_F = \tfrac 32$.
In the remainder of this proof it will be crucial that the constants $B_F$ and $c_0$ are independent of $n$. For any $n\in\N$ we further define the approximate energy functional by 
    \begin{align}
    \label{DEF:En}
        \mathcal E_n:H^1(\Omega) \to \R, \quad
        \mathcal E_n(\varphi) := \int_{\Omega} A(\nabla\varphi) +  F_n(\varphi) \dx
            + \int_{\Gamma_w} G(\varphi) \dS.
\end{align}

    {\bf \textit{Step 2: A priori estimates for the sequence of approximate solutions.} }
    We now conclude from Theorem~\ref{THM:REGPOT} that for every $n\in\N$, there exists a weak solution $(\varphi_n,\mu_n)$ of the system \eqref{eq:diffint} to the potential $F_n$ in the sense of Definition~\ref{DEF:WS}. 
    In the following, the letter $C$ will denote generic positive constants that may depend on $\varphi_0$, $\kappa$ and the constants in \ref{ass:dom}--\ref{ass:M} but not on the approximation index $n$.

    As the weak solutions $(\varphi_n,\mu_n)$ satisfy the weak energy dissipation law \eqref{DEF:WS:DISS} written for $\mathcal E_n$, we deduce the estimate
    \begin{align*}
        &\frac 12 \intO \abs{\Grad\varphi_n(t)}^2 \dx - B_F |\Omega|
        + \frac 12 M_0 \int_0^t  \intO \abs{\Grad\mu_n(s)}^2 \dx\ds 
        \notag\\
        &\quad
        \le \mathcal E_n\big(\varphi_n(t)\big) 
            + \frac 12  \int_0^t \intO M\big(\Grad\varphi_n(s),\varphi_n(s)\big) \abs{\Grad\mu_n(s)}^2 \dx\ds 
        \notag\\
        &\quad 
        \le \mathcal E_n(\varphi_0)
        \le C \norm{\Grad\varphi_0}_{L^2(\Omega)}^2
            + C \norm{F_0}_{L^\infty([-1,1])}
        \le C
    \end{align*}
    for almost all $t\in[0,T]$ and all $n\in\N$.
    As $B_F$ is independent of $n$, we use Poincar\'e's inequality to conclude the uniform bound
    \begin{align}
        \label{EST:WS*:1}
        \norm{\varphi_n}_{L^\infty(0,T;H^1(\Omega))}
            + \norm{\mu_n}_{L^2(0,T;H^1(\Omega))}
            \le C.
    \end{align} 
    Integrating the weak formulation \eqref{DEF:WS:WF1} written for $(\varphi_n,\mu_n)$ with respect to time from $0$ to $T$, we now use \eqref{EST:WS*:1} to derive the uniform estimate
    \begin{align}
        \label{EST:WS*:2}
        \norm{\delt\varphi_n}_{L^2(0,T;H^1(\Omega)')}
            \le C.
    \end{align}
    Furthermore, Corollary~\ref{COR:REGPOT} provides the estimate
    \begin{align}
    \label{EST:F':N}
        \norm{F_n'(\varphi)}_{L^2(Q)}^2 
        \le \frac{c}{\kappa^2}
        \big( 1 + \norm{F_n}_{L^\infty([-R,R])}^2 \big),
    \end{align}
    where $ R= \abs{\mean{\varphi_0}_\Omega} + \frac{\kappa}{2} < 1$. Here the constant $c$ depends only on $\varphi_0$, $\mathcal E_n(\varphi_0)$, $c_0=1$ and the constants in \ref{ass:dom}--\ref{ass:M}.
    Since $F_n = F_0$ on $[-1,1]$, we know that $F_n(\varphi_0) = F_0(\varphi_0)$ for all $n\in\N$.
    Consequently, $\mathcal E_n(\varphi_0)$
    does not depend on $n$ and thus, $c$ is independent of $n$.
    We infer the uniform bound
    \begin{align}
    \label{EST:F'*}
        \norm{F_n'(\varphi_n)}_{L^2(Q)}^2 
        \le \frac{c}{\kappa^2} \big( 1 + \norm{F_0}_{L^\infty([-1,1])}^2 \big)
        \le C.
    \end{align}
    Using \eqref{EST:WS*:1}, we further get
    \begin{align}
    \label{EST:I'*}
        \norm{F_0'(\varphi_n)}_{L^2(Q)} 
        = \norm{\varphi_n}_{L^2(Q)} 
        \le C \norm{\varphi_n}_{L^\infty(0,T;L^2(\Omega))}
        \le C.
    \end{align}
    Combining \eqref{EST:F'*} and \eqref{EST:I'*}, we now conclude
    \begin{align}
    \label{EST:J'*}
        \norm{J'(\varphi_n)}_{L^2(Q)} 
        \le \frac 1n \big( \norm{F_0'(\varphi_n)}_{L^2(Q)} 
            + \norm{F_n'(\varphi_n)}_{L^2(Q)} \big)
        \le \frac Cn.
    \end{align}
    
    {\bf \textit{Step 3: Convergence to a weak solution.}}
    In view of the uniform estimates \eqref{EST:WS*:1} and \eqref{EST:WS*:2}, we now use the continuous embedding $H^1(\Omega)\emb L^4(\partial\Omega)$, the Banach--Alaoglu theorem, and the Aubin--Lions lemma along with the compact embeddings $H^1(\Omega)\emb L^2(\Omega)$  and $H^1(\Omega)\emb L^r(\partial\Omega)$ for $r\in[1,4)$ to conclude the existence of functions $\varphi$ and $\mu$ such that
    \begin{alignat}{2}
        \label{CONV:1*}
        \delt \varphi_n &\to \delt \varphi
        &&\quad\text{weakly in $L^2(0,T;H^1(\Omega)')$,} 
        \\
        \label{CONV:2*}
        \varphi_n &\to \varphi
        &&\quad\text{weakly-$^*$ in $L^\infty(0,T;H^1(\Omega))$ and in $L^\infty(0,T;L^4(\partial\Omega))$,}
        \notag\\
        &&&\qquad\text{strongly in $C([0,T];L^2(\Omega))$, a.e.~in $Q$,}
        \notag\\
        &&&\qquad\text{strongly in $C([0,T];L^r(\partial\Omega))$
        and a.e.~on $\Sigma$,}
        \\
        \label{CONV:3*}
        \mu_n &\to \mu
        &&\quad\text{weakly in $L^2(0,T;H^1(\Omega))$},
    \end{alignat}
    for all $r\in[1,4)$ as $n\to\infty$, after extraction of a subsequence. 
    Using the uniform bound \eqref{EST:J'*} along with the Banach--Alaoglu theorem as well as the pointwise--a.e.~convergence stated in \eqref{CONV:2*}, we deduce 
    \begin{subequations}
    \label{CONV:4*}
    \begin{alignat}{2}
        \label{CONV:4.1*}
        J'(\varphi_n) &\to 0
        &&\quad\text{strongly in $L^2(Q)$},
        \\
        \label{CONV:4.2*}
        J'(\varphi_n) &\to J'(\varphi)
        &&\quad\text{a.e.~in $Q$},
    \end{alignat}
    \end{subequations}
    as $n\to\infty$.
    As the strong limit in $L^2(Q)$ and the pointwise limit coincide, we have $J'(\varphi) = 0$ a.e.~in $Q$. Since $J'(s)=0$ if $\abs{s}\le 1$ and $\abs{J'(s)}>0$ if $\abs{s}>1$, we conclude 
    \begin{align*}
        \abs{\varphi}\le 1
        \qquad\text{a.e.~in $Q$.}
    \end{align*}
    As $F_0'(\varphi_n) = - \varphi_n$, the convergence
    \begin{alignat}{2}
        \label{CONV:5*}
        F_0'(\varphi_n) &\to F_0'(\varphi)
        &&\quad\text{weakly in $L^2(Q)$ and a.e.~in $Q$},
    \end{alignat}
    follows directly from \eqref{CONV:2*}.
    Moreover, using the growth condition on $G'$ (see \ref{ass:G}), \eqref{CONV:2*} and Lebesgue's general convergence theorem (see \cite[Section~3.25]{Alt}), we obtain
    \begin{alignat}{2}
        \label{CONV:6*}
        G'(\varphi_n) &\to G'(\varphi)
        &&\quad\text{strongly in $L^{4/3}(\Sigma)$ and a.e.~on $\Sigma$},
    \end{alignat}
    as $n\to\infty$, after another subsequence extraction.
    Arguing as in the proof of Theorem~\ref{THM:REGPOT}, we exploit the strong monotonicity condition on $A'$ from \ref{ass:A} to further derive the convergences
    \begin{alignat}{2}
        \label{CONV:7*}
        \Grad\varphi_n &\to \Grad\varphi
        &&\quad\text{strongly in $L^2(Q)$ and a.e.~in $Q$},
        \\
        \label{CONV:8*}
        A'(\Grad\varphi_n) &\to A'(\Grad\varphi)
        &&\quad\text{strongly in $L^2(Q;\R^d)$},
    \end{alignat}
    as $n\to\infty$, up to subsequence extraction. Combining \eqref{CONV:2*} and \eqref{CONV:7*}, we eventually get
    \begin{alignat}{2}
        \label{CONV:9*}
        M(\Grad\varphi_n,\varphi_n) &\to M(\Grad\varphi,\varphi)
        &&\quad\text{a.e.~in $Q$}.
    \end{alignat}
    
    Let now $n\in\N$ be arbitrary and let $\zeta\in H^1(\Omega)$ and $\eta\in L^2(0,T;H^1(\Omega))$ with $\abs{\eta}\le 1$ a.e.~in $Q$ be an arbitrary test functions.
    This already implies that $J'(\eta) = 0$ a.e.~in $Q$. Moreover, since $J$ is convex its derivative $J'$ is monotonically increasing. We thus have
    \begin{align}
        \label{EST:J*}
        J'(\varphi_n) (\varphi_n - \eta) 
        \ge J'(\eta) (\varphi_n - \eta)
        = 0
        \quad\text{a.e.~in $Q$}.
    \end{align}
    We now recall that the weak solution $(\varphi_n,\mu_n)$ satisfies the weak formulation \eqref{DEF:WS:WF} written for $(\varphi_n,\mu_n)$ instead of $(\varphi,\mu)$. 
    The weak formulation \eqref{DEF:WS:WF1} written for $(\varphi_n,\mu_n)$ and tested with $\zeta$ reads as
    \begin{alignat}{2}
        \label{WS*:1}
        \bigang{\delt\varphi_n}{\zeta}_{H^1(\Omega)} 
        &= - \intO M(\Grad\varphi_n,\varphi_n)\, \Grad \mu \cdot \Grad \zeta \dx.
    \end{alignat}
    Testing the weak formulation \eqref{DEF:WS:WF2} written for $(\varphi_n,\mu_n)$ with $\varphi_n-\eta$, integrating with respect to time from $0$ to $T$, and employing estimate \eqref{EST:J*}, we obtain
    \begin{align}
        \label{WS*:2}
        \iint_Q \mu_n\, (\varphi_n - \eta) \dx\dt
        &\ge \iint_Q A'(\Grad\varphi_n)\cdot (\Grad\varphi_n - \Grad\eta) + F_0'(\varphi_n)\, (\varphi_n - \eta) \dx\dt
        \notag\\
        &\qquad 
        + \intSw G'(\varphi_n) \, (\varphi_n - \eta) \dS\dt.
    \end{align}
    Using the convergences \eqref{CONV:1*}--\eqref{CONV:9*}, Lebesgue's dominated convergence theorem as well as the weak-strong convergence principle, we pass to the limit $n\to\infty$ in \eqref{WS*:1} and \eqref{WS*:2}. This proves that the pair $(\varphi,\mu)$ satisfies the weak formulation \eqref{DEF:WS*:WF1} for all $\zeta\in H^1(\Omega)$ as well as the variational inequality \eqref{DEF:WS*:WF2} for all $\eta\in L^2(0,T;H^1(\Omega))$ with $\abs{\eta}\le 1$ a.e.~in $Q$. Moreover, \eqref{CONV:2*} directly implies that $\varphi$ satisfies the initial condition \eqref{DEF:WS*:INI}. This means that all conditions of Definition~\ref{DEF:WS*}\ref{DEF:WS*:WEAK} are verified.
    
    Proceeding similarly as in Step 4 of the proof of Theorem~\ref{THM:REGPOT}, and using the weak formulation \eqref{DEF:WS*:WF1}, we can show a posteriori that $\varphi$ is Hölder continuous in time in the sense that $\varphi \in C^{0,1/4}([0,T];L^2(\Omega))$. In combination with \eqref{CONV:1*}--\eqref{CONV:3*}, this proves that all conditions of Definition~\ref{DEF:WS*}\ref{DEF:WS*:REG} are fulfilled.
    
    Recalling that $\abs{\varphi}\le 1$ a.e.~in $Q$, we use \eqref{CONV:2*} along with Lebesgue's general convergence theorem (see \cite[Section~3.25]{Alt}) to conclude
    \begin{align*}
        \intO F(\varphi) \dx
        = \intO F_0(\varphi) \dx
        = \underset{n\to\infty}{\lim}\; \intO F_0(\varphi_n) \dx
        \le \underset{n\to\infty}{\lim\inf}\; \intO F_n(\varphi_n) \dx
    \end{align*}
    a.e. in $[0,T]$. 
    Using the convergences \eqref{CONV:2*}, \eqref{CONV:3*}, \eqref{CONV:7*} and \eqref{CONV:9*}, we now proceed similarly as in Step~5 of the proof of Theorem~\ref{THM:REGPOT} (cf.~\eqref{IEQ:E}) to verify that the pair $(\varphi,\mu)$ satisfies the weak energy dissipation law \eqref{DEF:WS*:DISS}. This means that Definition~\ref{DEF:WS*}\ref{DEF:WS*:ENERGY} is also fulfilled.
    
    In summary, we conclude that the pair $(\varphi,\mu)$ is a weak solution to system \eqref{eq:diffint} (with $F$ being the 
    double-obstacle potential) in the sense of Definition~\ref{DEF:WS*}. Hence, the proof of Theorem~\ref{THM:OBST} is complete. \hfill$\Box$

\section{Numerical results}
\label{se:nr}
In this section, we present numerical comparisons between the 
diffuse-interface model \eqref{eqn:regmod} and its sharp-interface limit \eqref{eqn:silim}.

For the sharp-interface computations, ({\bf SI}), 
we employ the parametric finite element approximation from \cite{tjtrue},
which uses piecewise linear finite elements and
relies crucially on the stable approximation of the anisotropy
introduced in \cite{triplejANI,ani3d}, see also \cite{BaoZ20preprint}. 
Here, we recall that this stable approximation is designed for anisotropy
functions of the form
\begin{equation} \label{eq:bgn}
\gamma(\vec p)=\sum_{\ell=1}^L \sqrt{\Lambda_{\ell} \vec p \cdot \vec p}, 
\end{equation}
where $\Lambda_\ell$, $\ell=1,\ldots,L$ 
are symmetric and positive definite matrices. We refer to 
\cite{triplejANI,ani3d,Barrett10cluster,BaoZ20preprint,tjtrue} for details.
\revised{Clearly, for \eqref{eq:bgn} the assumption \ref{ass:A} is satisfied, recall
Remark~\ref{rem:ass}.}

For the diffuse-interface approximations, ({\bf DI}), we adapt the 
finite element discretizations from \cite{vch} to the system
\eqref{eqn:regmod}.
To this end, we assume that $\Omega$ is a polyhedral domain and
let $\mathcal{T}_{h}$ be a regular triangulation of $\Omega$ into disjoint open simplices. 
Associated with $\mathcal{T}_h$ is the piecewise linear finite element space
\begin{align*}
S^{h} = \left \{\zeta \in C^{0}(\overline\Omega) : \,  \zeta_{\vert_{o}} \in P_{1}(o) \, \forall o \in \mathcal{T}_{h} \right \},
\end{align*}
where we denote by $P_{1}(o)$ the set of all affine linear functions on $o$,
cf.\ \cite{Ciarlet78}. 
We also let $(\cdot,\cdot)$ denote the $L^{2}$-inner product on $\Omega$,
and let $(\cdot,\cdot)^{h}$ be the usual mass lumped $L^{2}$-inner product on
$\Omega$ associated with $\mathcal{T}_{h}$.
In a similar fashion, we let $\langle \cdot,\cdot \rangle_{\Gamma_w}^h$
denote the mass lumped $L^{2}$-inner product on $\Gamma_w$.
Finally, $\Delta t$ denotes a chosen uniform time step size.

Our fully discrete finite element approximation of \eqref{eqn:regmod} is then
given as follows. For $n\geq 0$, let $\varphi_h^{n} \in S^h$ be given.
Then find $(\varphi_h^{n+1}, \mu^{n+1}_h) \in S^h \times S^h$ such that
\begin{subequations} \label{eq:FEA}
\begin{align} 
\alpha 
\left(\frac{\varphi_h^{n+1} - \varphi_h^n}{\Delta t}, \chi\right)^h
+ \eps^{-1}\bigl(m^\eps(\varphi_h^n) \beta^\eps(\nabla \varphi_h^n) \nabla \mu_h^{n+1},
\nabla \chi\bigr) &= 0, \\
\eps \bigl(B(\nabla \varphi_h^{n})\nabla \varphi_h^{n+1}, \nabla \eta\bigr)
+ \eps^{-1} \bigl(F'(\varphi_h^{n+1}), \eta\bigr)^h
+ \cF\sigma\langle G'(\varphi_h^{n+1}), \eta \rangle_{\Gamma_w}^h\label{eq:muhimpl}
& = ( \mu_h^{n+1}, \eta)^h 
\end{align}
\end{subequations}
for all $(\chi, \eta) \in S^h \times S^h$. The above scheme utilizes the
linearization $B(\vec p) \vec p = A'(\vec p)$ for anisotropies of the form
\eqref{eq:bgn}, which was first introduced in \cite{BGN2013stable}.  
In particular, the symmetric positive definite matrices $B$ are defined by
\begin{equation} \label{eq:BGNB}
B(\vec p) =
\begin{cases}
\gamma(\vec p) \displaystyle\sum_{\ell =1}^L 
\frac{\Lambda_\ell}{\sqrt{\Lambda_{\ell} \vec p \cdot \vec p}}
& \vec p \neq \vec 0, \\
L \displaystyle\sum_{\ell =1}^L \Lambda_\ell & \vec p =\vec 0.
\end{cases}
\end{equation}
\revised{%
We stress that the induced semi-implicit discretization of $A'(\nabla\varphi)$ in 
\eqref{eq:muhimpl} ensures that our numerical method is stable. In fact, using the techniques in
\cite{BGN2013stable,vch}, and on employing semi-implicit approximations of $F'(\varphi)$
and $G'(\varphi)$ based on convex/concave splittings of $F$ and $G$, an unconditional stability
result can be shown. However, for the purposes of this paper we prefer the simpler
approximation \eqref{eq:FEA}. We also note that extending the scheme \eqref{eq:FEA} to
the case of the double-obstacle potential \eqref{eq:doF}, when \eqref{eq:muhimpl} needs to be replaced with a
variational inequality, is straightforward. We refer to \cite{BGN2013stable,vch} for the precise details.
}

We implemented the scheme \eqref{eq:FEA}, \revised{and its obstacle potential variant}, with the help of the finite element 
toolbox ALBERTA, see \cite{Alberta}. 
To increase computational efficiency, we employ adaptive meshes, which have a 
finer mesh size $h_{f} = \frac{\sqrt{2}}{N_f}$ within the diffuse interfacial 
regions and a coarser mesh size $h_{c} = \frac{\sqrt{2}}{N_c}$ away from them, 
with $N_f, N_c \in \N$, see \cite{voids3d,voids}
for a more detailed description. The nonlinear systems of equations arising 
from \eqref{eq:FEA} at each time step are solved with a Newton method, where we
employ the sparse factorization package UMFPACK, see \cite{Davis04}, 
for the solution of the linear systems at each iteration.
\revised{In the case of the double-obstacle potential, we employ the solution method from \cite{voids3d,vch}.}

In all our computations we fix the mobility $D(\bs{\nu}) =1$
and, up to possible rotations, use the anisotropy 
\begin{equation}
    \gamma(\vec p) = \sum_{\ell = 1}^d\sqrt{(1-\delta^2) p_\ell^2 + \delta^2|\vec p|^2},\qquad \vec p = (p_1,\cdots,p_d)^T,
    \label{eq:asyform}
\end{equation}
which can be regarded as a smoothed $\ell^1$-norm, with a small regularization
parameter $\delta>0$. \revised{Note that \eqref{eq:asyform} is a special case of \eqref{eq:bgn}.}
For the ({\bf DI}) computations we choose \revised{for} the potential $F$ 
\revised{either \eqref{eq:dwF}}, so that $\cF = \revised{\frac43}$,
\revised{or \eqref{eq:doF}, so that $\cF = \frac\pi2$}.
\revised{We let} $G$ \revised{be} defined by \eqref{DEF:G},
while the regularized mobility functions are defined via
\eqref{DEF:ME} and \eqref{DEF:BE}, with $r=2$ and $\gamma_0 = d_1 = 1$.
We also choose
$\alpha=\frac{\cFsqu}4$ 
so that \eqref{eq:anisf} is consistent with \eqref{eq:sharpv}.
\revised{Finally, unless otherwise stated we use the smooth potential \eqref{eq:dwF} for our ({\bf DI}) computations.} 

\subsection{2d results}

 \begin{figure}[!ht]
\centering
\includegraphics[angle=-0,width=0.8\textwidth]{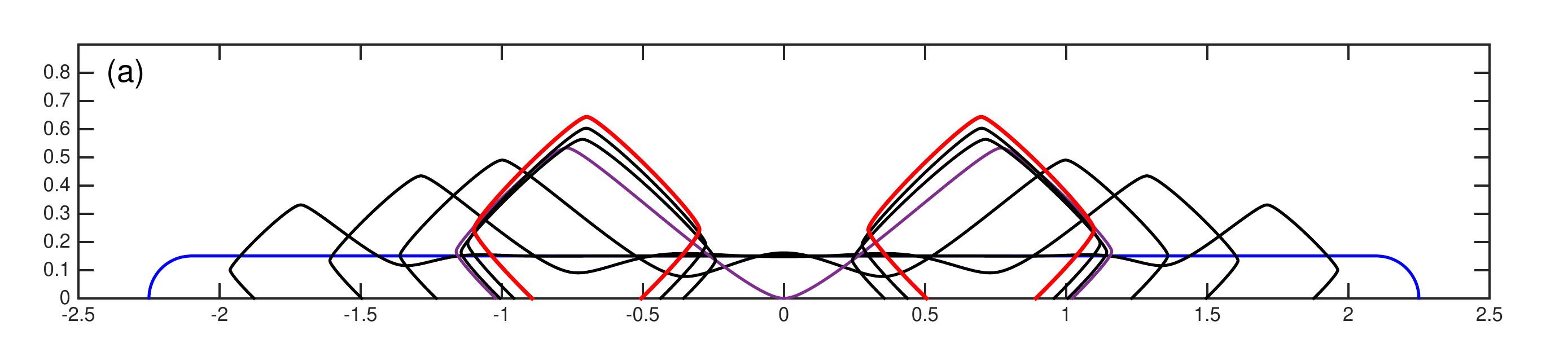}
\includegraphics[angle=-0,width=0.8\textwidth]{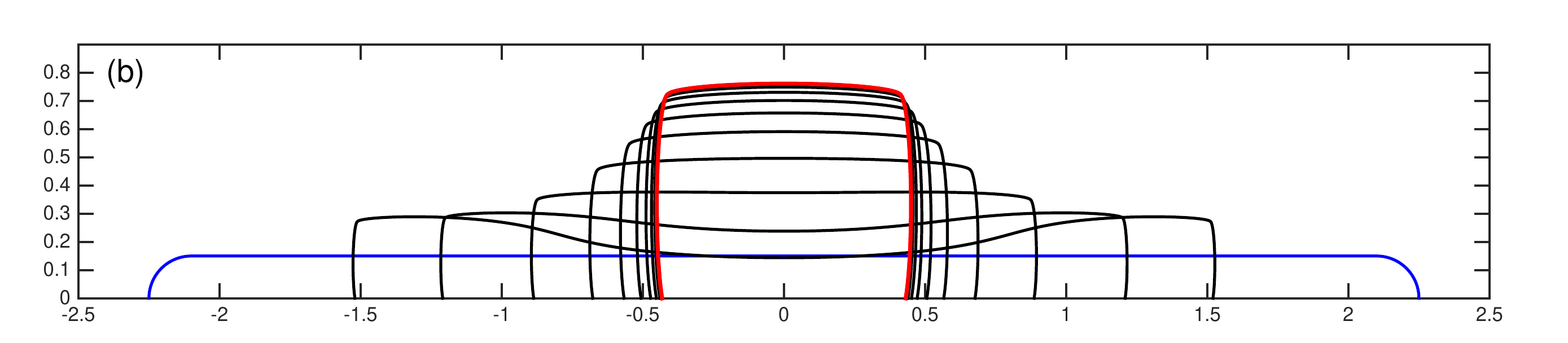}
\includegraphics[angle=-0,width=0.8\textwidth]{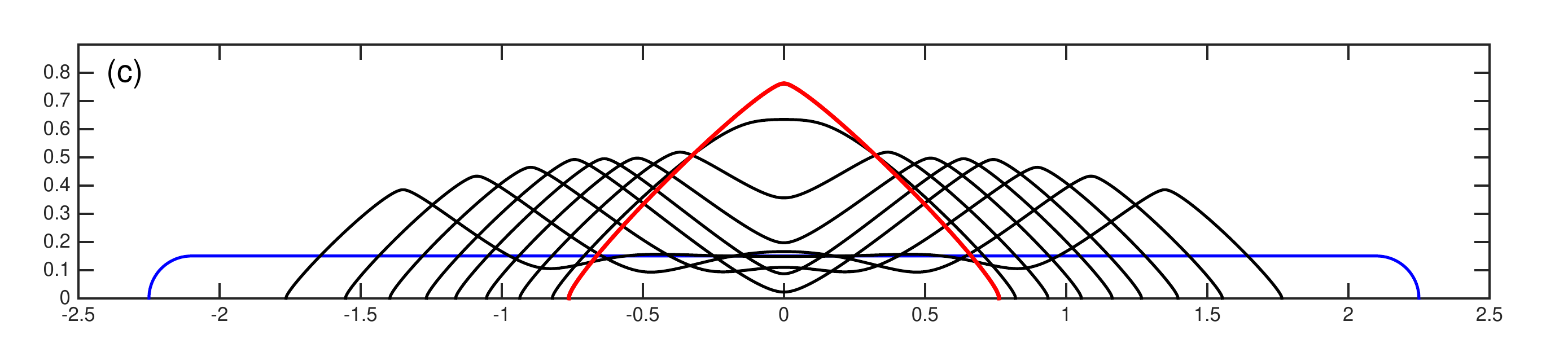}
\includegraphics[angle=-0,width=0.8\textwidth]{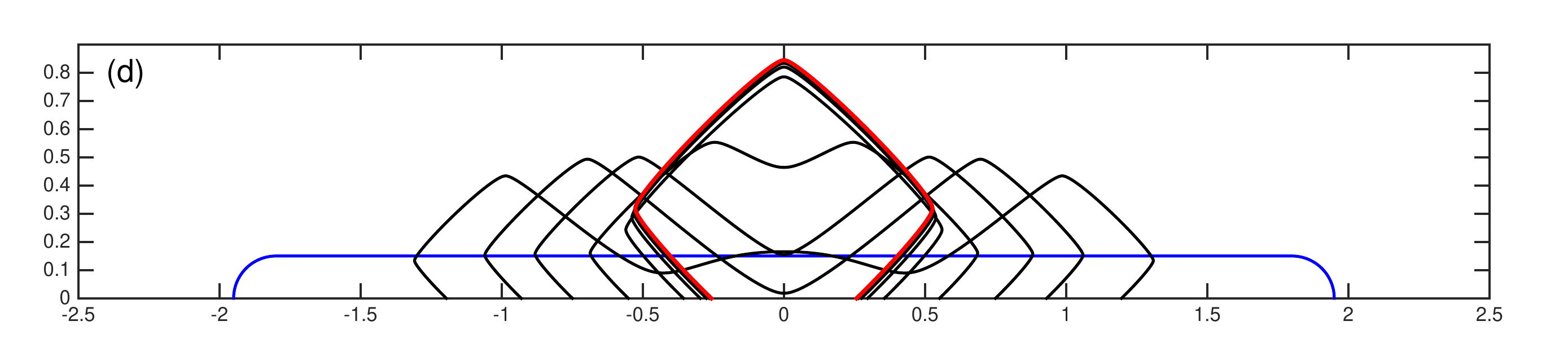}
\caption{ Evolution of small island films towards the equilibrium (red line) for the SI approximations.  (a) Plots at $t=0,0.002,0.01,0.02,0.030878,0.0319, 0.0339, 0.1$, where the island occurs pinch-off at $t=0.030878$; (b), (c) and (d) are the plots at $t=0,0.01,0.02,\cdots,0.1$. }
\label{fig:2disland}
\end{figure}

In numerical simulations of solid-state dewetting problems it is often of
interest whether a thin film of material breaks up into islands. For example,
in two space dimensions and in the isotropic case with a $90^\circ$ contact
angle condition it has been observed that 
elongated films undergo pinch-off once the aspect ratio of length versus height
goes beyond a critical value $R_0 \approx 127.9$, \cite{DBDLE06, Wang15}.
For nonzero values of $\sigma$, the critical value behaves like
$R_0 \approx 96.6/\sin(\frac{1}{2}\arccos\sigma)-8.66$, \cite{DBDLE06}.

It turns out that the anisotropy $\gamma$ can have a dramatic influence on the
critical value $R_0$. To investigate this numerically, we simulate the 
evolution of small thin films, starting from an initial 
interface in the form of the upper half of a tube with aspect ratio 
$R = L / H$, and fix $H = 0.3$. We consider the following four example 
situations:
\begin{itemize}
    \item [(a)]an island of $R=15$ with anisotropy $\gamma(\mR(\frac{\pi}{4})\vec p)$ and $\sigma =\cos\frac{5\pi}{6}$; 
    \item[(b)] an island of $R=15$ with anisotropy $\gamma(\vec p)$ and $\sigma =\cos\frac{5\pi}{6}$;
    \item[(c)] an island of $R=15$ with anisotropy $\gamma(\mR(\frac{\pi}{4})\vec p)$ and $\sigma =\cos\frac{\pi}{2}$;
    \item[(d)] an island of $R=13$ with anisotropy $\gamma(\mR(\frac{\pi}{4})\vec p)$ and $\sigma =\cos\frac{5\pi}{6}$,
\end{itemize}
where $\mR(\theta)$ is the rotation matrix with an angle $\theta$, and $\gamma(\vec p)$ is given by \eqref{eq:asyform} with $d=2$, $\delta = 0.1$. 
\revised{We note that anisotropies with a four-fold symmetry like our choices above are often used
in two-dimensional models for materials with a cubic crystalline surface energy \cite{LiuM93,McFaddenCS00,ZhangG03}.}

 Plots of the interface profiles for the SI approximations are presented in Fig.~\ref{fig:2disland}(a)-(d) for the four examples, respectively, where the approximated polygonal curve consists of 2048 line segments, and the time step size is fixed as $10^{-6}$. From these figures, we can observe the influence of the anisotropy $\gamma$, the contact energy density difference $\sigma$, and the aspect ratio $R$ of the thin film on the evolution.  
In particular, comparing the evolutions in Fig.~\ref{fig:2disland}(a) and (d) 
we see that the critical value $R_0$ for break-up to occur appears to satisfy
$13 < R_0 \leq 15$, which is much smaller than in the isotropic case.
Moreover, we see that either rotating the anisotropy, 
Fig.~\ref{fig:2disland}(b), or changing the contact angle, 
Fig.~\ref{fig:2disland}(c), ensures that no break-up occurs, meaning that
$R_0 > 15$ in both cases.

Let us remark that the pinch-off observed in Fig.~\ref{fig:2disland}(a)
represents a singularity for the parametric description on which the SI
approximations are based. Hence we perform a heuristical topological change,
from a single curve to two separate curves, once an inner vertex of the
polygonal curve touches the substrate.
In what follows we will use the computations in Fig.~\ref{fig:2disland} as
reference solutions for our DI approximations, in order to empirically
confirm our theoretical results from Section~\ref{se:asymptotic}.
 
\begin{figure}[!htp]
\centering
\includegraphics[angle=-0,width=0.4\textwidth]{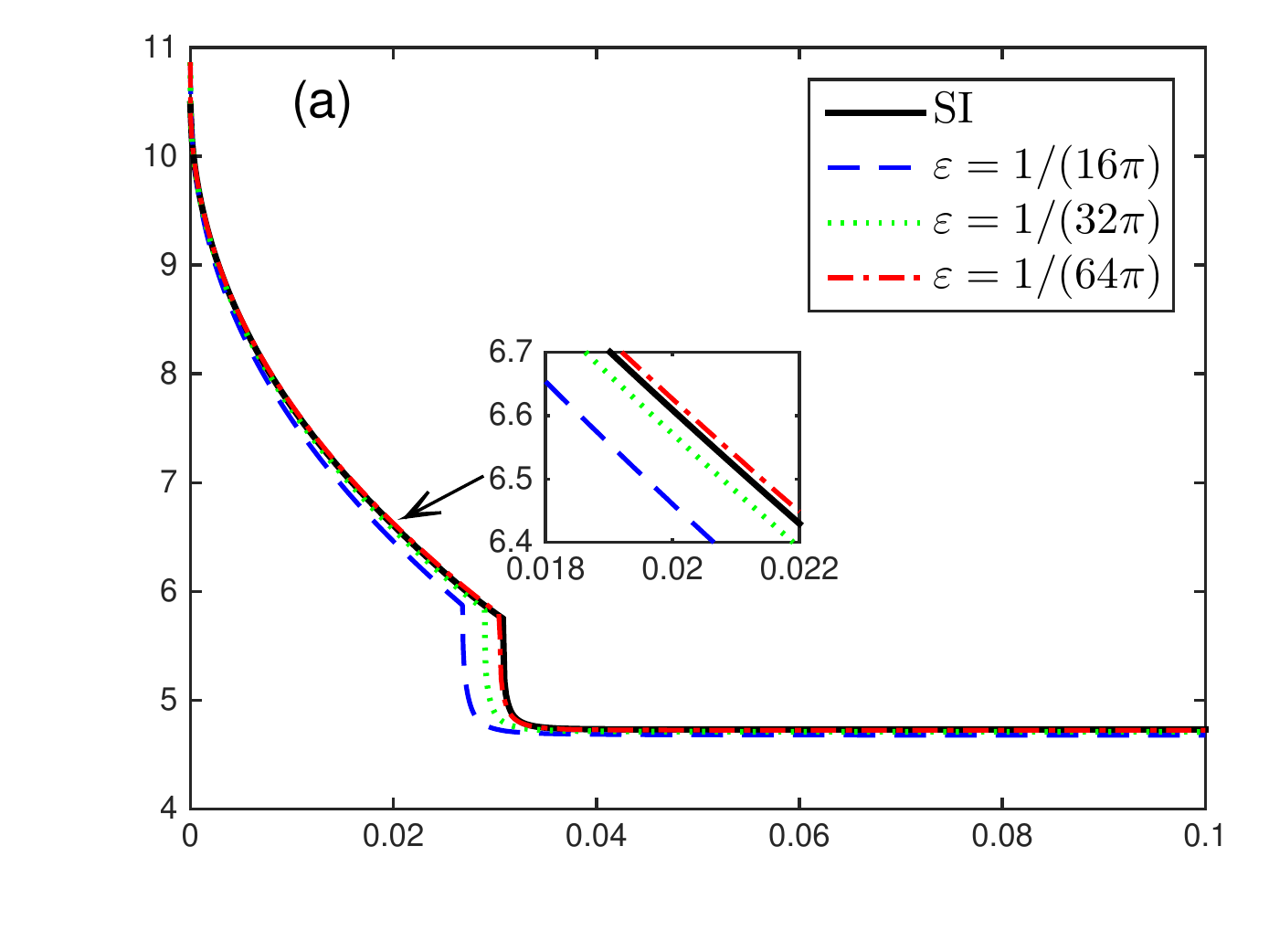}
\includegraphics[angle=-0,width=0.4\textwidth]{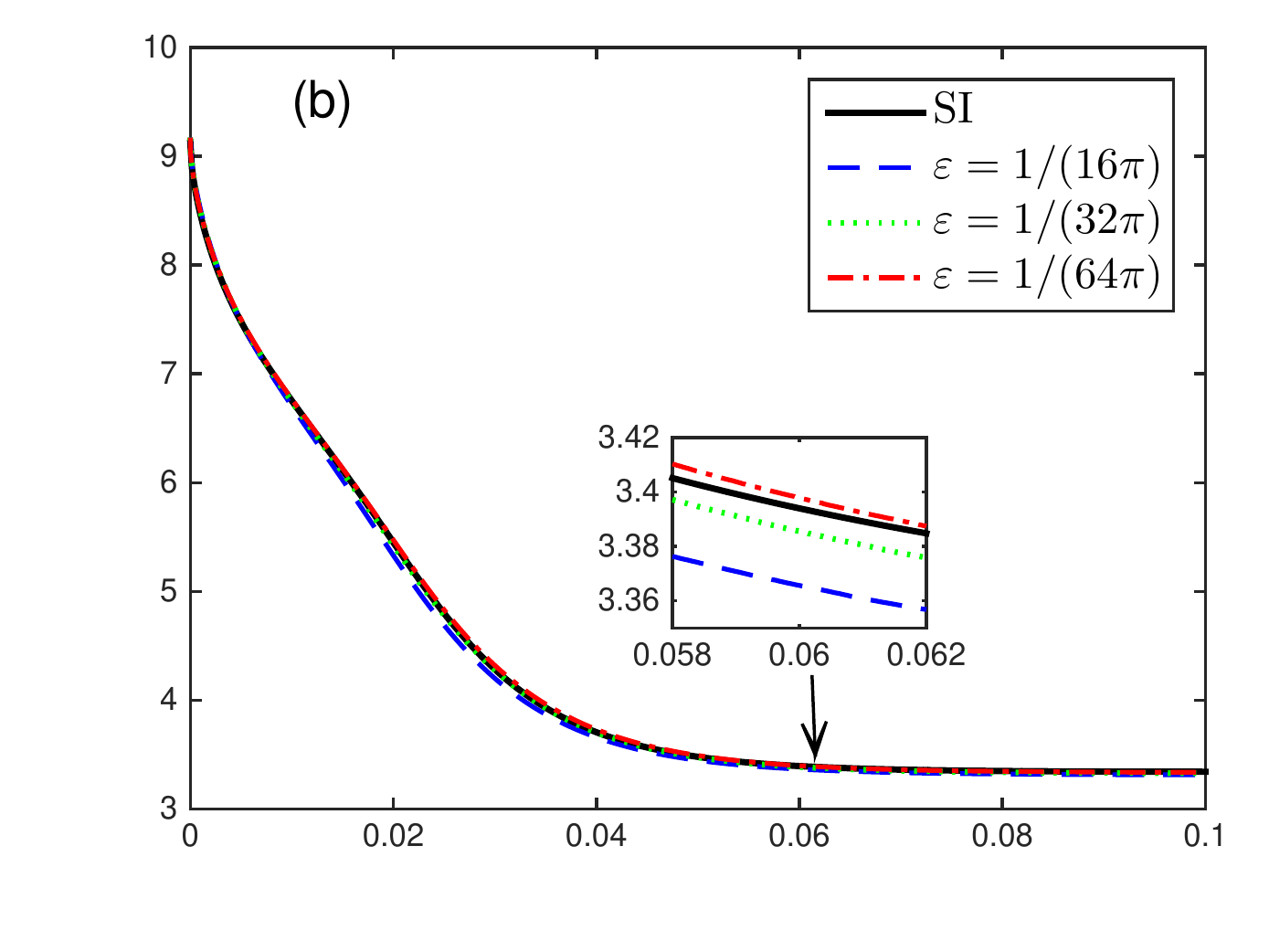}
\includegraphics[angle=-0,width=0.4\textwidth]{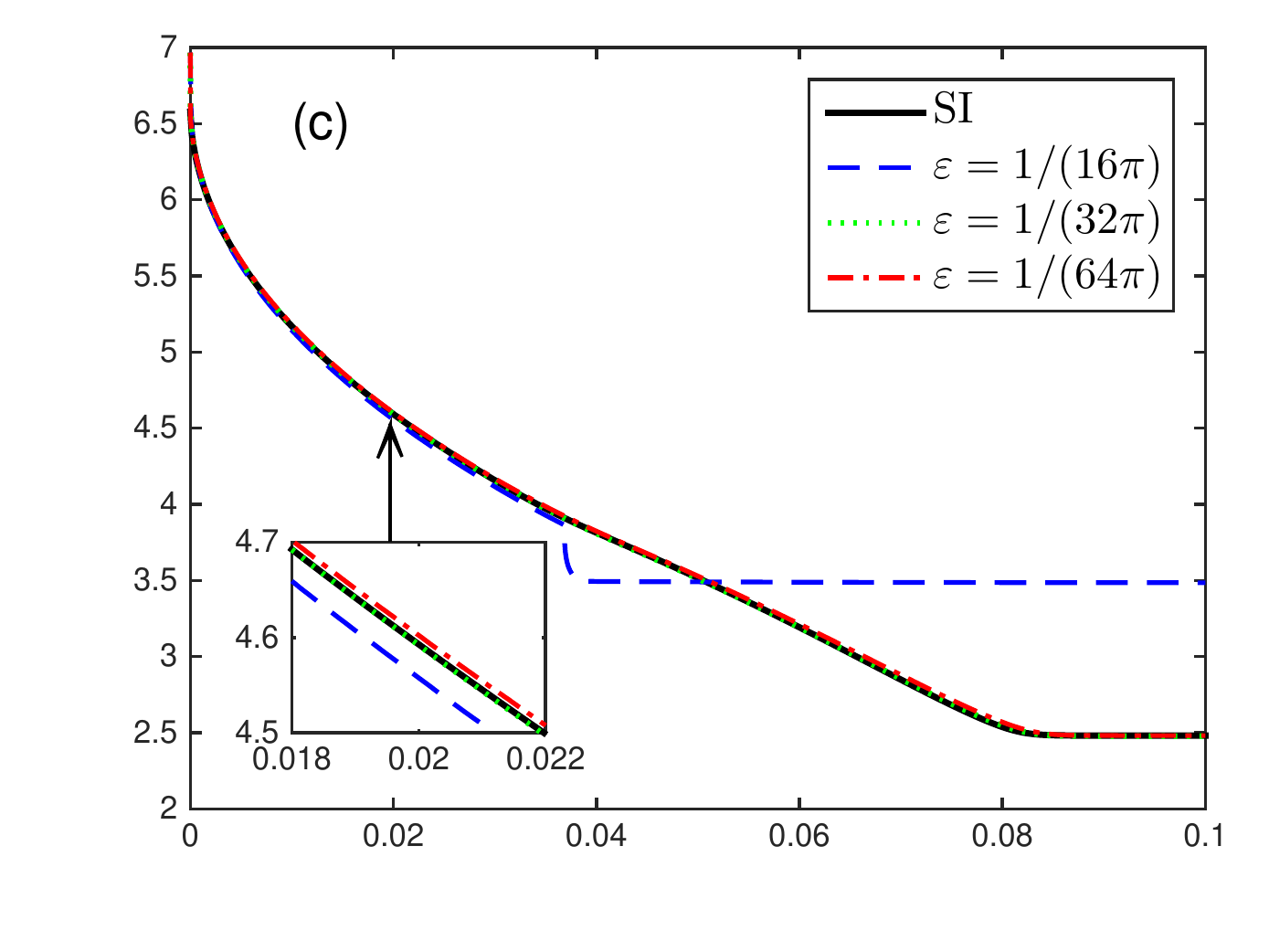}
\includegraphics[angle=-0,width=0.4\textwidth]{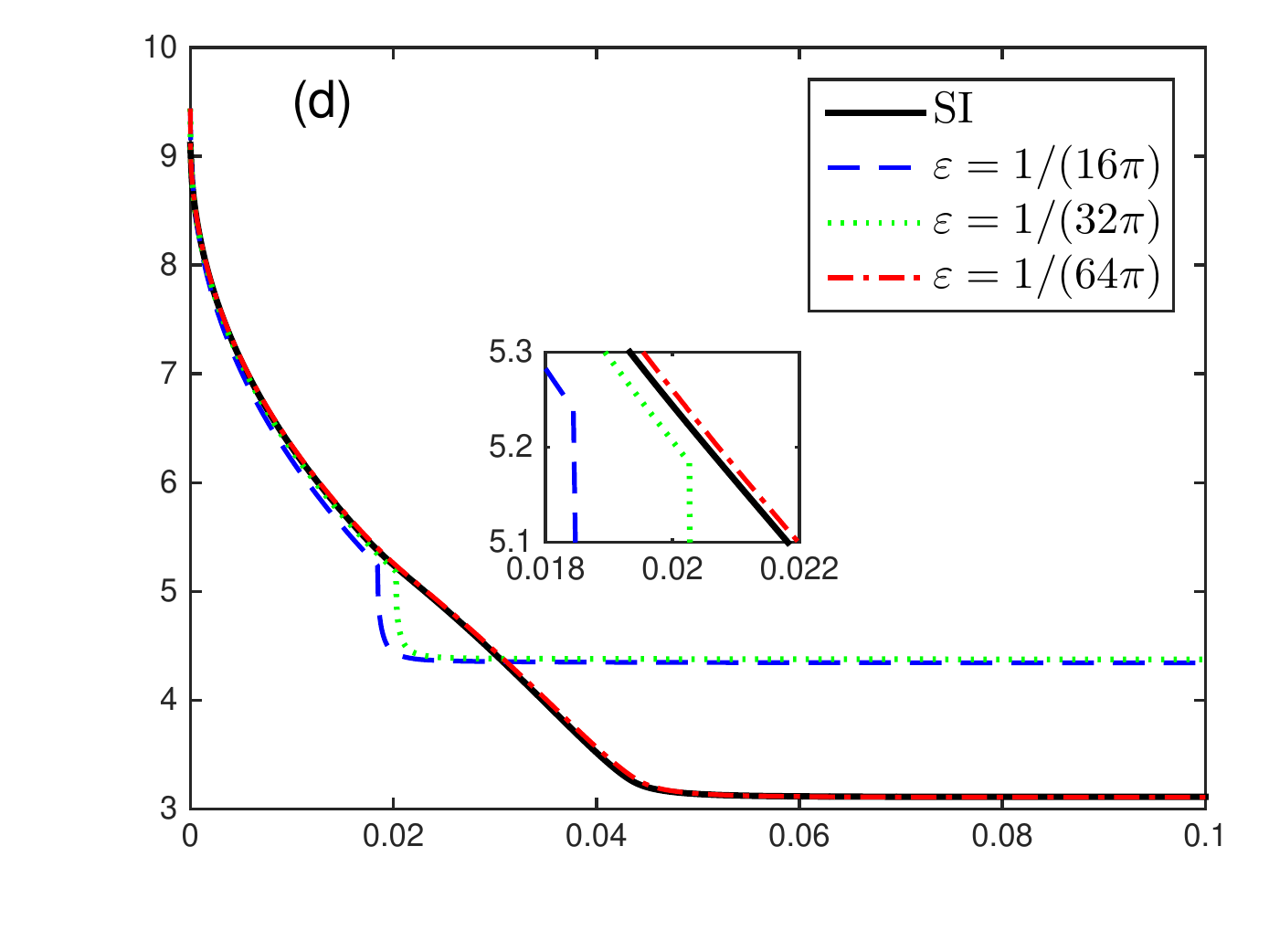}
\caption{The time history of the energy for the DI and SI approximations in the four different examples using the double smooth potential. }
\label{fig:2dEnergy}
\end{figure}

\begin{figure}
    \centering
    \includegraphics[width = 0.4\textwidth]{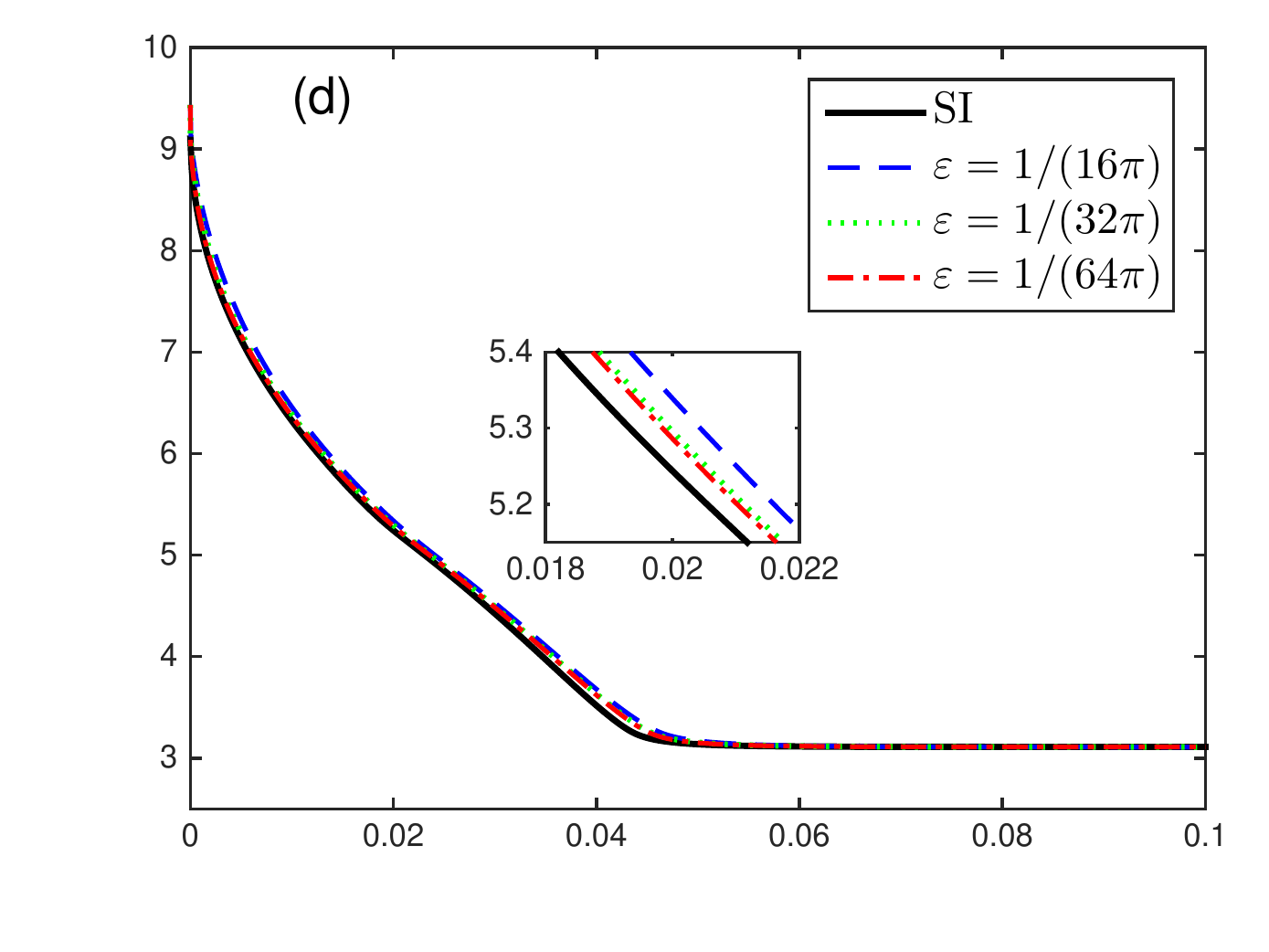}
    \includegraphics[width = 0.4\textwidth]{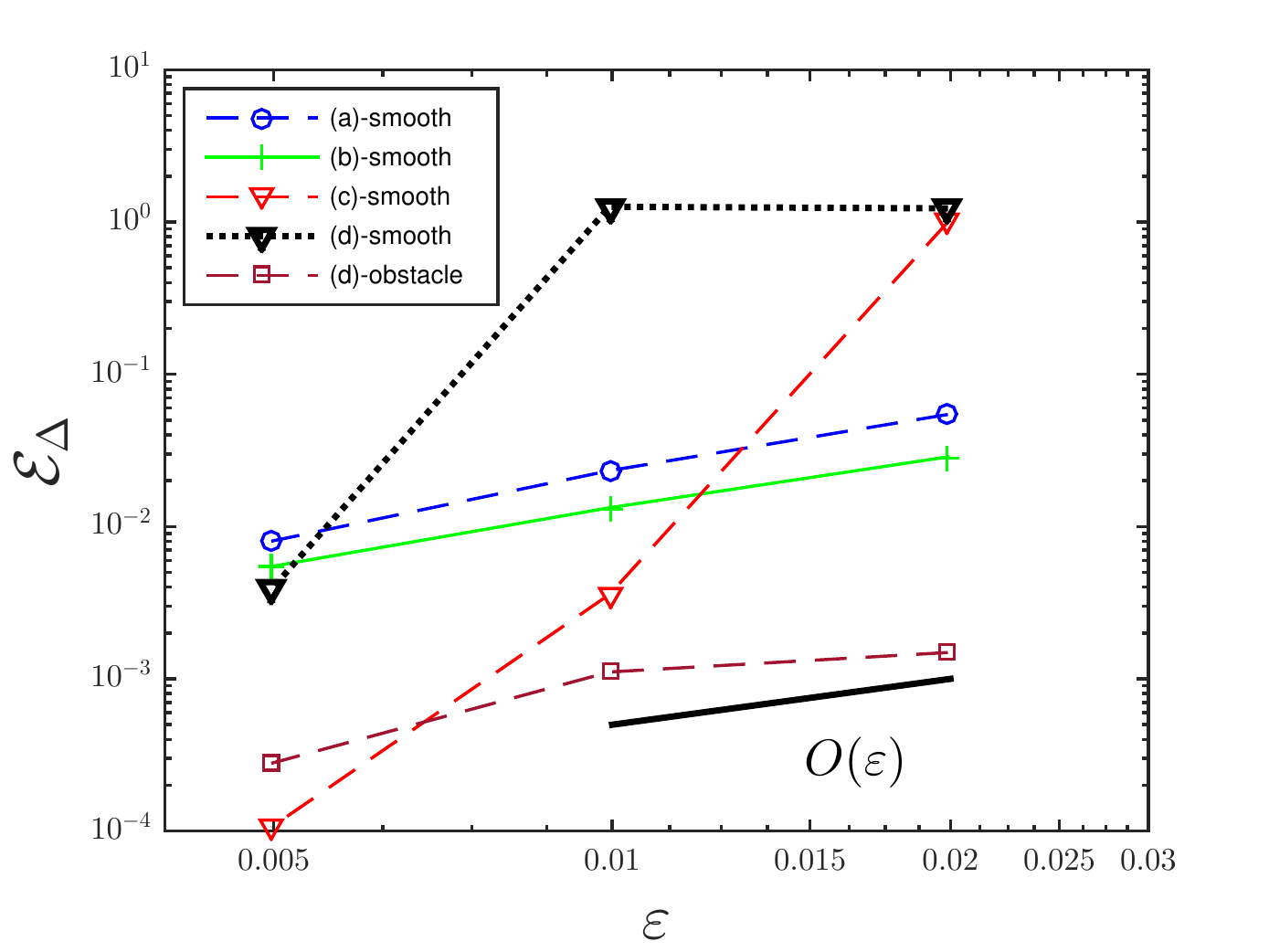}
    \caption{\revised{Left panel: The time history of the energy for the DI and SI approximations in Example (d) using the obstacle potential \eqref{eq:doF}. Right panel: The errors $\mathcal{E}_\Delta$ of the energy at the final time $T=0.1$ between the SI and DI approximations plotted against $\eps$. Here  ``(d)-obstacle'' refers to Example (d) with the obstacle potential, while the remaining graphs are for Examples (a)--(d) with the smooth potential.}}
    \label{fig:2dEnergyob}
\end{figure}

For our DI approximations we consider the computational domain 
$\Omega = [0,3]\times[0,1]$, on which for symmetry reasons we only compute the
right half of the evolving thin film. As interfacial parameters we consider
$\eps = 1 / (2^{4+i}\pi)$, for $i=0,\ldots,2$, 
and choose the discretization parameters as
$N_f = 2^{8+i}$, $N_c = 2^{5+i}$, $\Delta t = 10^{-3} / 2^{4+2i}$.
These spatial adaptive discretization parameters allow for a sufficient
resolution of the diffuse interface, while the temporal 
discretization parameters yield an excellent agreement with the SI
approximations. In fact, in Fig.~\ref{fig:2dEnergy} we show the energy plots 
of the DI approximations and compare them with the corresponding SI 
approximations for the four different examples from Fig.~\ref{fig:2disland}.
We observe that for sufficiently small values of $\eps$ there is excellent
agreement between the SI and DI evolutions, in line with our asymptotic 
analysis in Section~\ref{se:asymptotic}.
What is interesting to note is that for Example~(a) the pinch-off time 
predicted by the DI computations is too early when $\eps$ is not small, and
this can be explained by the fact that the wider interfacial region ``sees''
contact with the substrate earlier, leading to the break-up into two islands.
For the same reason, in Examples~(c) and (d) the DI computations for 
$\eps=1/(16\pi)$ erroneously predict a pinch-off, leading to a larger final
energy. But once $\eps$ is sufficiently small, no pinch-off occurs, in
agreement with the SI evolutions.

\revised{We note that using the double-obstacle potential \eqref{eq:doF} leads to
very similar results. As an example we show the evolution of the discrete energies
for Example (d) in Fig.~\ref{fig:2dEnergyob}. In addition, in order to also have a 
quantitative comparison between our SI and DI computations, in the same figure we
also present plots of the energy difference $\mathcal{E}_\Delta$ between the
final SI and DI solutions against $\eps$. The presented results suggest that the DI
energies of the final states approach the corresponding SI energy with $O(\eps)$.
Note that the three instances where $\mathcal{E}_\Delta \geq 10^{-1}$ correspond to cases where
the DI computations wrongly predict a pinch-off. Moreover, in practice, we observe that the contact angles between DI and SI at the final time agree very well, with the error being of order $10^{-3}$ throughout.}

\begin{figure}[!htp]
\centering
\includegraphics[angle=-0,width=0.45\textwidth]{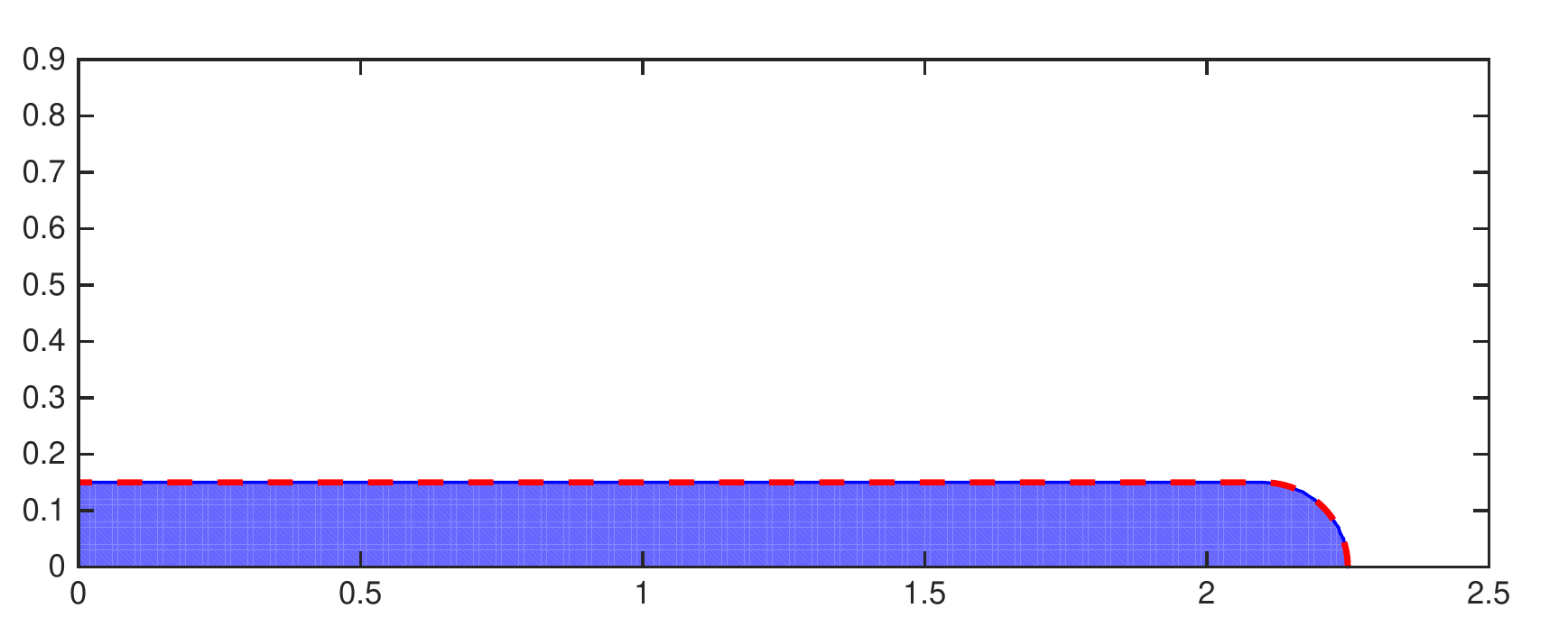}\;\;
\includegraphics[angle=-0,width=0.45\textwidth]{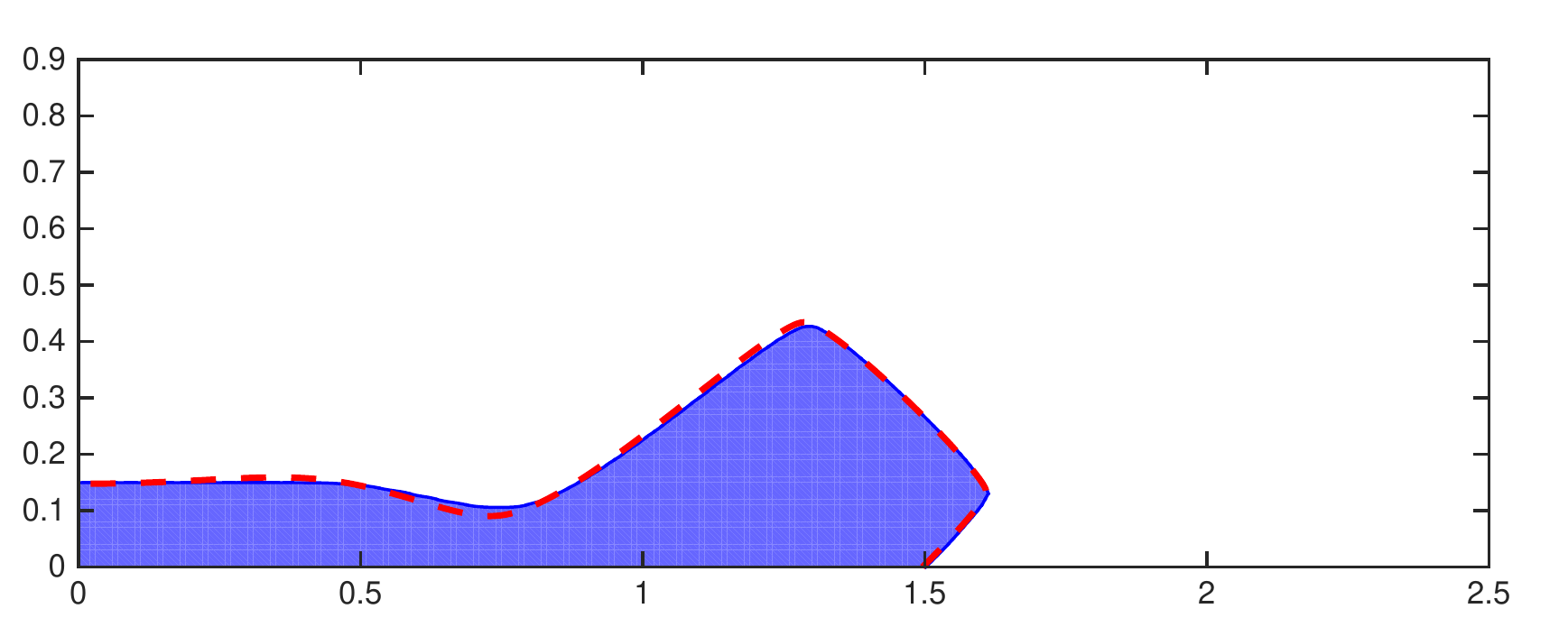}
\includegraphics[angle=-0,width=0.45\textwidth]{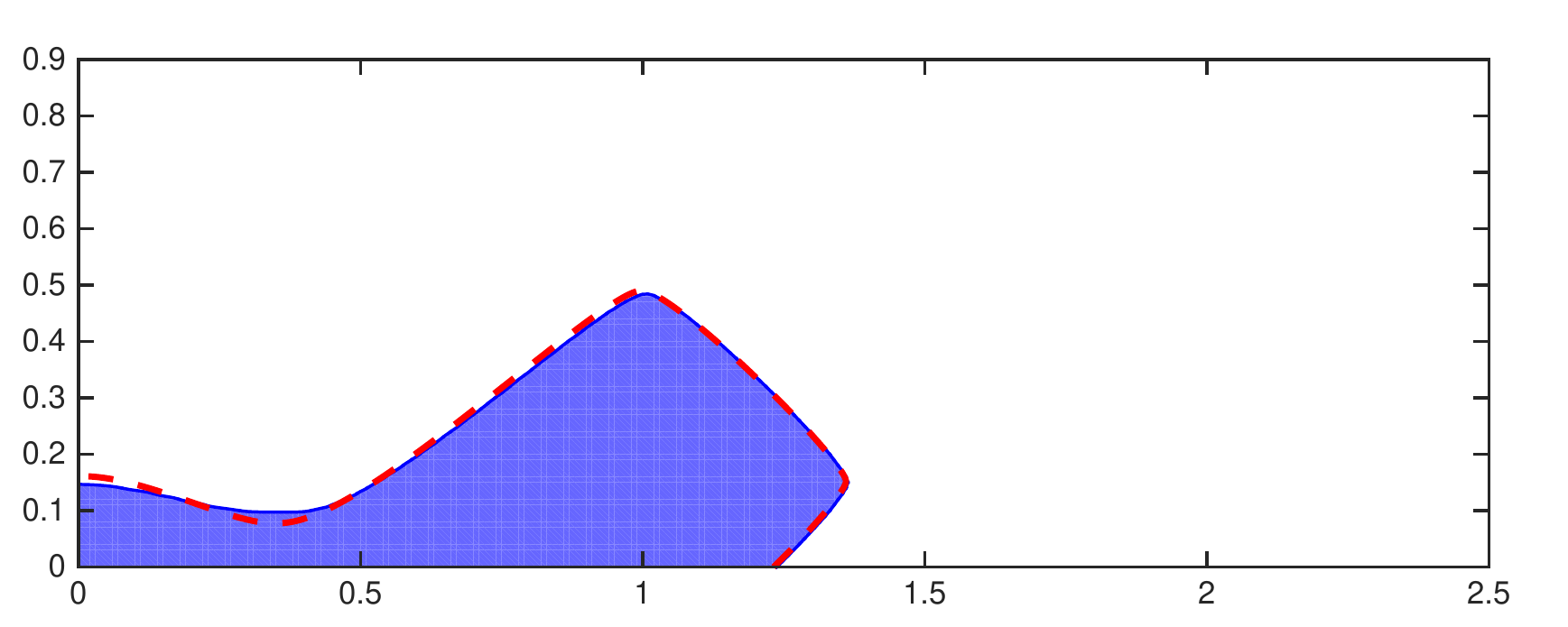}\;\;
\includegraphics[angle=-0,width=0.45\textwidth]{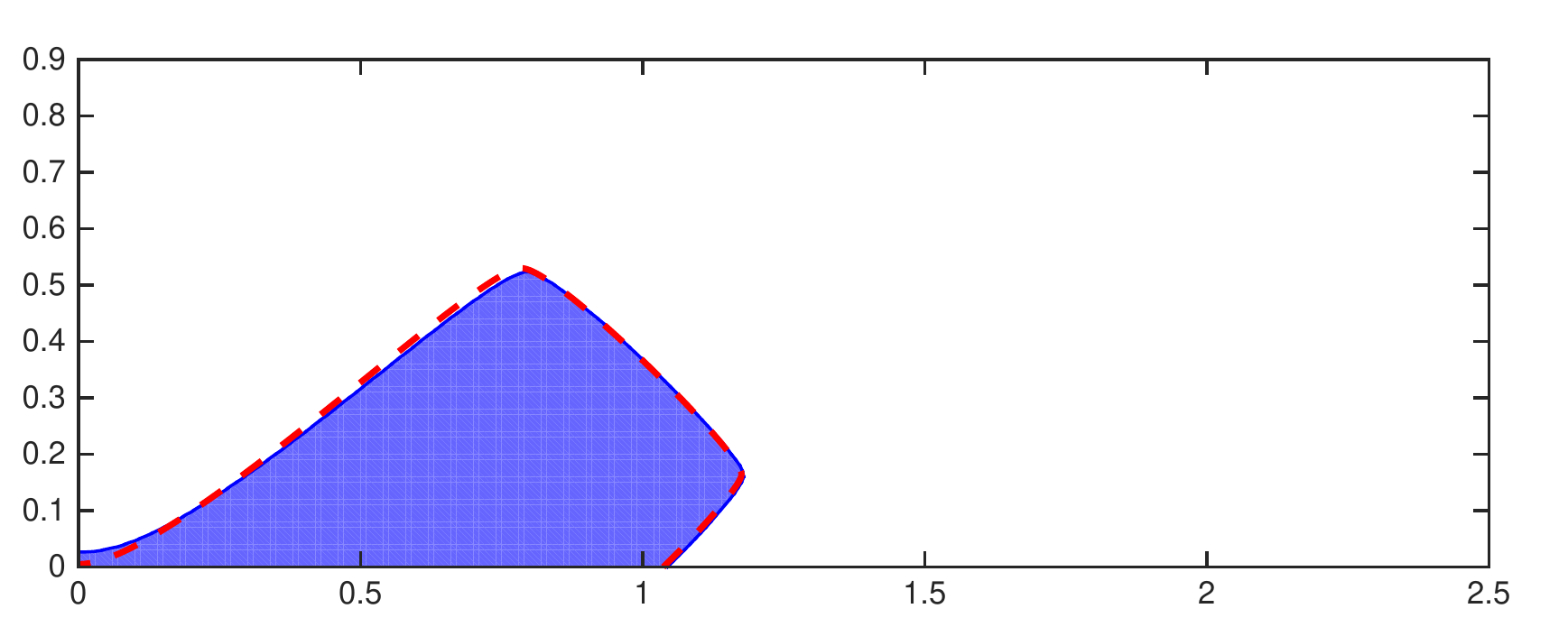}
\includegraphics[angle=-0,width=0.45\textwidth]{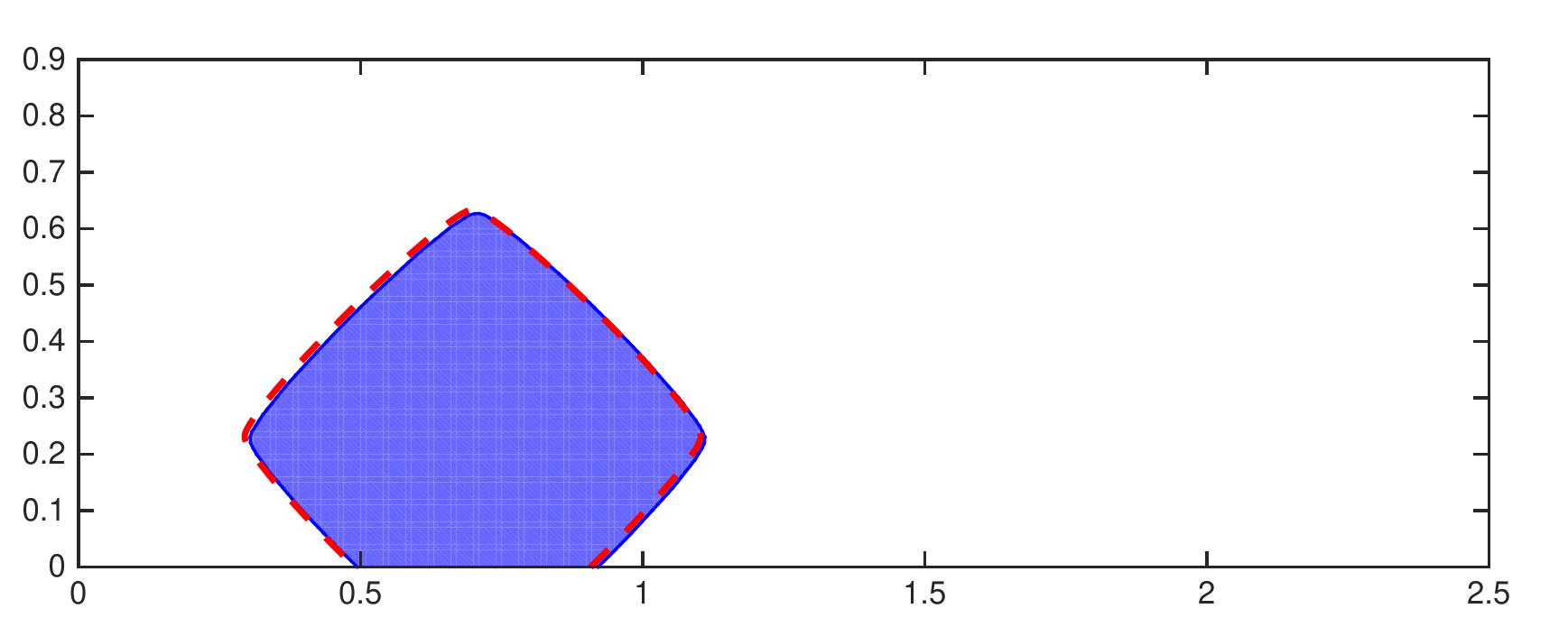}\;\;
\includegraphics[angle=-0,width=0.45\textwidth]{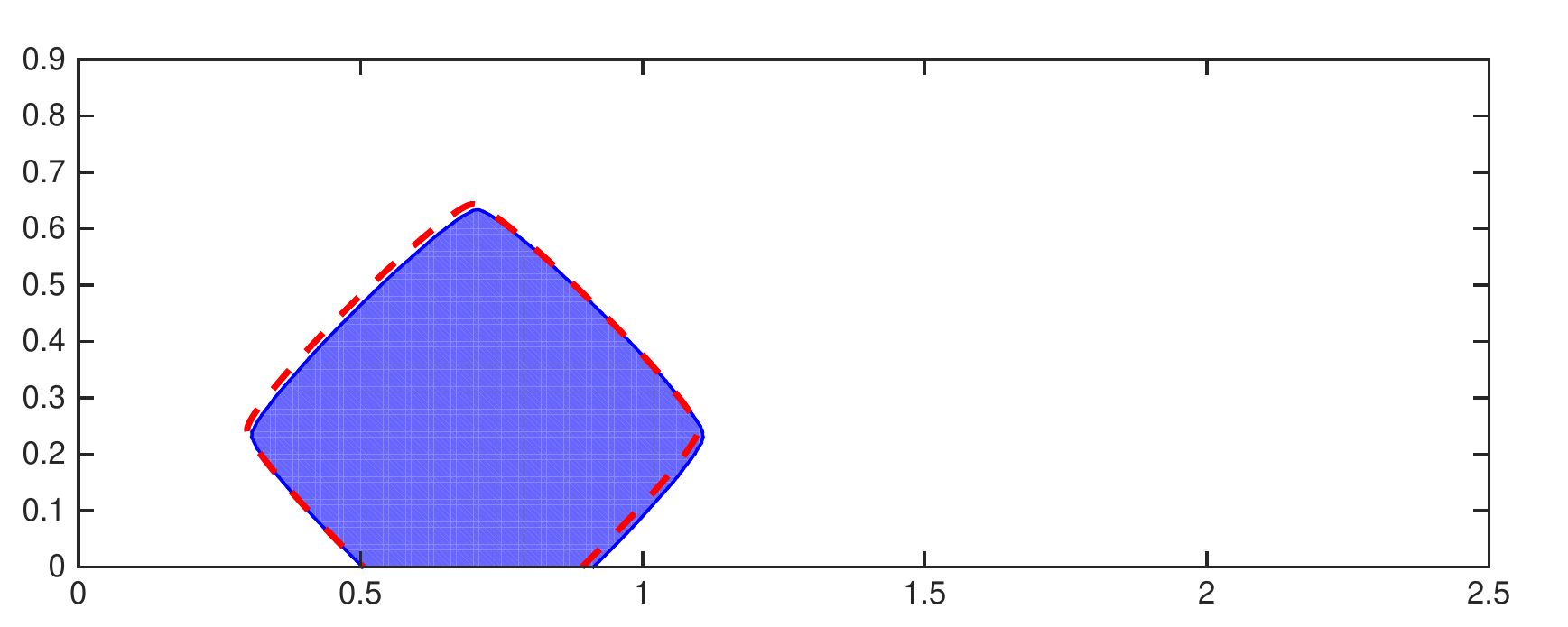}
\caption{[\;{\bf Example~(a)}\;] Interface profiles at times $t=0,0.01,0.02,0.03,0.04,0.1$ for the DI approximations with $\eps = 1/(64\pi)$, and the red dash line represents the SI approximations.} 
\label{fig:R15}
\end{figure}

\begin{figure}[!htp]
\centering
\includegraphics[angle=-0,width=0.45\textwidth]{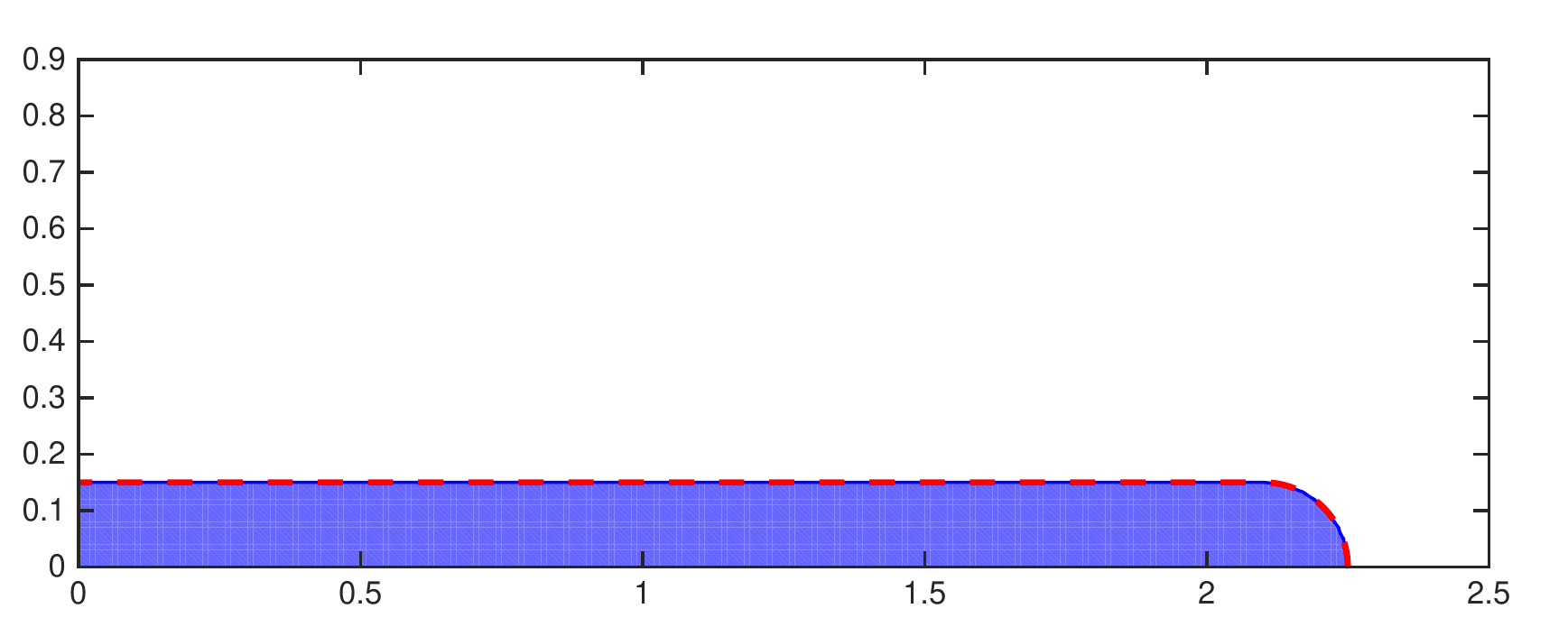}\;\;
\includegraphics[angle=-0,width=0.45\textwidth]{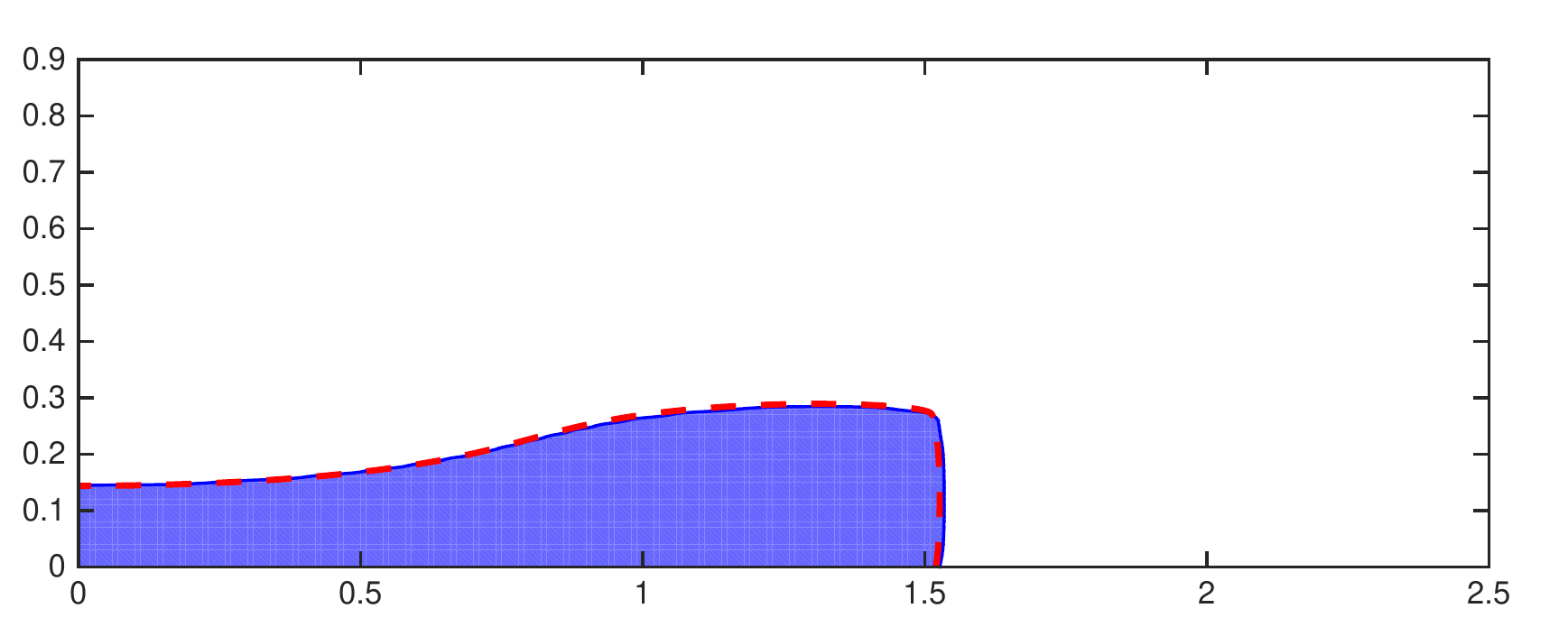}
\includegraphics[angle=-0,width=0.45\textwidth]{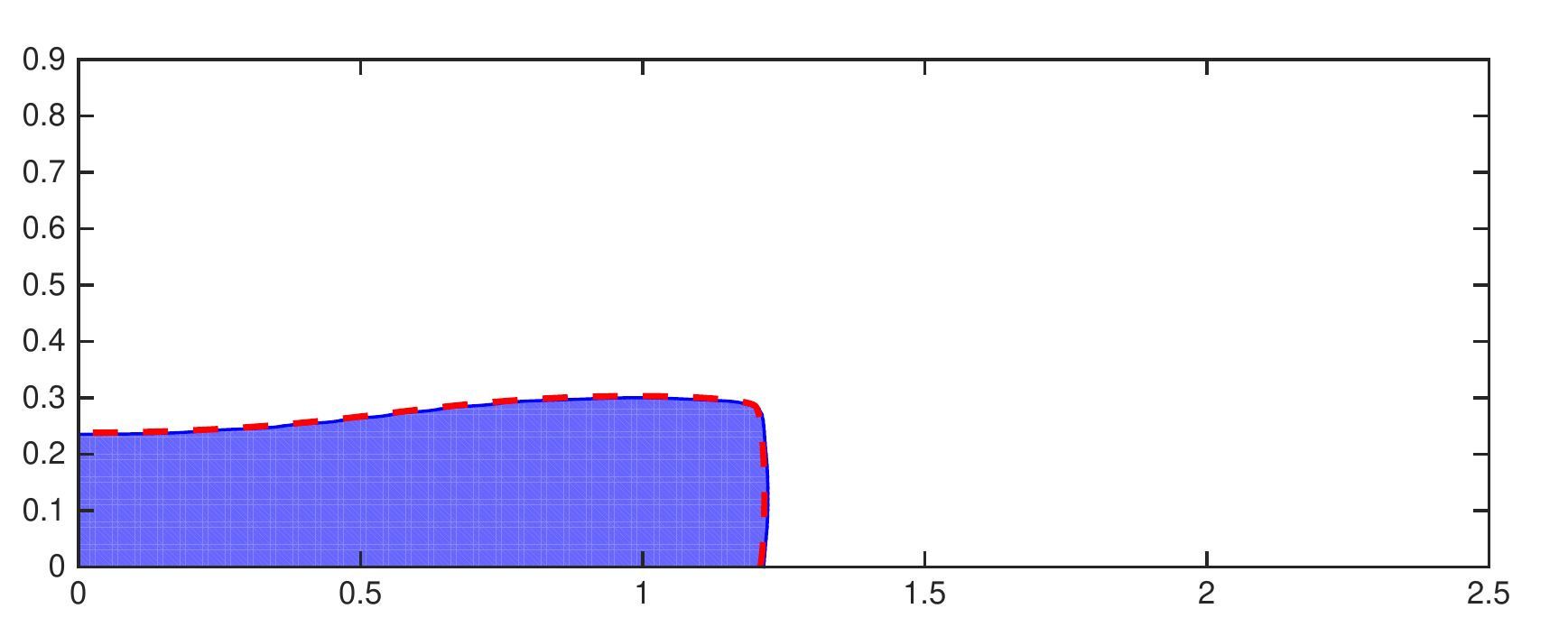}\;\;
\includegraphics[angle=-0,width=0.45\textwidth]{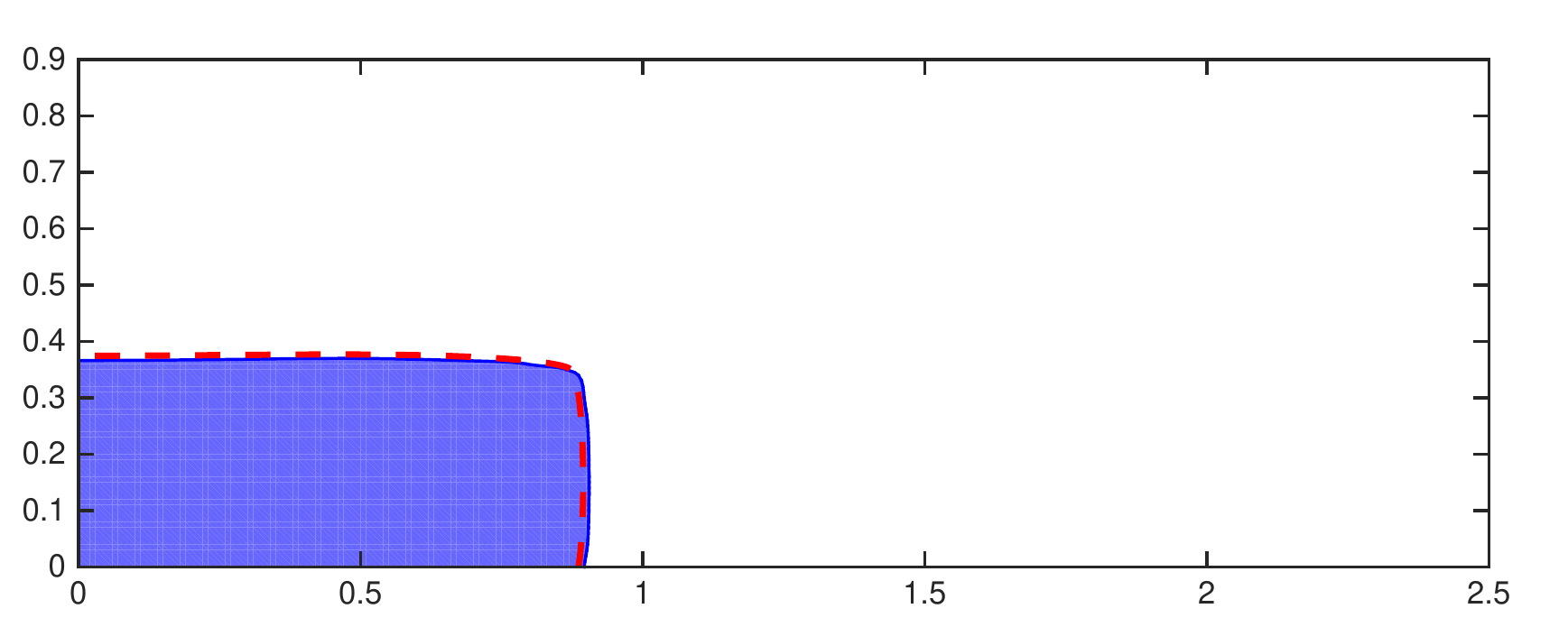}
\includegraphics[angle=-0,width=0.45\textwidth]{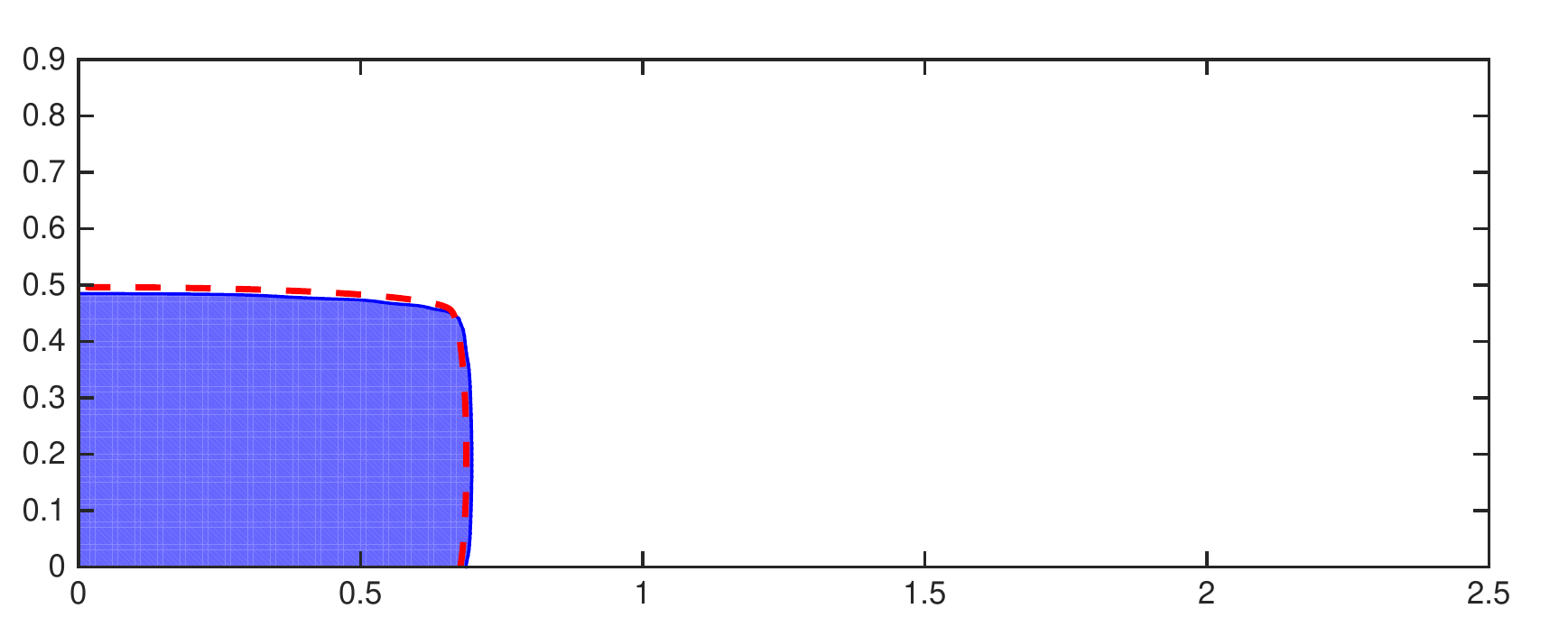}\;\;
\includegraphics[angle=-0,width=0.45\textwidth]{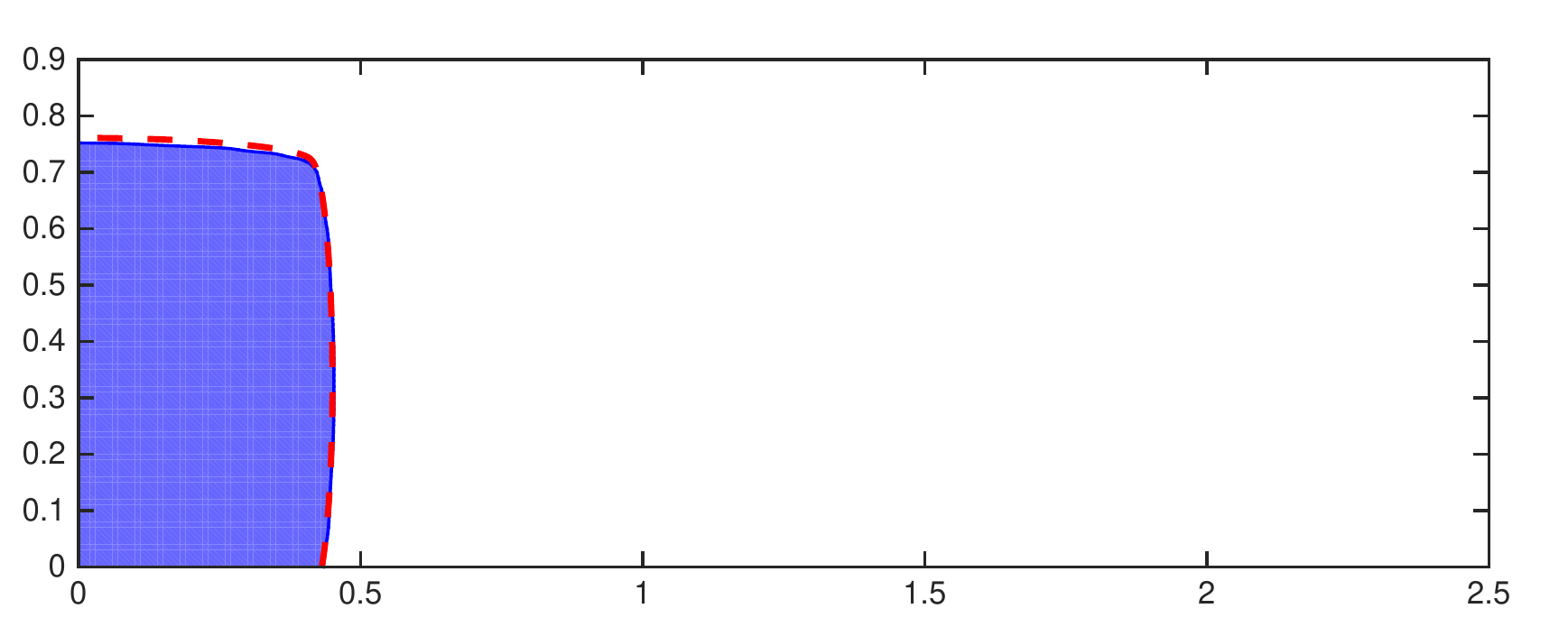}
\caption{[\;{\bf Example~(b)}\;] Interface profiles at times $t=0,0.01,0.02,0.03,0.04,0.1$ for the DI approximations with $\eps = 1/(64\pi)$, and the red dash line represents the SI approximations.} 
\label{fig:R15NR}
\end{figure}

\begin{figure}[!htp]
\centering
\includegraphics[angle=-0,width=0.45\textwidth]{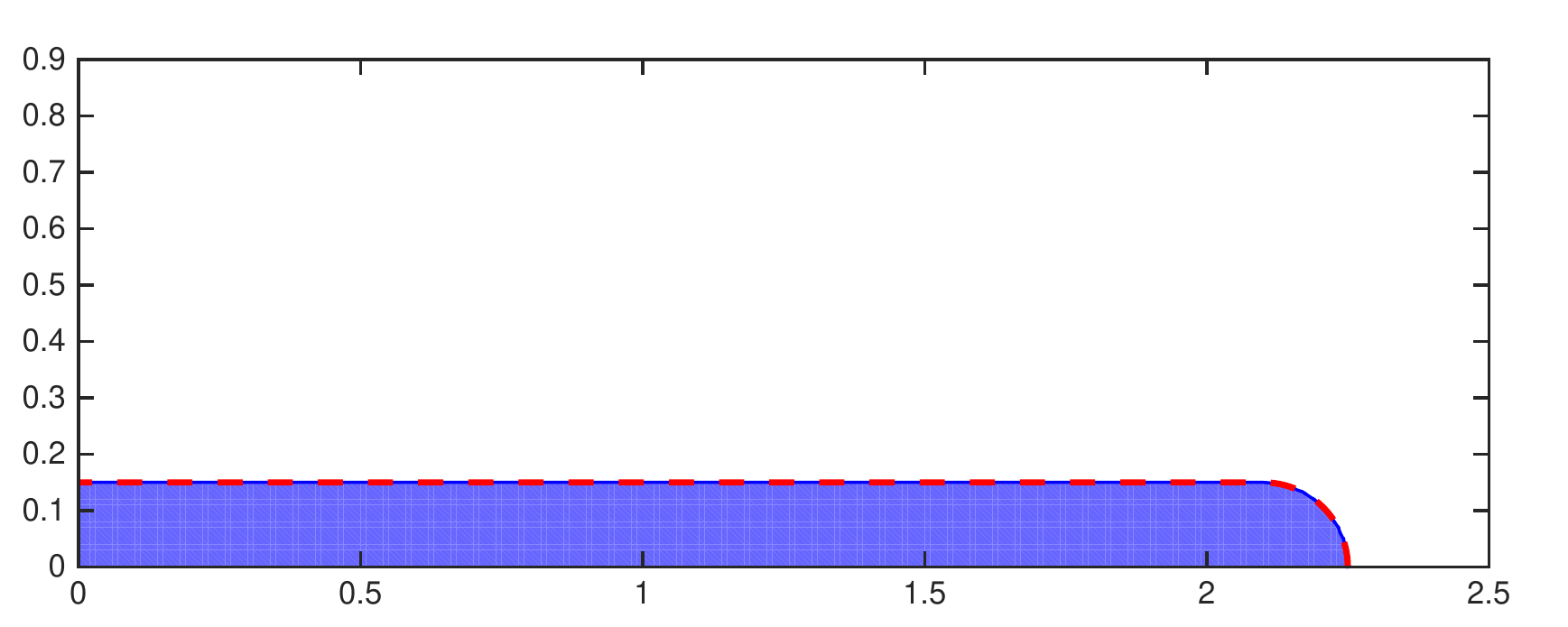}\;\;
\includegraphics[angle=-0,width=0.45\textwidth]{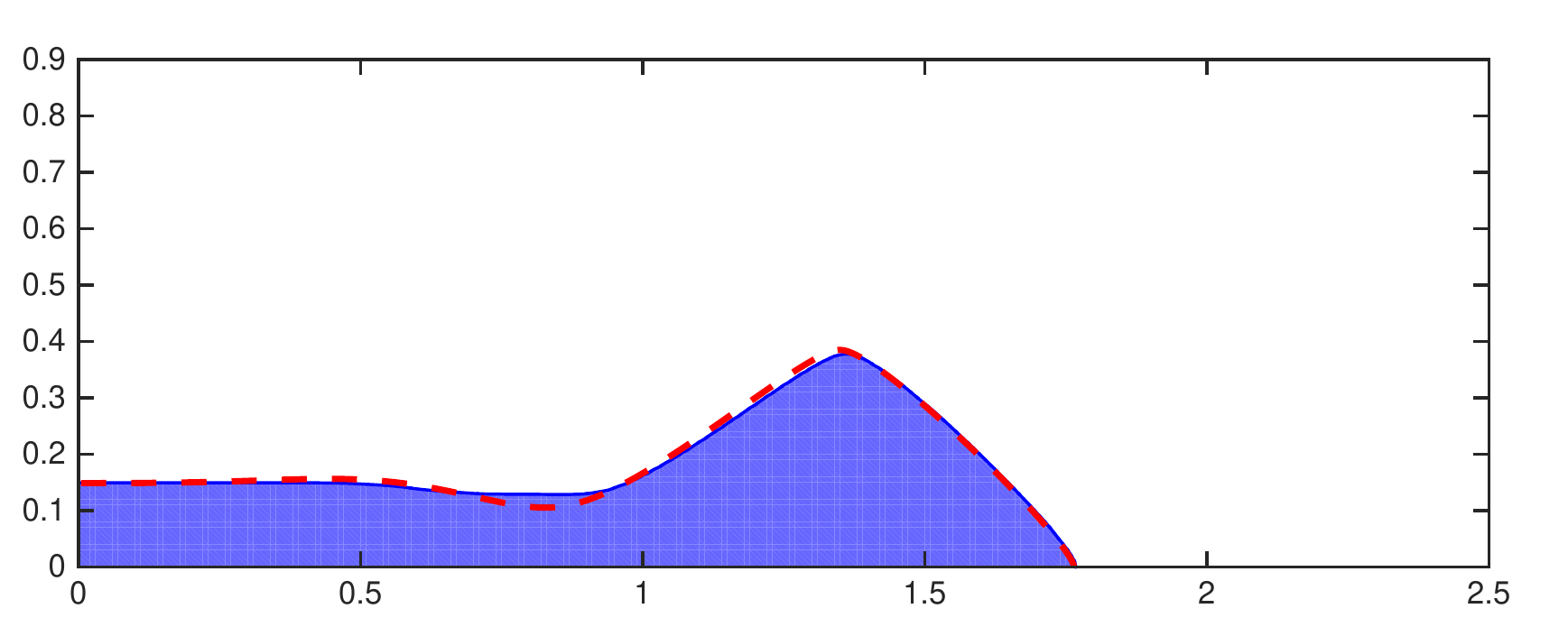}
\includegraphics[angle=-0,width=0.45\textwidth]{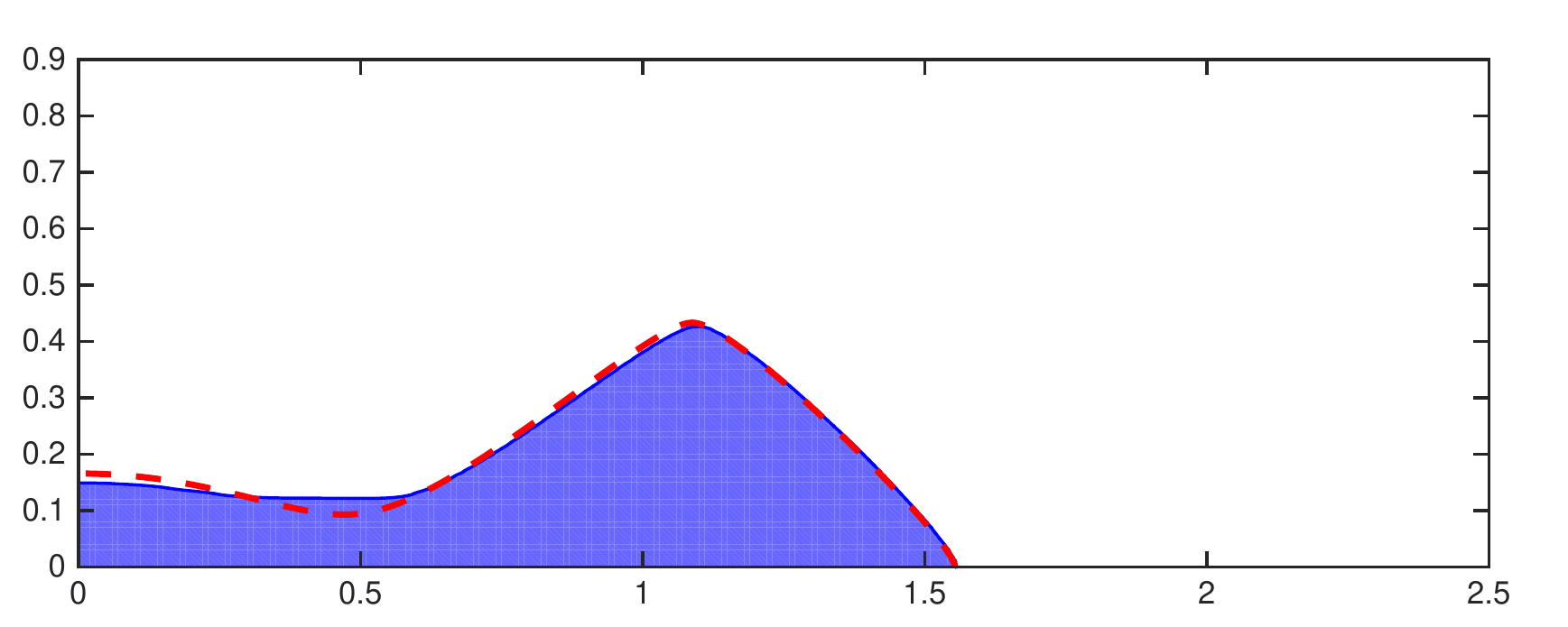}\;\;
\includegraphics[angle=-0,width=0.45\textwidth]{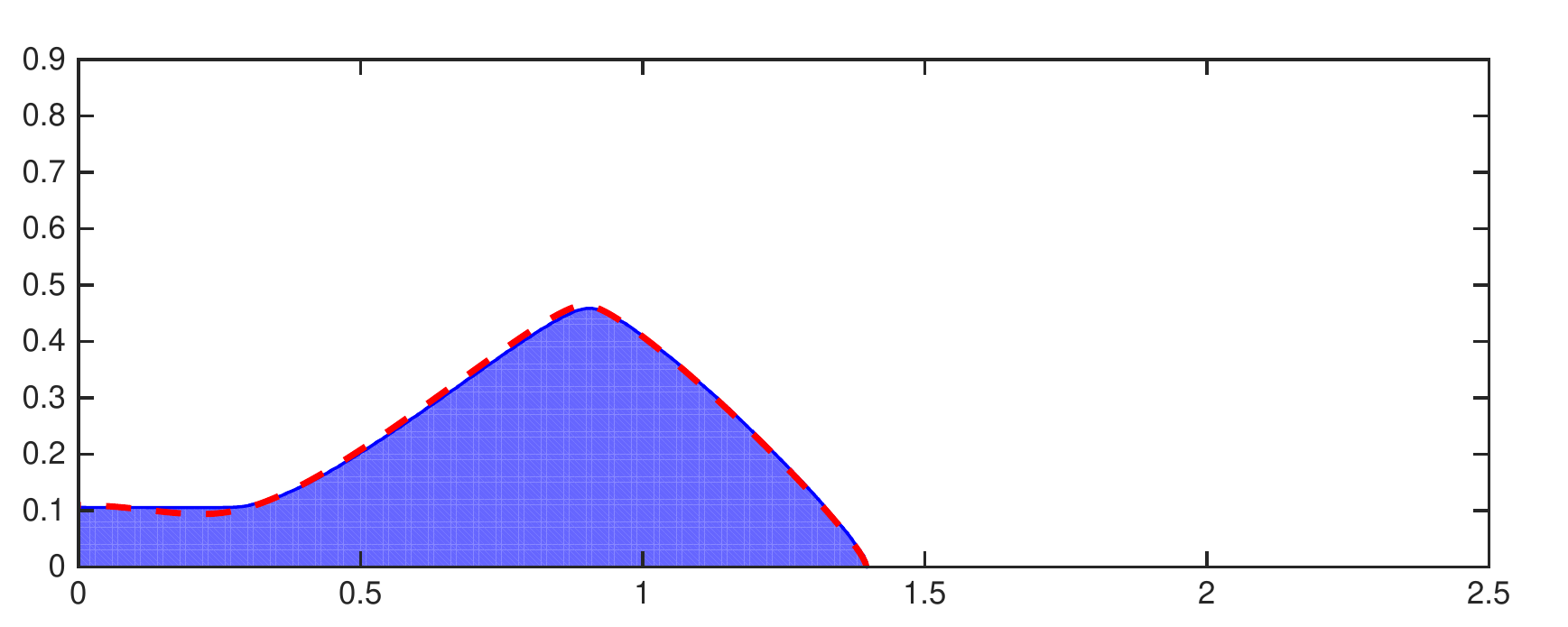}
\includegraphics[angle=-0,width=0.45\textwidth]{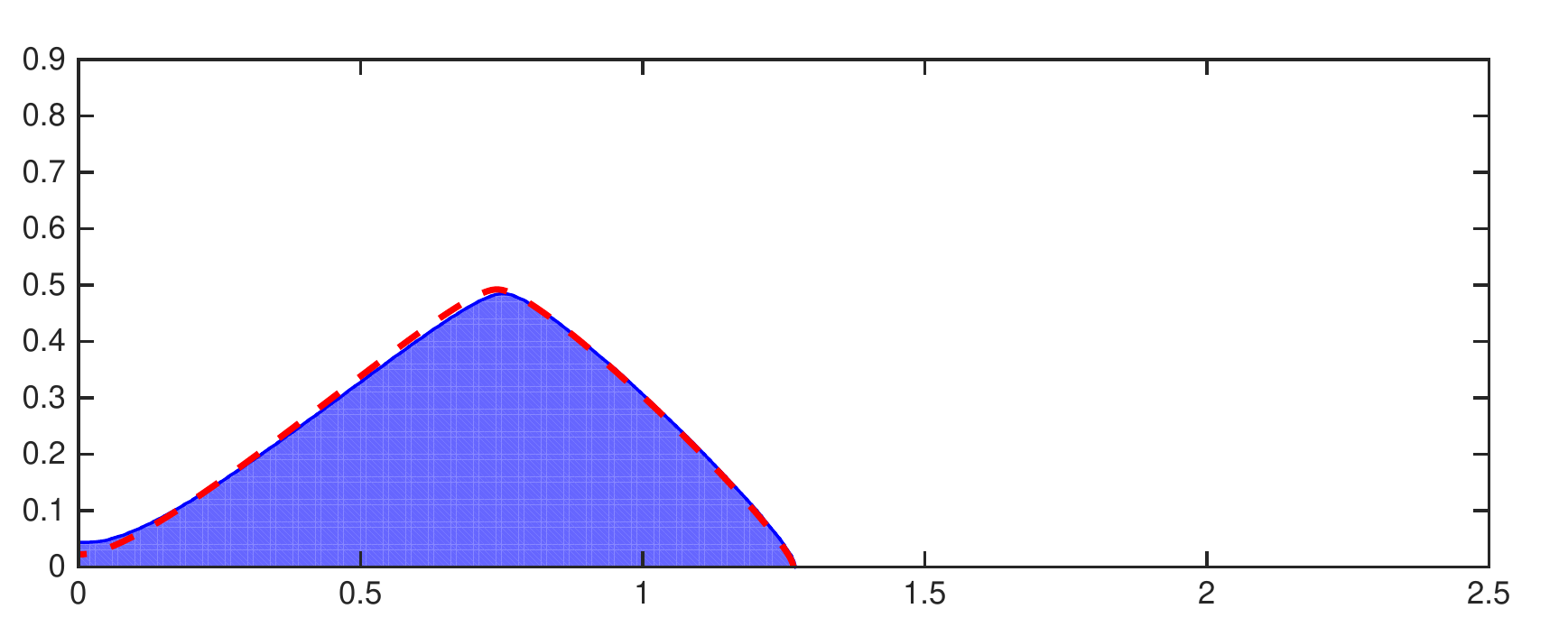}\;\;
\includegraphics[angle=-0,width=0.45\textwidth]{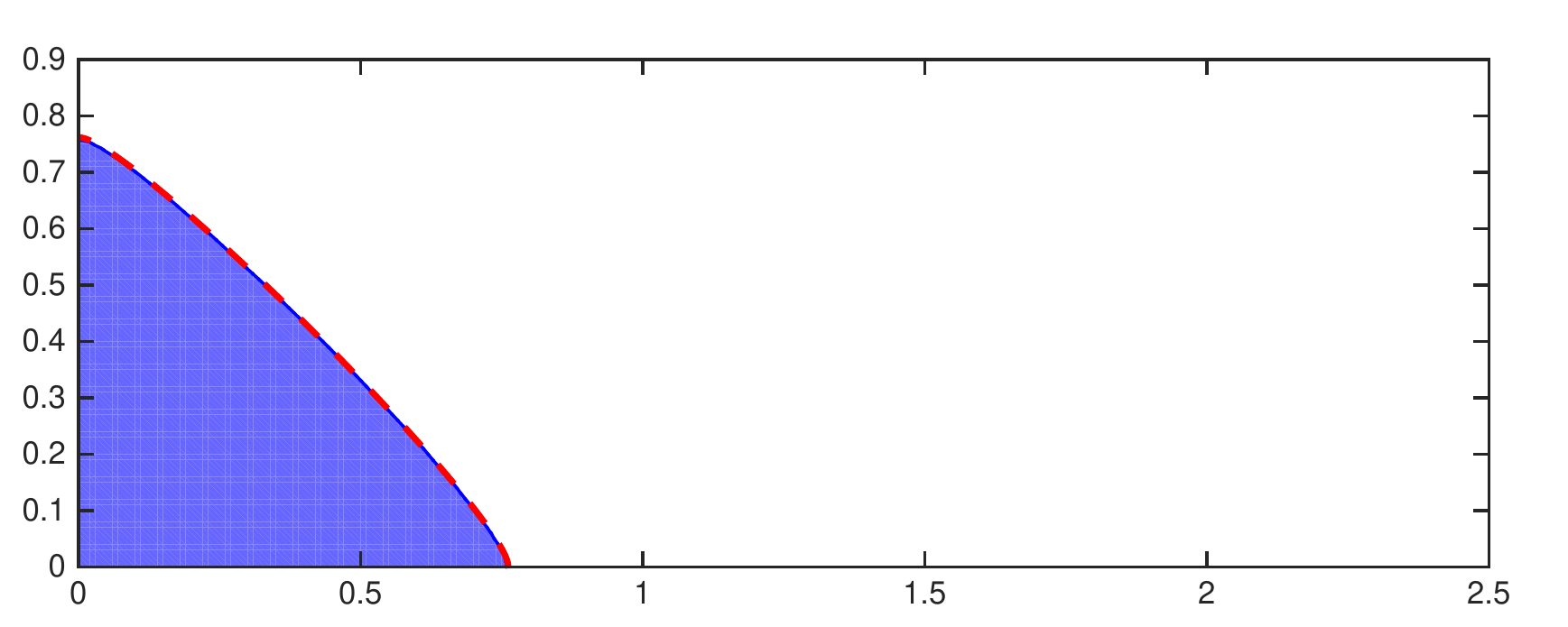}
\caption{[\;{\bf Example~(c)}\;] Interface profiles at times $t=0,0.01,0.02,0.03,0.04,0.1$ for the DI approximations with $\eps = 1/(64\pi)$, and the red dash line represents the SI approximations.} 
\label{fig:R150}
\end{figure}

\begin{figure}[!htp]
\centering
\includegraphics[angle=-0,width=0.45\textwidth]{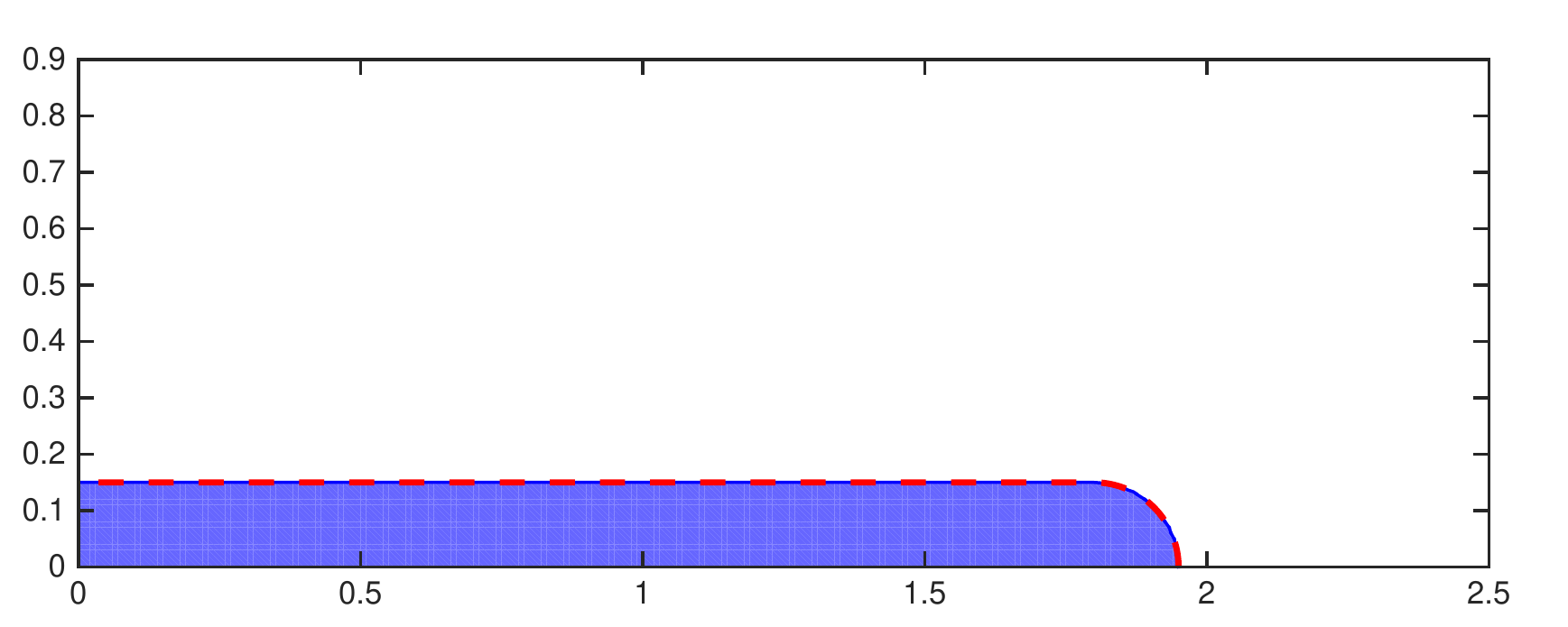}\;\;
\includegraphics[angle=-0,width=0.45\textwidth]{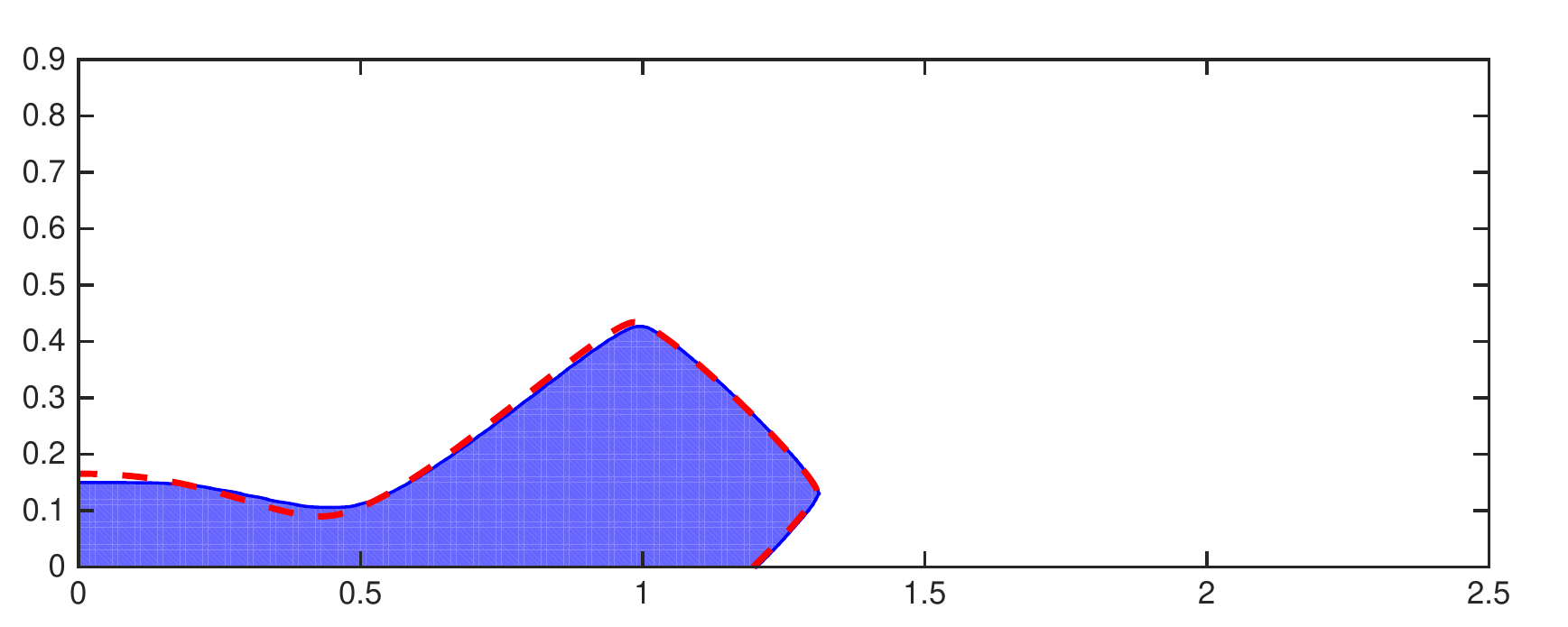}
\includegraphics[angle=-0,width=0.45\textwidth]{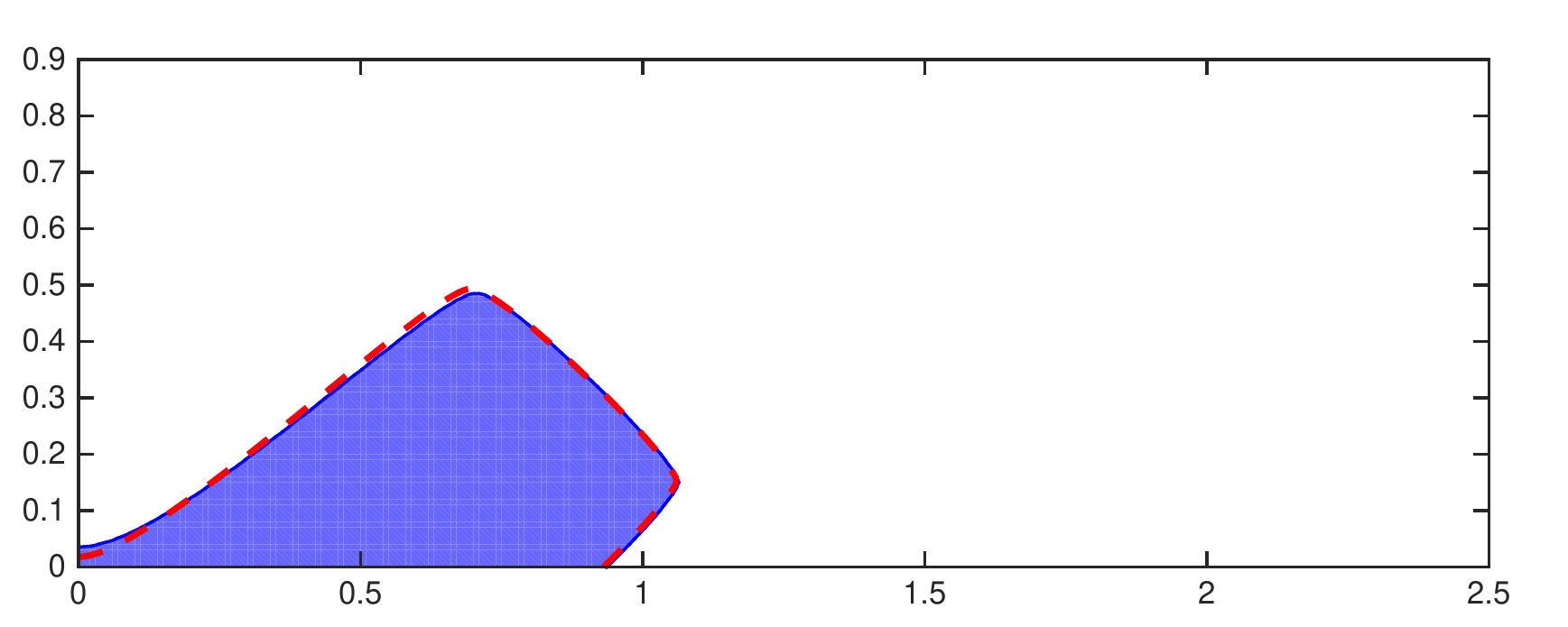}\;\;
\includegraphics[angle=-0,width=0.45\textwidth]{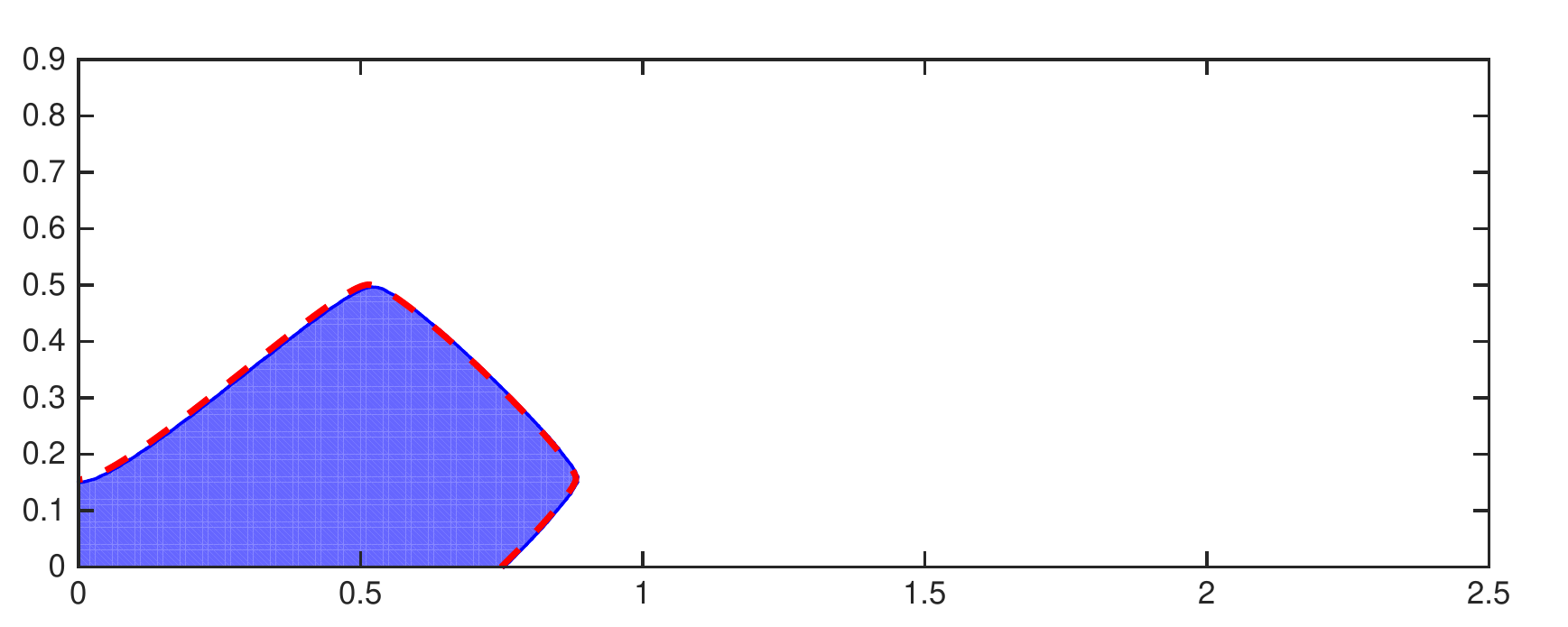}
\includegraphics[angle=-0,width=0.45\textwidth]{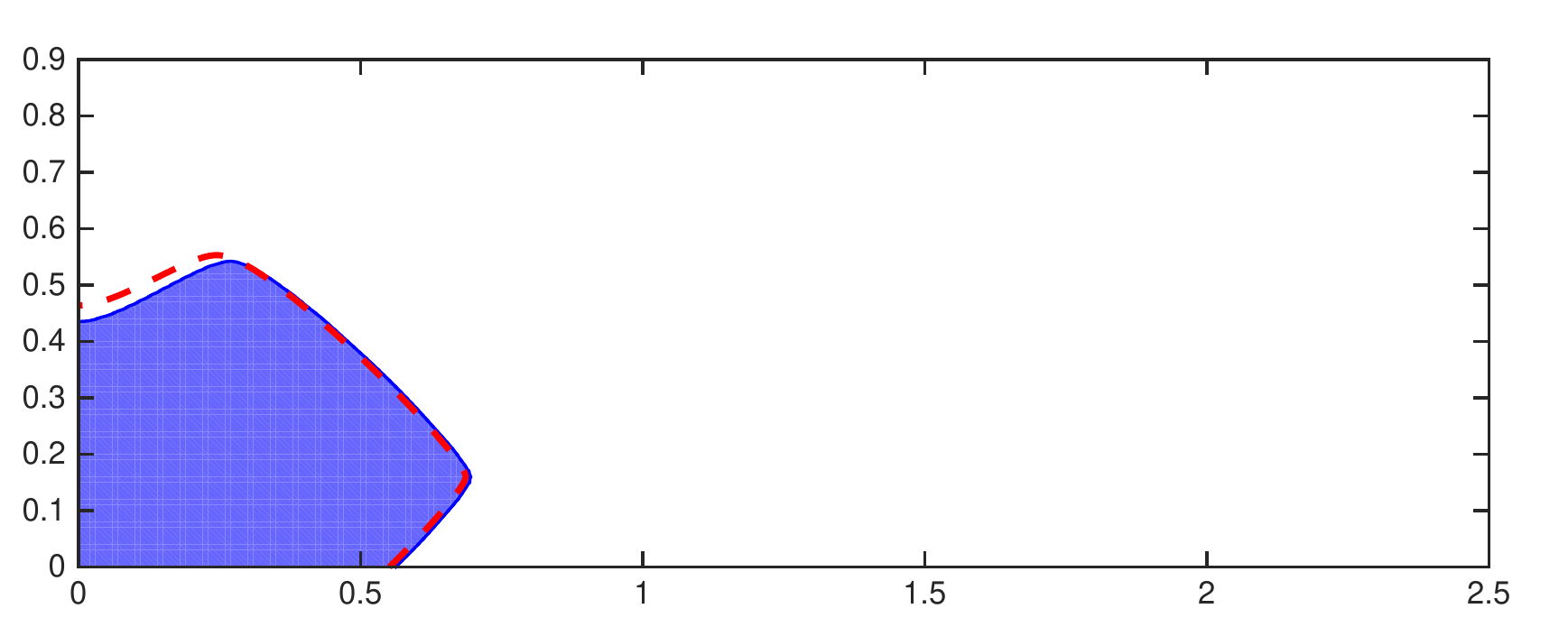}\;\;
\includegraphics[angle=-0,width=0.45\textwidth]{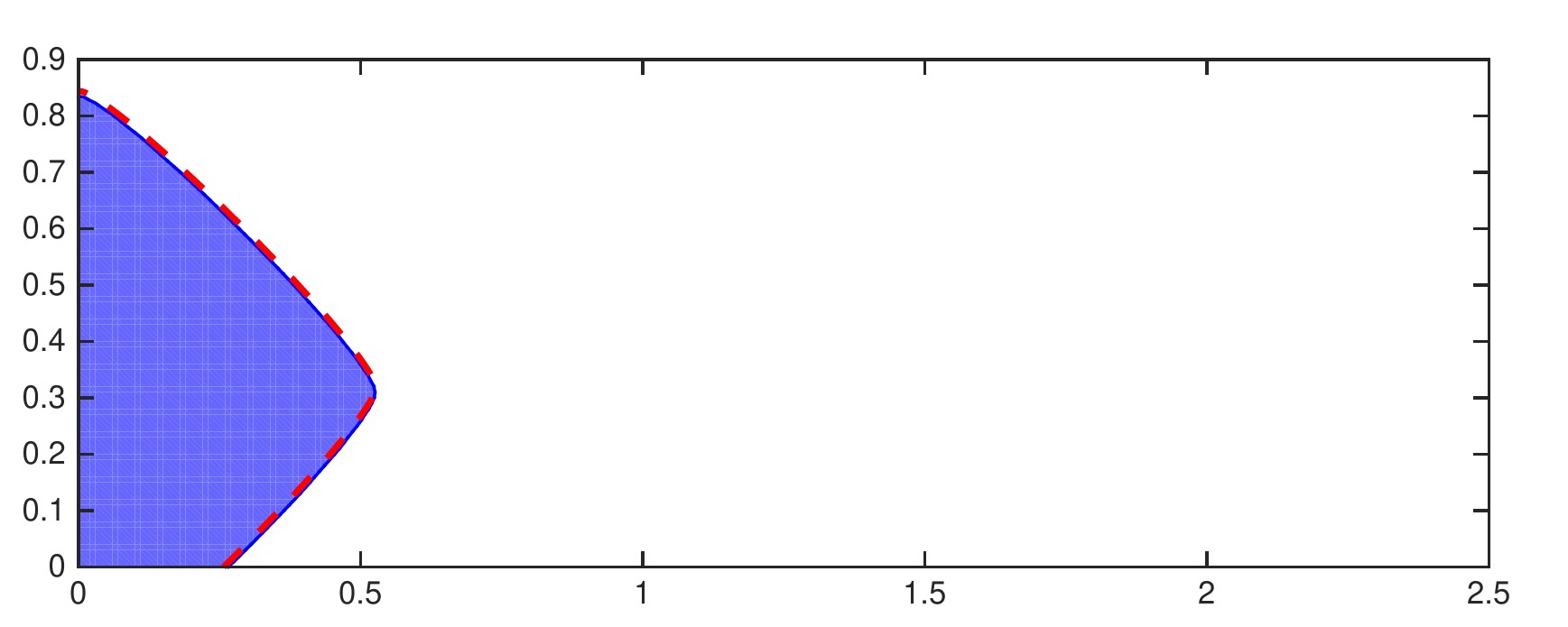}
\caption{[\;{\bf Example~(d)}\;] Interface profiles at times $t=0,0.01,0.02,0.03,0.04,0.1$ for the DI approximations with $\eps = 1/(64\pi)$, and the red dash line represents the SI approximations.} 
\label{fig:R14}
\end{figure}

The qualitative behaviour of the DI and SI approximations is compared in
Figs.~\ref{fig:R15}--\ref{fig:R14}. In all four examples we note an excellent
agreement between the two different approaches. This is particularly noteworthy in Example (a) with the occurrence of a topological change,  which is not covered by our asymptotic analysis.

\subsection{3d results}

In 3d, we compare our SI and DI approximations for the evolution of an initially
spherical island for the anisotropy $\gamma(\mR_x(\frac{\pi}{4})\mR_y(\frac{\pi}{4})\vec p)$, where $\gamma(\vec p)$ is given by \eqref{eq:asyform} with 
$d=3$, $\delta=0.1$, and where 
$\mR_x(\theta), \mR_y(\theta)$ are rotation matrices which rotate a vector through an angle $\theta$ within the $(y,z)$- and $(x,z)$-planes,
respectively. The initial interface is chosen to be a semisphere of radius
$0.4$, attached to the $(x,y)$-plane, and we let 
$\sigma = \cos(\frac{5\pi }{6})$.

\begin{figure}[!htp]
    \centering
    \includegraphics[width=0.45\textwidth]{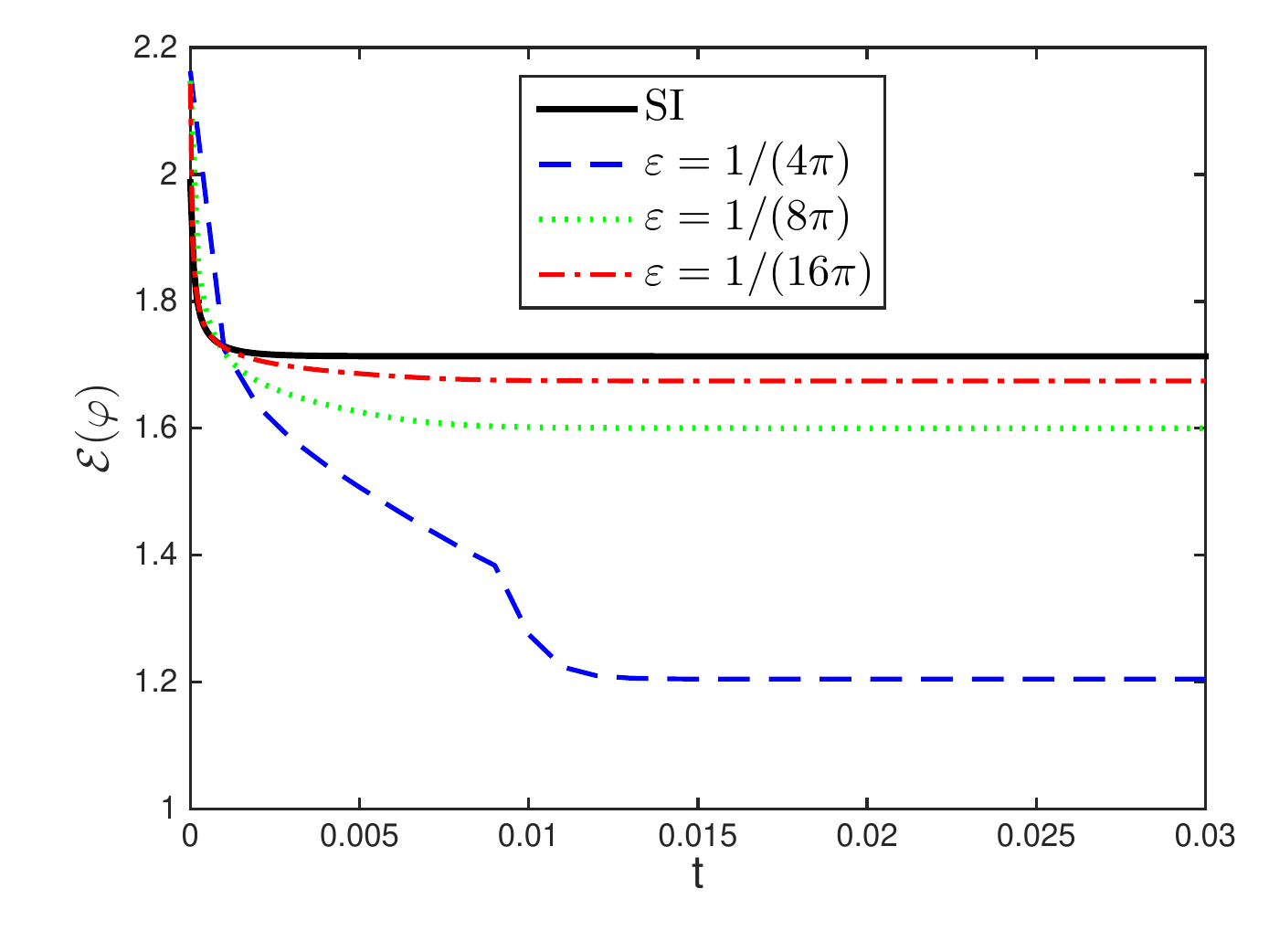}
    \includegraphics[width=0.45\textwidth]{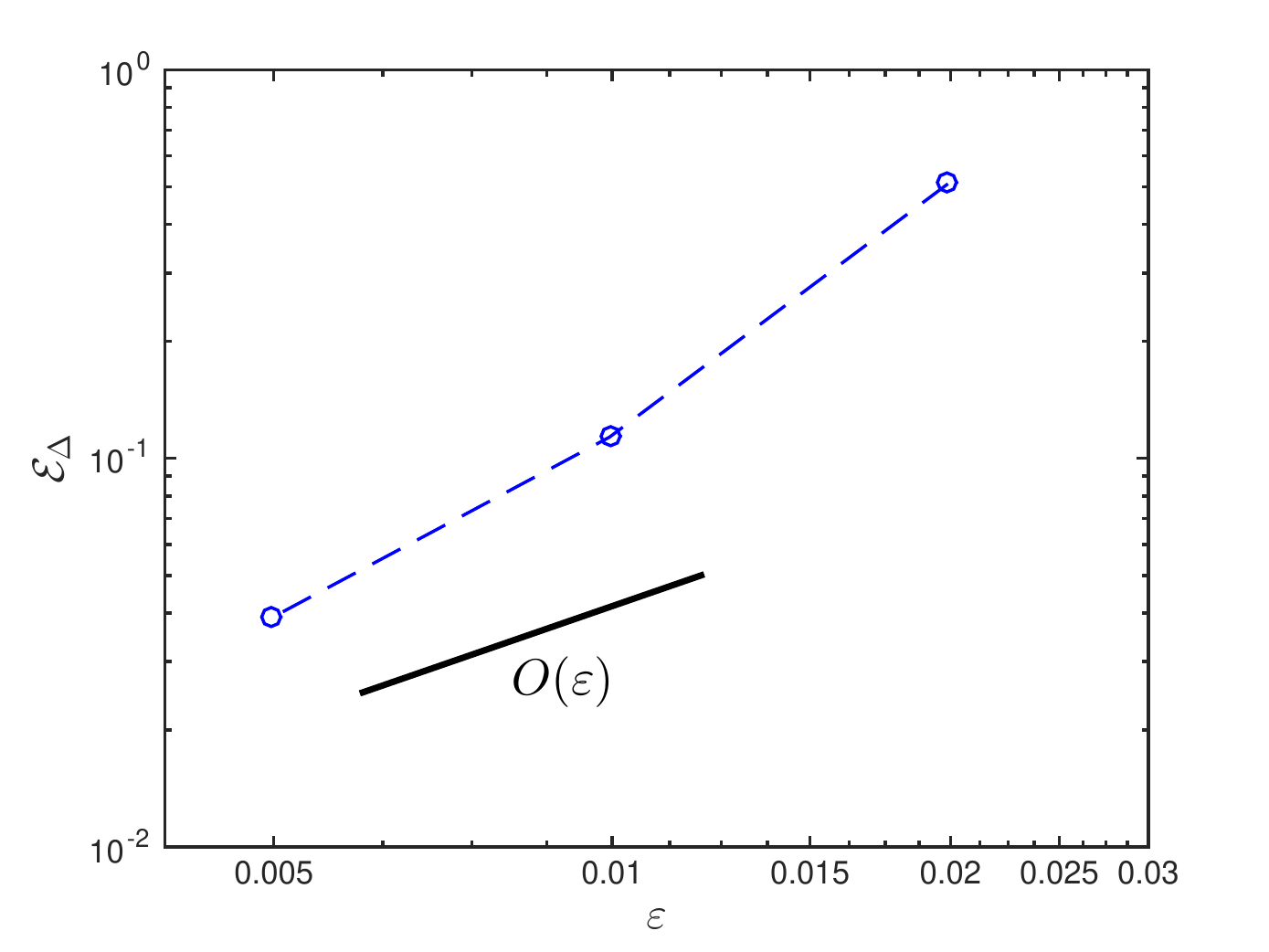}
    \caption{\revised{Left panel: The time history of the energy for the DI and SI approximations for the semisphere experiment in 3d. Right panel: The error $\mathcal{E}_\Delta$ of the energy at the final time between the DI and SI approximations plotted against $\eps$.}}
    \label{fig:3denergy}
\end{figure}

For the SI computation, we consider a polyhedral surface with 8256 triangles and
4225 vertices, and a time step size $10^{-4}$.
For our DI approximations, on the other hand, 
we consider the computational domain $\Omega = [-\frac12,\frac12]^3$ and 
as interfacial parameters consider
$\eps = 1 / (2^{2+i}\pi)$, for $i=0,\ldots,2$, 
with the corresponding discretization parameters
$N_f = 2^{5+i}$, $N_c = 2^{2+i}$, $\Delta t = 10^{-3} / 2^{2i}$.
In Fig.~\ref{fig:3denergy} we show the energy plots 
of the DI approximations and compare them with the corresponding SI 
simulation, noting once again an excellent agreement when $\eps$ is
sufficiently small. \revised{We also present a plot of the error in the energy between the DI and SI approximations against $\eps$. Note that the large error for $\eps = 1/(4\pi)$ is due to that DI
simulation wrongly predicting a pinch-off.}

\begin{figure}[!htp]
\center
\includegraphics[angle=-0,width=0.3\textwidth]{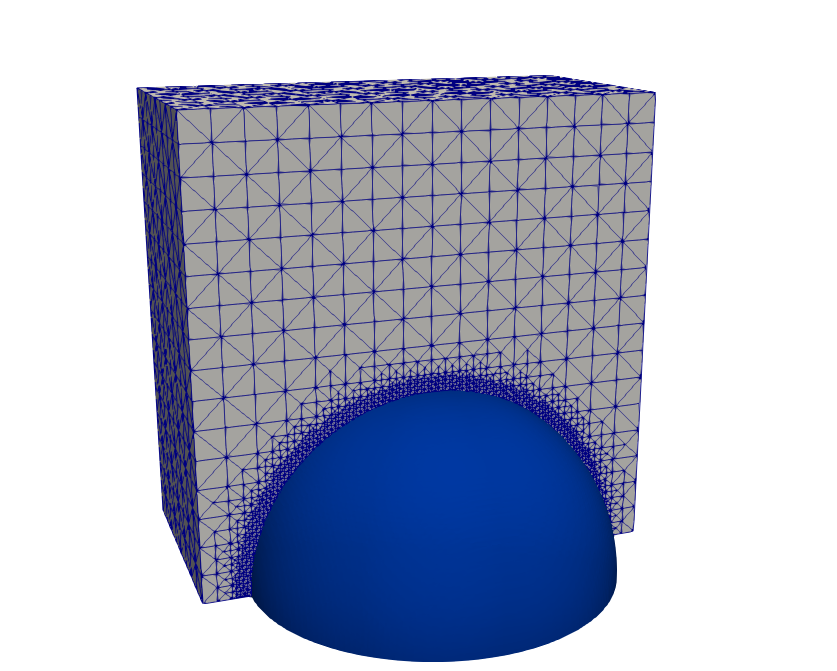}
\includegraphics[angle=-0,width=0.3\textwidth]{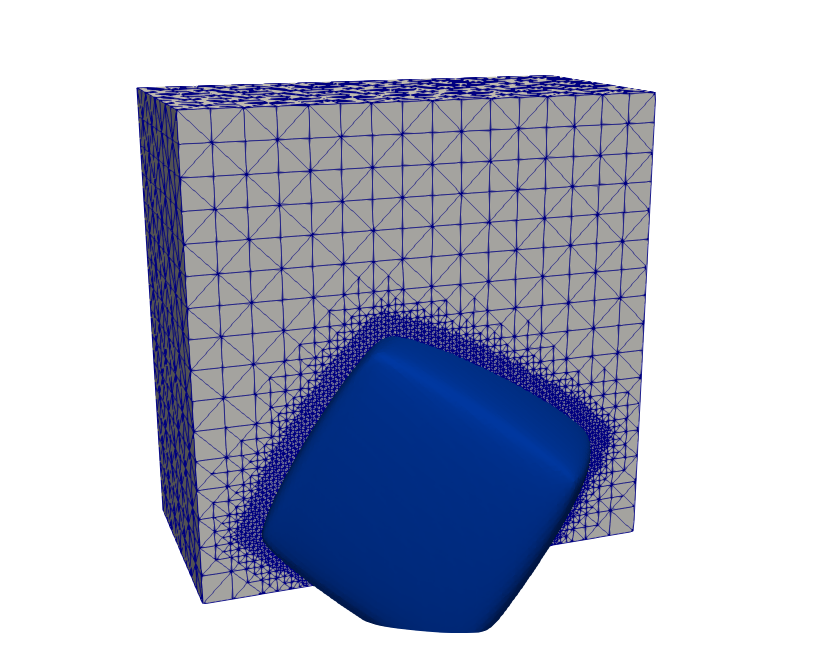}
\includegraphics[angle=-0,width=0.3\textwidth]{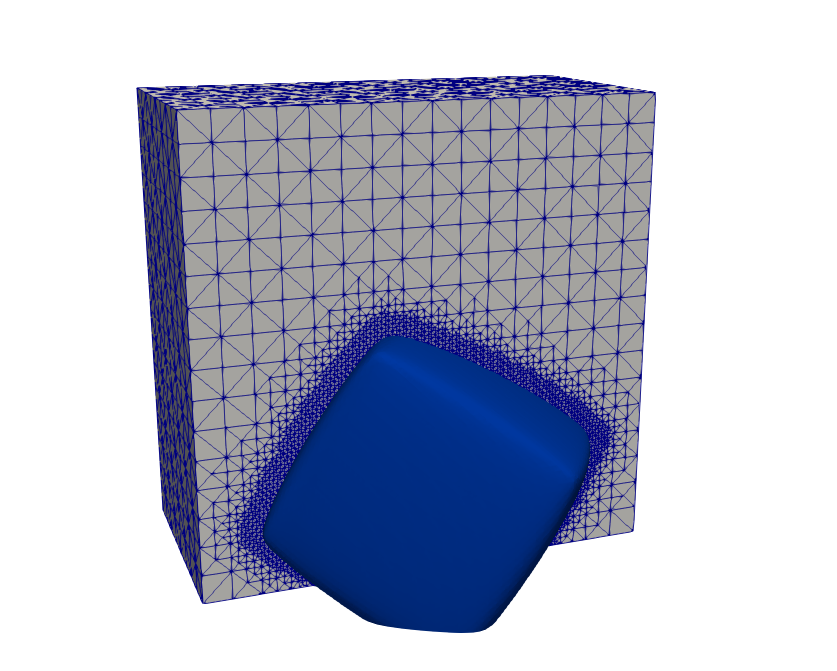} 
\includegraphics[angle=-0,width=0.45\textwidth]{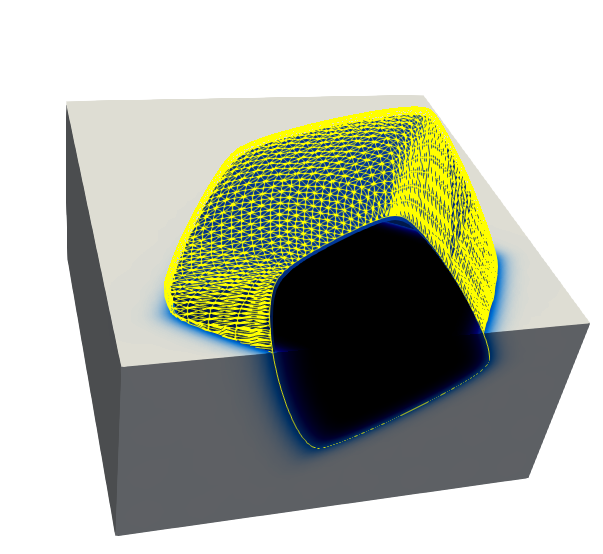}
\caption{A visualization of the zero level sets of the DI approximations
for $\eps=(16\pi)^{-1}$ at times $t=0,0.01,0.1$, together with a slice
through the adaptive mesh. Below a comparison \revised{between the DI and} the SI computation
at time $t=0.01$.
}
\label{fig:3disland2}
\end{figure}

Moreover, a qualitative comparison between the evolutions
of the interface for both approaches is shown in Fig.~\ref{fig:3disland2}.
In particular, at the bottom of Fig.~\ref{fig:3disland2} we see that the sharp
interface approximation agrees very well with the zero level set from the DI
computation, underlining once more our asymptotic 
analysis in Section~\ref{se:asymptotic}.


\section*{Acknowledgement}
We acknowledge the support from the RTG 2339 ``Interfaces, Complex Structures, and Singular Limits''
of the German Science Foundation (DFG) (Garcke, Knopf) and the Alexander von Humboldt Foundation (Zhao). 

\section*{Author contributions}
All authors contributed equally to the research presented in this article as well as to the preparation and revision of the manuscript.

\section*{Data availability}
The datasets generated during and/or analyzed during the current study are available from the corresponding author on reasonable request.

\section*{Conflict of interests}
The authors do not have any financial or non-financial interests that are directly or indirectly related to the work submitted for publication.

\end{document}